\documentclass[12pt]{article}
\usepackage{ifpdf}
\usepackage{subfigure}
\usepackage[english]{babel}
\usepackage{graphicx}
\usepackage{mathrsfs}
\usepackage{amsfonts}
\usepackage{amsmath}
\usepackage{moreverb}

\ifx\pdfoutput\undefined
\usepackage[hypertex]{hyperref}
\else
\usepackage[pdftex,hypertexnames=false]{hyperref}
\fi

\begin{document}
\thispagestyle{empty}
\vfill
\vspace{0.5cm}
\begin{center} {\Large A para-model agent for dynamical systems} 

\vspace{1.5cm}
{\it *** 10 years of model-free control methodology ***}
\end{center}

\vspace{0.5cm}
\begin{center}{\sc Lo\"ic MICHEL} \end{center}
\vspace{1.5cm}

\begin{quotation}              
\begin{center} {\bf Abstract} \end{center}
\noindent
Consider a dynamical system $ u \mapsto x, \dot{x} = f_{nl}(x,u)$ where $f_{nl}$ is a nonlinear (convex or nonconvex) function, or a combination of nonlinear functions that can eventually switch.
We present, in this preliminary work\footnote{This work is distributed under CC license \url{http://creativecommons.org/licenses/by-nc-sa/4.0/}. Email of the corresponding author : 
loic.michel54@gmail.com}, a generalization of the standard model-free
control, that can either {\it control} the dynamical system, given an output reference trajectory, or {\it optimize} the dynamical system as a derivative-free optimization based 
"extremum-seeking" procedure. Multiple applications\footnote{The control of the Epstein frame (described in \S 3.4) has been experimentally validated and the results have been presented at the 
French Symposium of 
Electrical Engineering in Grenoble, Jun. 2016 \url{http://sge2016.sciencesconf.org/}.} are presented and the robustness of the proposed method is studied in simulation.
\end{quotation}

\newpage

\tableofcontents

\newpage
\setcounter{page}{1}
\section{Introduction}
Based on the model-free control methodology, originally proposed by Fliess $\&$ Join \cite{fliess_original} \cite{esta} \cite{fliess2} ten years ago, which is referred to as a self-tuning controller in \cite{kumar} and which has been widely and successfully 
applied to many mechanical
and electrical processes, the para-model agent (PMA) aims to generalize the model-free control by not only controlling nonlinear system, but also performing an "extremum-seeking" control.
On the one hand, we studied the dynamic performances when controlling generic switched minimum phase, non-minimum phase systems (e.g. \cite{astrom} \cite{isidori} \cite{chen} 
\cite{Benosman} \cite{Gurumoorthy} \cite{Karagiannis} \cite{Barkana}) 
as well as the control of some nonlinear systems taken from applications in physics. On the other hand, we present how the PMA can be used to find the optimum of some classes of 
nonlinear functions.
The proposed para-model agent\footnote{A justification of the proposed name "para-model" is given in the note of the conclusion.} is a simple derivation of the discrete model-free control law \cite{Michel}. The last progresses result in two contributions: first, the substitution of the computation of 
the numerical derivatives in the original model-free control approach \cite{fliess2}, by an initialization function that makes 
the controller more robust when stabilizing for example switched processes\footnote{The first steps toward the elaboration of the proposed algorithm were to extend the capabilities of the model-free
control regarding the control of switched non-minimum phase systems \cite{michel_nm}.}. Then, we propose to extend the properties of the model-free controller to include the extremum seeking control 
of nonlinear systems without any computation of derivative or gradient. In this case, instead of tracking a working point of nonlinear systems, an appropriate choice of the output 
reference may stabilize nonlinear systems to their extremum.

The paper is structured as follows. Section 2 presents the general principle of the proposed para-model agent. In Section 3, some examples illustrate the control of generic 
switched linear systems, the control of a three-phase motor, the control of a ballistic fire, the control of a biological system and the control of the measure of magnetic hysteresis in the framework of magnetic materials characterization 
(this last application is also referred to as the control of nonlinear switched systems). 
Section 4 presents how the proposed PMA approach can be used as derivative-free optimization / "extremum-seeking" control.

\section{General Principle}

We consider a nonlinear SISO dynamical system to control:

\begin{equation}\label{eq:gen_sys}
u \mapsto y,  \quad \left\{ \begin{matrix}
\dot{x} = f_{nl}(x,u) \\
y = C x
\end{matrix} \right.
\end{equation}

\noindent
where $f_{nl}$ is a nonlinear system, the para-model agent is an application $(y^*, y) \mapsto u$ whose purpose is to control the output $y$ of (\ref{eq:gen_sys}) 
following an output reference $y^*$. {\it In simulation, the system (\ref{eq:gen_sys}) is controlled in its "original formulation" without any modification / linearization.}

\subsection{Definition of the closed-loop}

\begin{figure}[!h]
\centering
\includegraphics[width=11cm]{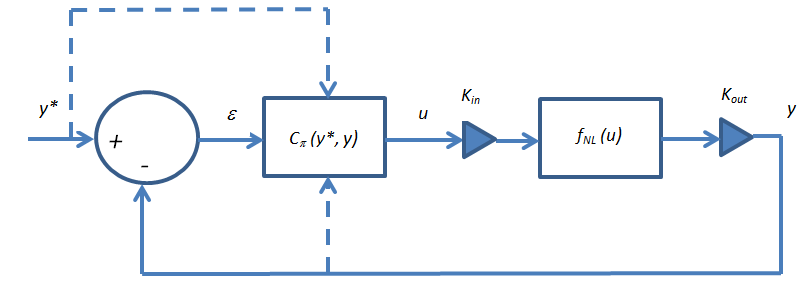}
\caption{Proposed PMA scheme to {\it control} or {\it optimize} a nonlinear system.}
\label{fig:CSM_gen}
\end{figure}

Consider the control scheme depicted in Fig. \ref{fig:CSM_gen} where $\mathcal{C}_{\pi}$ is the proposed PMA controller.

\subsection{Definition of the PMA algorithm}

For any discrete moment $t_k, \, k \in \mathbb{N}^*$, one defines the discrete controller $\mathcal{C}_{\pi}$ such that symbolically:
\begin{equation} \label{eq:iPI_discret_nm_eq}
\mathcal{C}_{\pi} : \begin{array}{c}
    \mathbb{R}^2 \rightarrow \mathbb{R}\\
    \displaystyle{  (y, y^*) \mapsto u_k =  \left. \int_0^t K_i \varepsilon_{k-1} d \, \tau \right|_{k-1} \underbrace{\left\{ u_{k-1}^i + {K_p} ( k_\alpha e^{-k_\beta k} - y_{k-1}) \right\}}_{u_{k}^i}}
     \end{array}
\end{equation}
where: $y^\ast$ is the output reference trajectory; $K_p$ and $K_I$ are real positive tuning gains; $\varepsilon_{k-1} = y^\ast_{k} - y_{k-1}$ is the tracking error; 
$u_k^i =  u_{k-1}^i + {K_p} ( k_\alpha e^{-k_\beta k} - y_{k-1})$ is an {\it internal} recursive term where $k_\alpha e^{-k_\beta k} - y_{k-1}$ is the associated {\it exponential tracking error} in which
$k_\alpha e^{-k_\beta k}$ is an {\it initialization function} where $k_\alpha$ and $k_\beta$ are real constants;
practically, the integral part is discretized using e.g. Riemann sums. The internal recursion\footnote{We refine the definition of the PMA algorithm, for which we aim to 
optimize the construction; in particular, further investigations concern the study of a direct 
recursion on $u_k$ taking into account the $K_i$-integration and thus comparing internal recursion 
(involving $u_k^i, u_{k-1}^i$) vs external recursion (involving $u_k, u_{k-1}^i$)...} 
on $u_{k}^i$ is defined such as: $u_{k}^i = u_{k-1}^i + {K_p} ( k_\alpha e^{-k_\beta k} - y_{k-1})$. 

\noindent
The set of $\mathcal{C}_{\pi}$-parameters of the controller, that needs to be adjusted by the user, is defined as the set of coefficients $\{K_p, K_i, k_\alpha, k_\beta \}$.

\paragraph{Practical algorithm}

A possible algorithmic implementation of the simple definition (\ref{eq:iPI_discret_nm_eq}) is given below (for all \texttt{ii} $\in \mathbb{N}^*$):

\vspace{1cm}

\begin{boxedverbatim}
y_int(ii) = k_alpha*exp(-k_beta*ii); % exp. init. function

para_exp_err = y_int(ii-1) - y(ii-1);    % exp. tracking error

para_stand_err(ii) = y_ref(ii) - y(ii-1); % stand. tracking error

para_u(ii) = para_u(ii-1) + Kp*para_exp_err; % internal recursion

para_G(ii) = Kint*para_stand_err(ii);  % def. of the integral part

para_tr(ii) = para_tr(ii-1) + h*(para_G(ii) + para_G(ii-1))/2; 
% trapezoidal integration

para_u_output = para_u(ii)*para_tr(ii); 
% final product (integrator X internal recursion)
\end{boxedverbatim}

\vspace{1cm}

\noindent
where:

\begin{itemize}
\item \texttt{ii} is the index of the sample in the (optional) vectorized process;
\item \texttt{exp} is the exponential function;
\item \texttt{para$\_$exp$\_$err} is the exponential tracking error;
\item \texttt{para$\_$stand$\_$err} is the (standard) output tracking error;
\item \texttt{para$\_$u} is the "internal" recursion;
\item \texttt{para$\_$G} constitutes the discrete integrator;
\item \texttt{para$\_$u$\_$output} is the output of the controller that corresponds to the final product between the discrete integrator and the internal recursive function. 
\end{itemize}

\noindent
and \texttt{k$\_$alpha}, \texttt{k$\_$beta}, \texttt{Kp} and \texttt{Kint} are real constants.

\paragraph{Remark} The proposed PMA algorithm could been seen as a "deformed" integrator since the internal recursion $u_k^i$ multiplies directly the integrator $\int_0^t K_i \varepsilon_{k-1} d \, \tau$.

\newpage
\subsection{Performances in simulation}\label{BFO_paragraph}

\subsubsection{Optimization of the closed-loop}

\begin{figure}[!h]
\centering
\includegraphics[width=12.5cm]{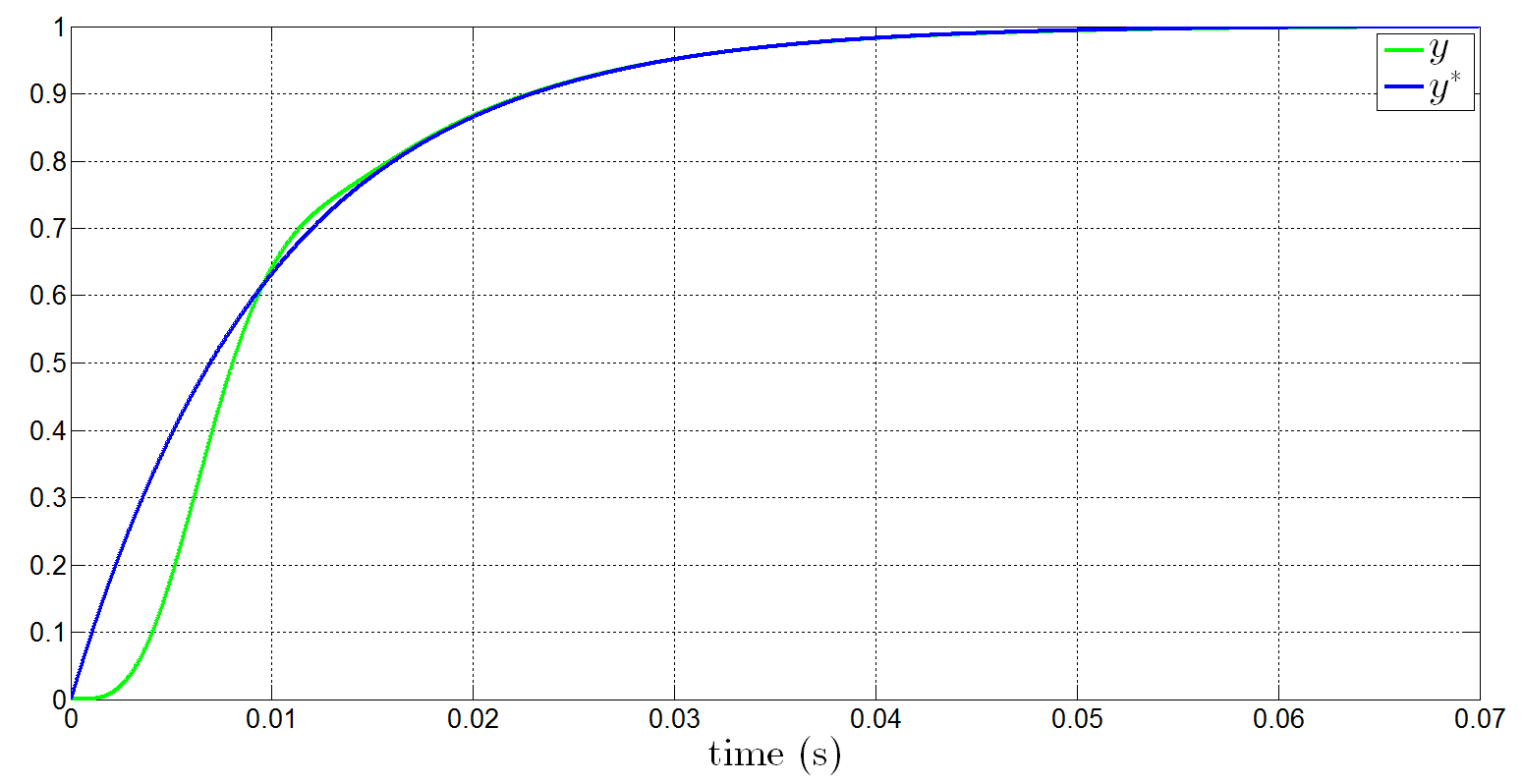}
\caption{Simulation with a set of $\mathcal{C}_{\pi}$-parameters arbitrary fixed to ensure at least asymptotic stability.}
\label{fig:BFO_non_opt_1}
\end{figure}

\begin{figure}[!h]
\centering
\includegraphics[width=12.5cm]{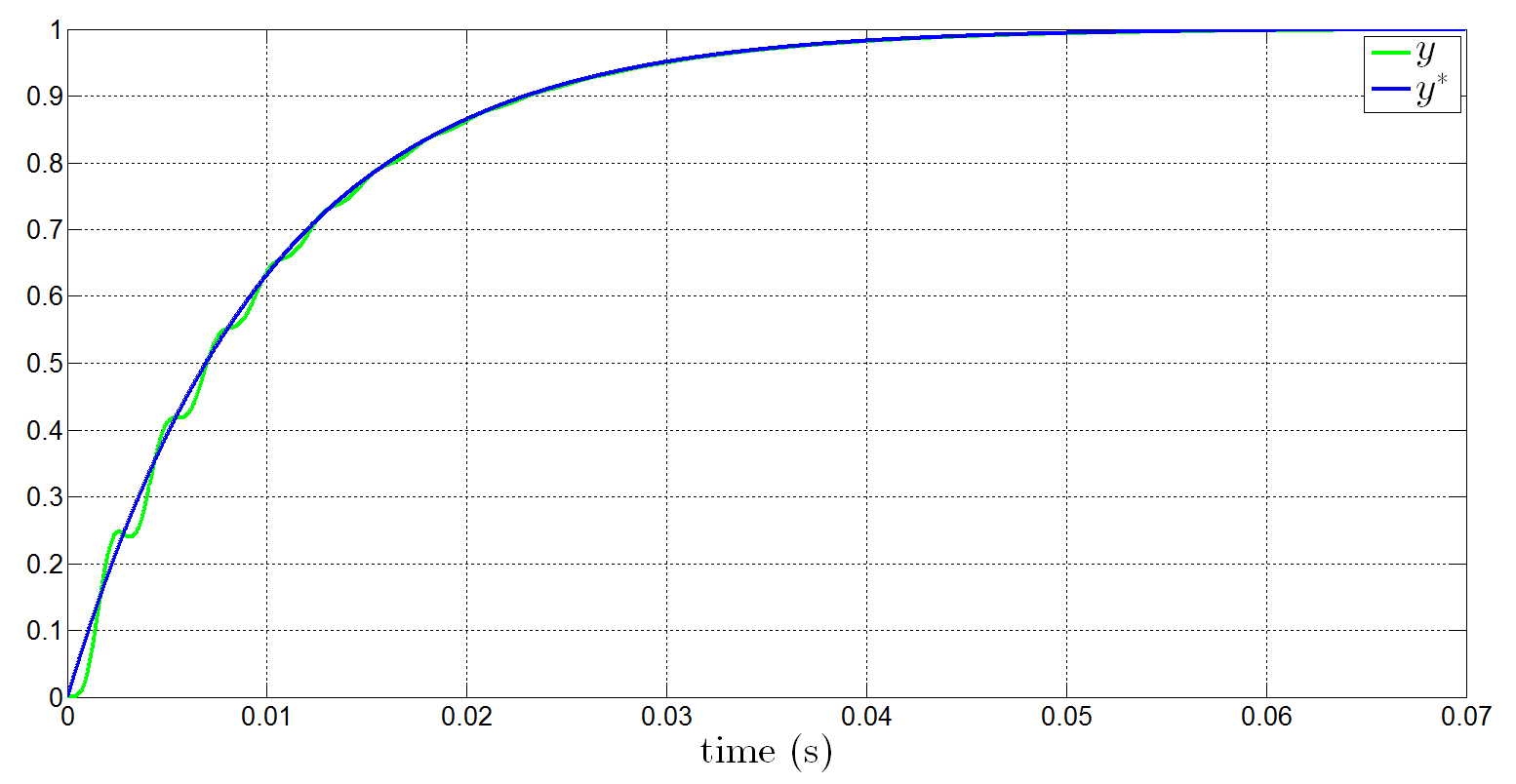}
\caption{Simulation with a set of $\mathcal{C}_{\pi}$-parameters BFO-optimized to minimize a transient performance index.}
\label{fig:BFO_opt_2}
\end{figure}

Optimizing the performances in simulation means that we want to solve the problem of finding the most appropriate set of $\mathcal{C}_{\pi}$-parameters relating to the minimization of some 
performances index \footnote{The classical performance index that are available are IAE, ISE, ITAE, and ITSE.
 see e.g. (\url{http://www.mathworks.com/matlabcentral/fileexchange/18674-learning-pid-tuning-iii--performance-index-optimization/content/html/optimalpidtuning.html})
 for a quick review (in the context of PID tuning).}, that may quantify the dynamical performances of the closed-loop. 
 This problem is thus equivalent to an optimization problem for which any optimization solver can be {\it a priori} used.
In particular, meta-heuristic solvers or derivative-free optimization solvers are preferred due to the pretty complexity of the closed-loop nonlinear form (in general). 
We are interested in using the "Brute Force Optimization" (BFO) solver \cite{BFO} that is very convenient and efficient to use. Figure \ref{fig:BFO_non_opt_1} illustrates 
a closed-loop first order system, whose the controlled transient has been BFO-optimized in Fig. \ref{fig:BFO_opt_2}.

\subsubsection{Sobol-based sensibility of the controlled transient}

To investigate the interactions between the $\mathcal{C}_{\pi}$-parameters that influence the minimization of the performances index, we propose a preliminary
study of the sensibility of the $\mathcal{C}_{\pi}$-parameters using the Sobol index methodology \cite{sobol_2}. Consider a controlled nonlinear system, for which the ISE index is evaluated under 
strict conditions w.r.t. the management of the $\mathcal{C}_{\pi}$-parameters\footnote{case 1 : the value of the evaluated index is bounded to 100 and a tolerance of $\pm 10 \%$ is permitted 
on the $\mathcal{C}_{\pi}$-parameters; case 2 : the value of the evaluated index is not bounded and a tolerance of $\pm 50 \%$ is permitted 
on the $\mathcal{C}_{\pi}$-parameters. The results are very similar between the two cases and the Sobol index for $\{K_p, K_i, k_\alpha, k_\beta \}$ are respectively $\{ 0.3324, 0.3318, 0.3326, 0.0002 \}$.}, 
the Sobol-based sensibility is evaluated only during the initialization / transient of the closed-loop.

\begin{figure}[!h]
\centering
\includegraphics[width=12.5cm]{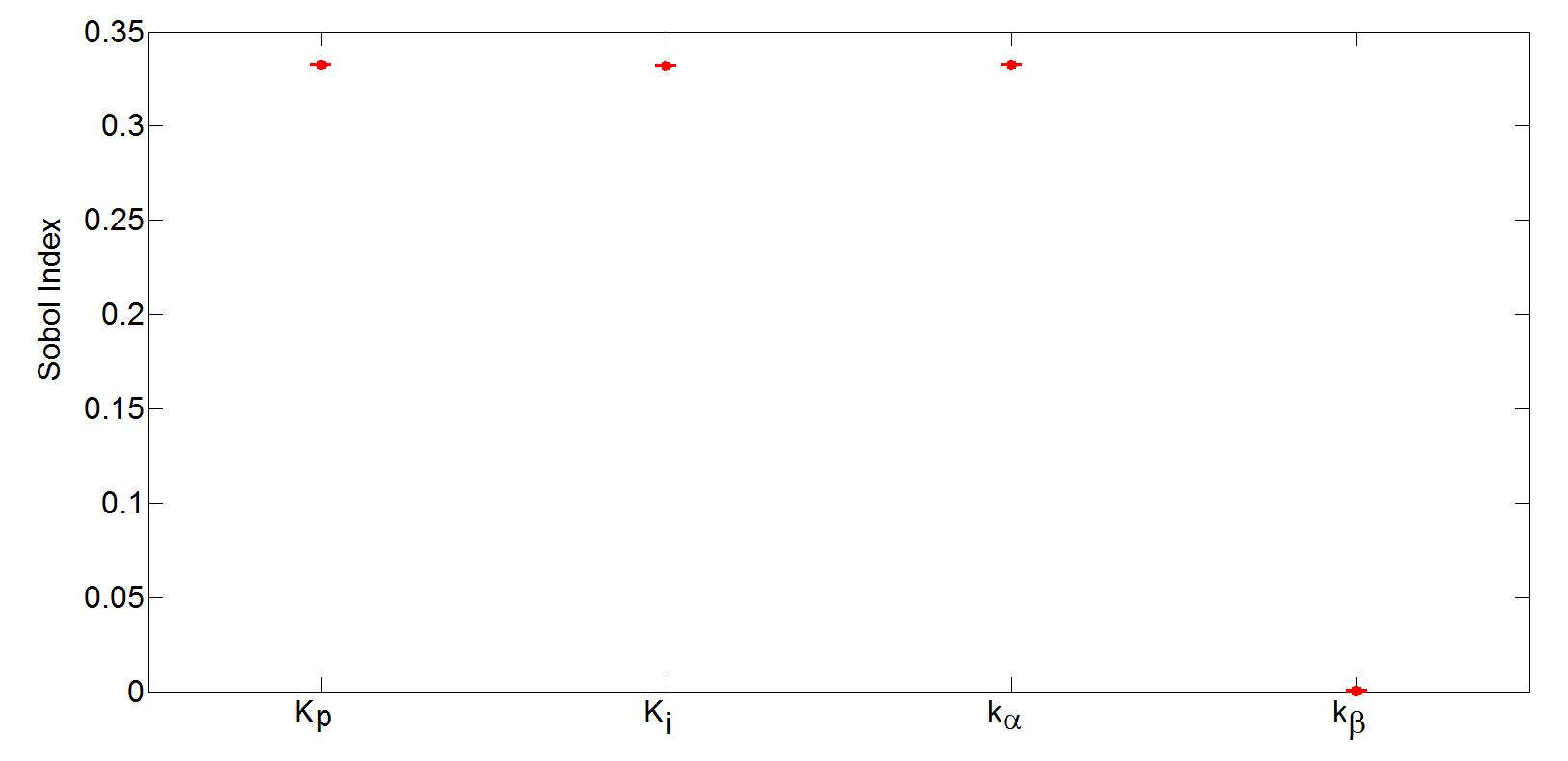}
\caption{Sobol-based analysis of the ISE index during the closed-loop initialization.}
\label{fig:fig_sobol}
\end{figure}

\noindent
Figure \ref{fig:fig_sobol} shows that {\it a priori} the coefficient $k_{\beta}$ does not influence the dynamical performances during the transient.

\clearpage
\section{Applications of the $\mathcal{C}_{\pi}$-control}

\subsection{Motor control in the $dq$ frame}

\begin{figure}[!h]
\centering
\includegraphics[width=12cm]{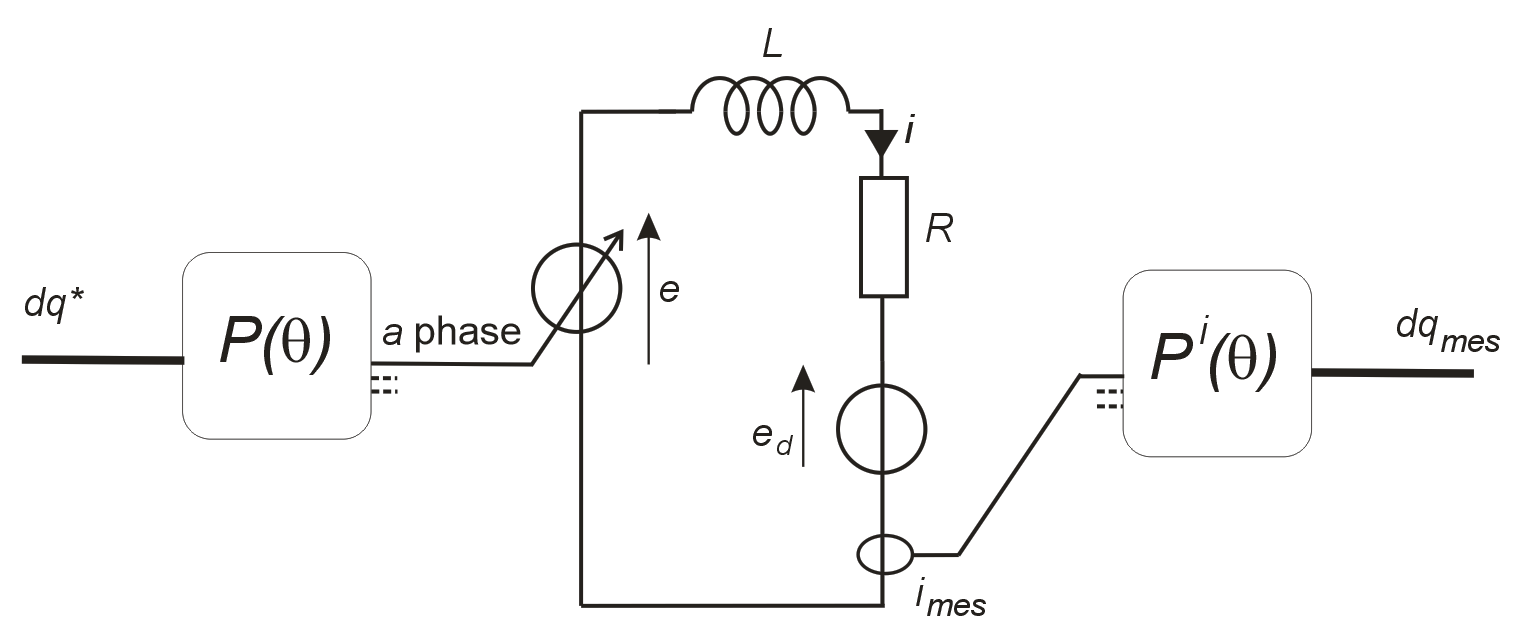}
\caption{Model of the motor in the $dq$ frame including an explicit disturbance $e_d$.}
\label{fig:motor_model_dq}
\end{figure}

Consider a three-phase motor driven in the $dq$ frame; the motor is supplied by a $e$ voltage source rotating at an $\omega$ angular frequency and modeled by a simple $RC$ circuit with 
an additional voltage source $e_d$ that acts as an (internal) disturbance. 
Figure \ref{fig:motor_model_dq} depicts the proposed model (a single phase is represented) where $P(\theta)$ and $P^i(\theta)$ are respectively the Park and the inverse Park transform. 
The purpose is to control simultaneously the $d$ and $q$ axis with an {\it a priori} unknown disturbance $e_d$ considering also that the angular frequency $\omega$ is increasing according to the time. 

\noindent
The disturbance $e_d$ is of the general form:

\begin{equation}
e_d := \left\{ \begin{array}{l}
 e_{d1} = A_1(t) \sin( \omega_d(t) t) \\
 e_{d2} = A_2(t) \sin \left( \omega_d(t) t - \frac{2\pi}{3} \right) \\
 e_{d3} = A_3(t) \sin \left(\omega_d(t) t + \frac{2\pi}{3 } \right)  \\
\end{array} \right.
\end{equation}

\noindent
where the amplitude $A_i$ and the angular frequency $\omega_d$ of each phase $i$ could depend on the time. The control structure is composed of two $\mathcal{C}_{\pi}$ controllers:
the $d$ axis is "maintained" close to zero ($d^*$ denotes the output ref. and $d_{mes}$, the controlled output) and the $q$ axis tracks a specific reference 
($q^*$ denotes the output ref. and $q_{mes}$, the controlled output).

\begin{figure}[!h]
\centering
\includegraphics[width=16.1cm]{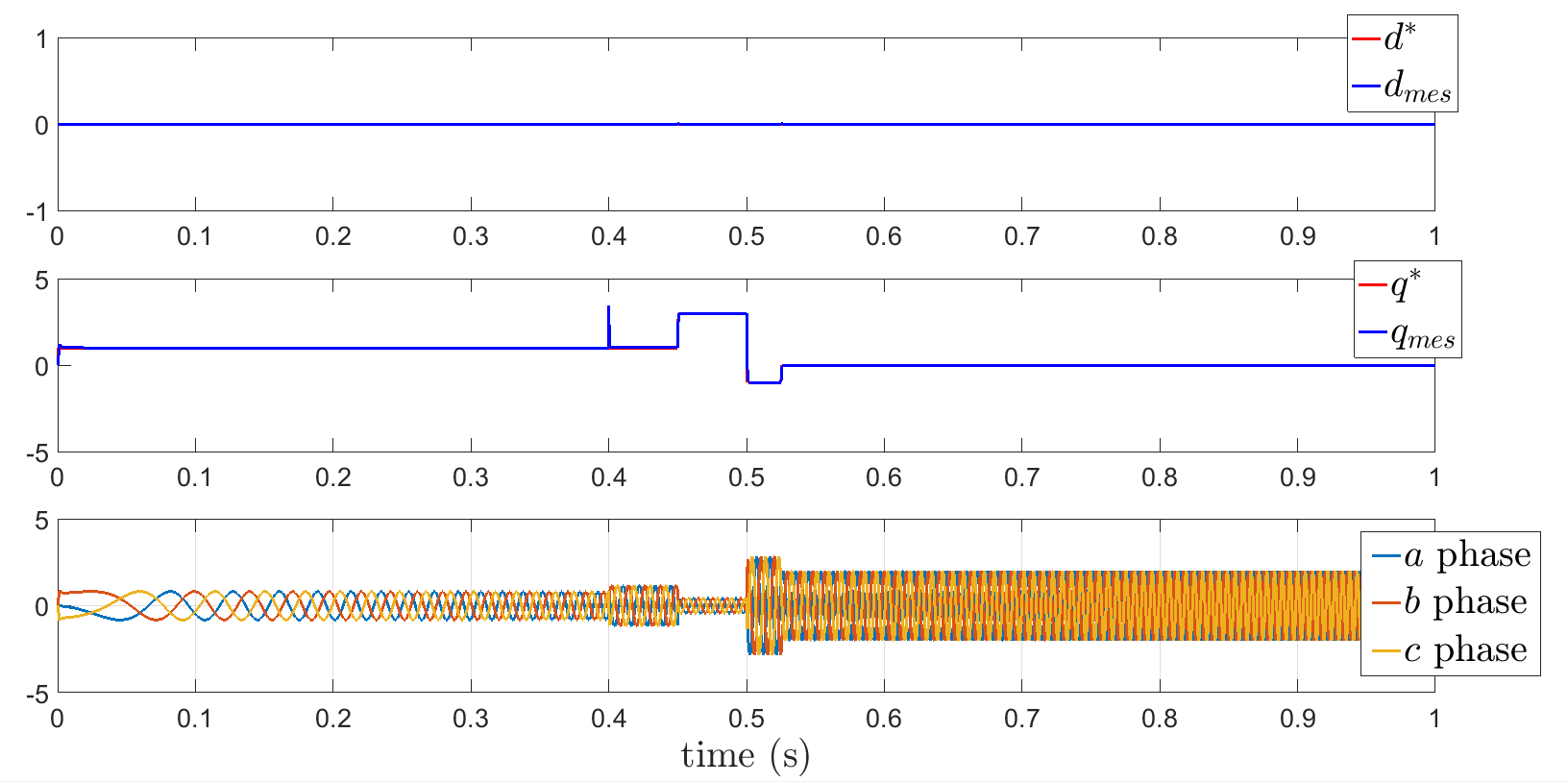}
\caption{Control in the $dq$ frame with "simple" disturbance $e_d$.}
\label{fig:motor_model_dq_fig1}
\end{figure}

\begin{figure}[!b]
\centering
\includegraphics[width=16.1cm]{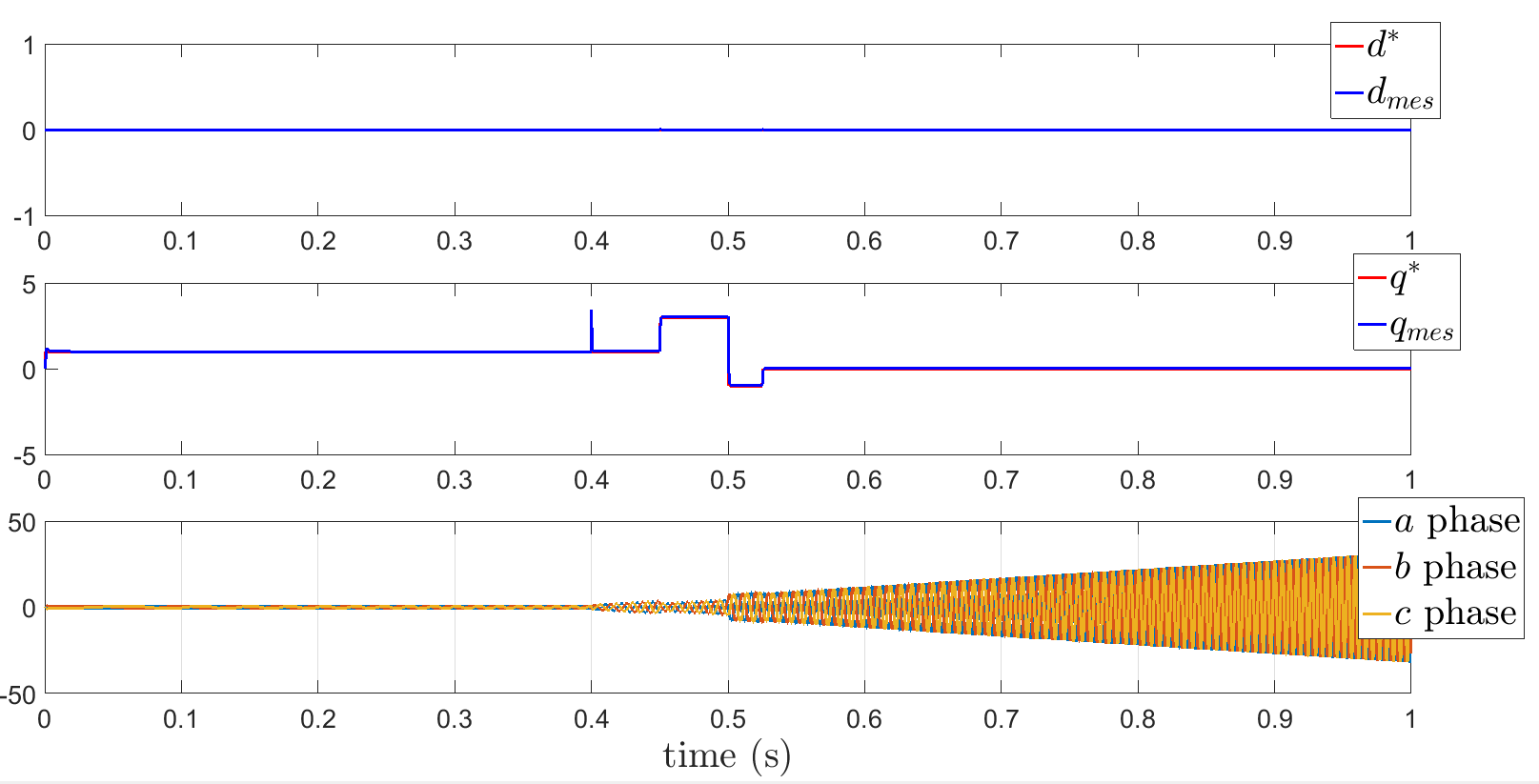}
\caption{Control in the $dq$ frame with an increasing amplitude of each component of $e_d$.}
\label{fig:motor_model_dq_fig2}
\end{figure}

\begin{figure}[!h]
\centering
\includegraphics[width=16cm]{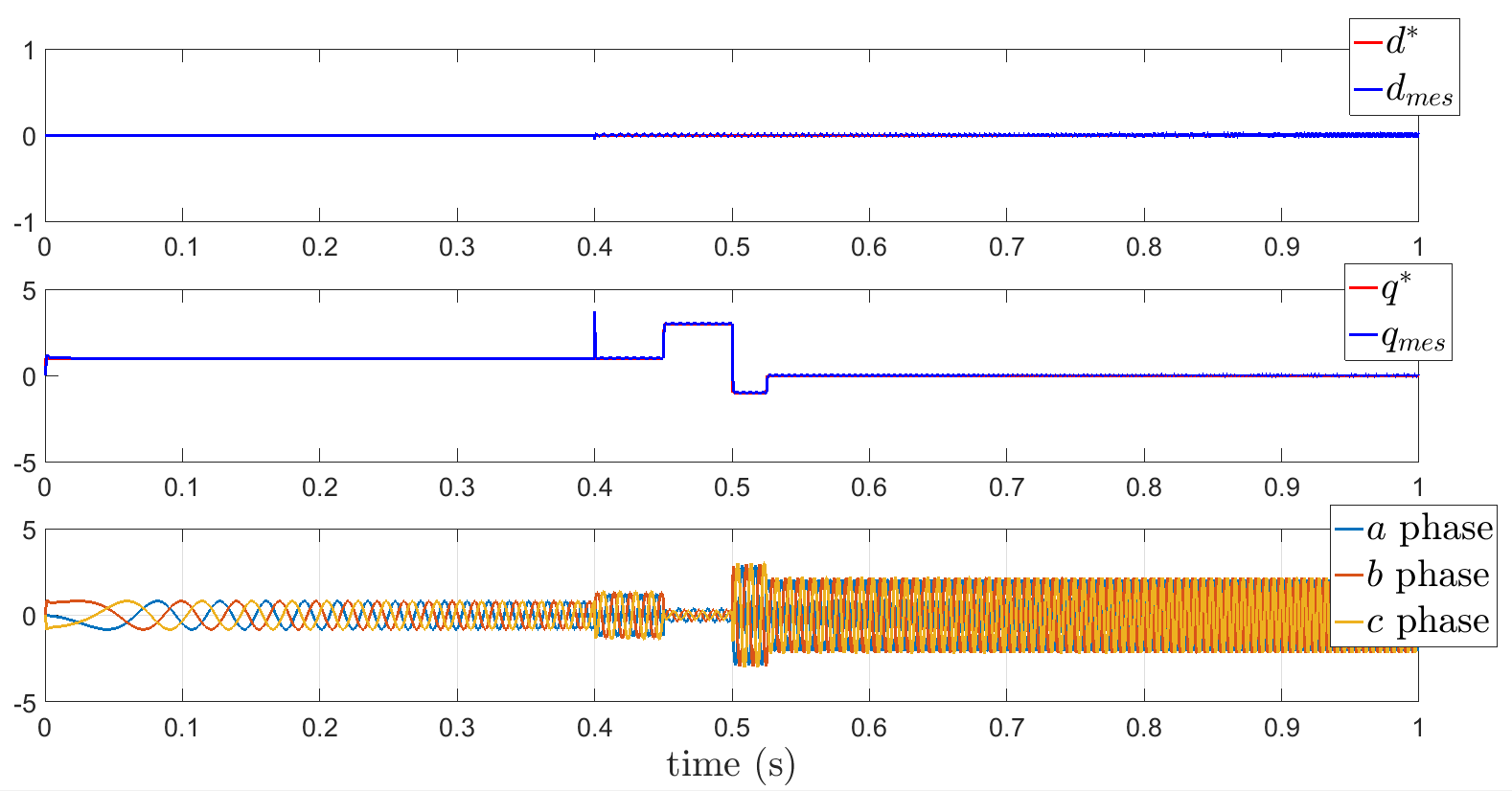}
\caption{Control in the $dq$ frame with variation of amplitude of each component of $e_d$.}
\label{fig:motor_model_dq_fig3}
\end{figure}

\begin{figure}[!b]
\centering
\includegraphics[width=16cm]{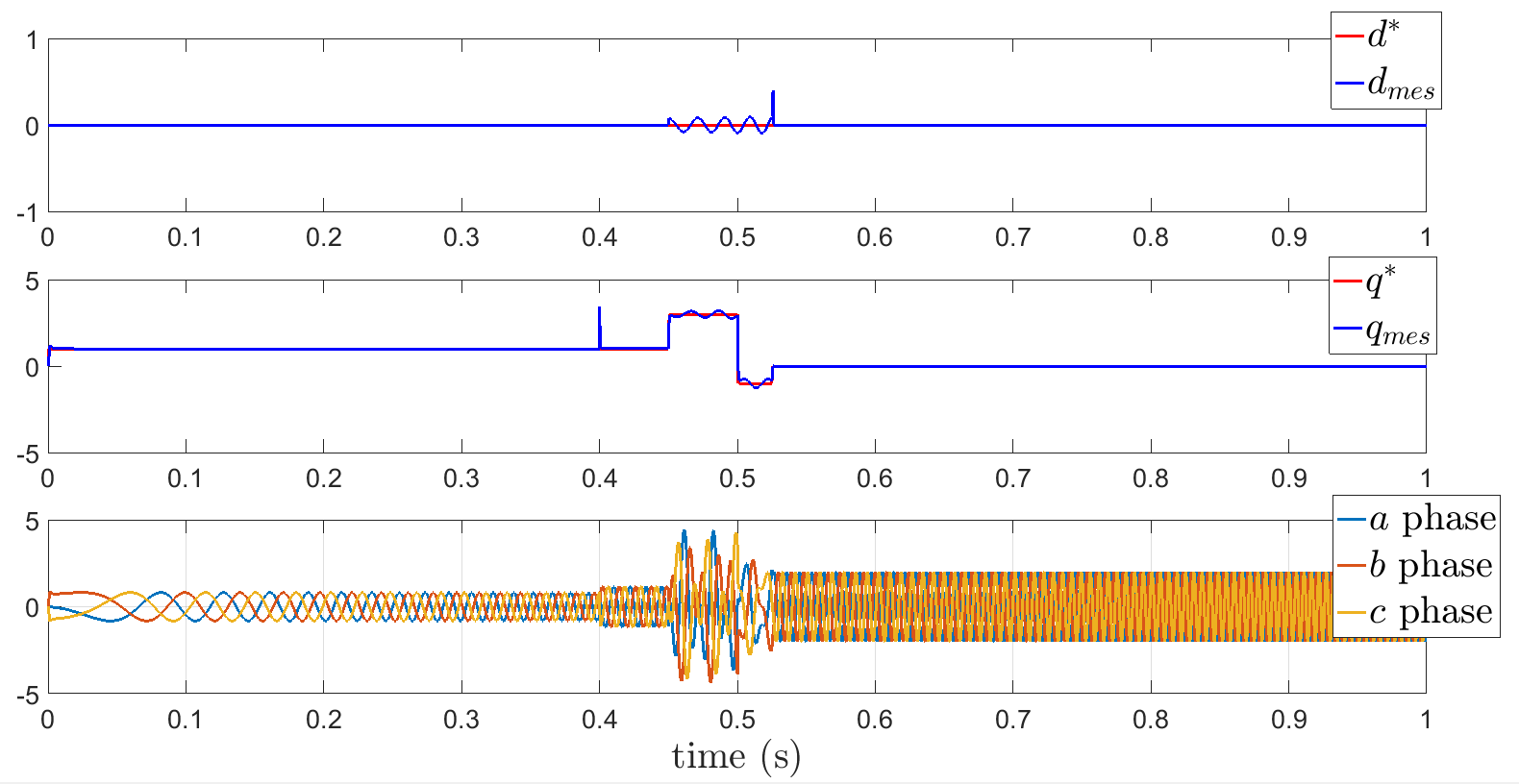}
\caption{Control in the $dq$ frame with variation of the $\omega_d$ frequency of $e_d$.}
\label{fig:motor_model_dq_fig4}
\end{figure}

The following figures illustrate some cases with different "behavior" of the disturbance $e_d$. Figure \ref{fig:motor_model_dq_fig1} presents the most simple case where 
$A_1 = A_2 = A_3 = \mathrm{Cst}$ and $\omega_d = \mathrm{Cst}$; 
in Fig. \ref{fig:motor_model_dq_fig2} is considered a time-varying disturbance where
$A_1 = A_2 = A_3$ are increasing according to the time; in Fig. \ref{fig:motor_model_dq_fig3}, small variations of amplitudes of $A_1, A_2$ and $A_3$ are considered 
(symbolically, $A_1 \sim A_2 \sim A_3$), and finally, Fig. \ref{fig:motor_model_dq_fig4} depicts the case where $\omega_d$ is time-varying only over a short period of time.



\clearpage
\subsection{Control of switched non-minimum and minimum phase systems}\label{para_sw_sys}

\begin{figure}[!h]
\centering
\includegraphics[width=13.5cm]{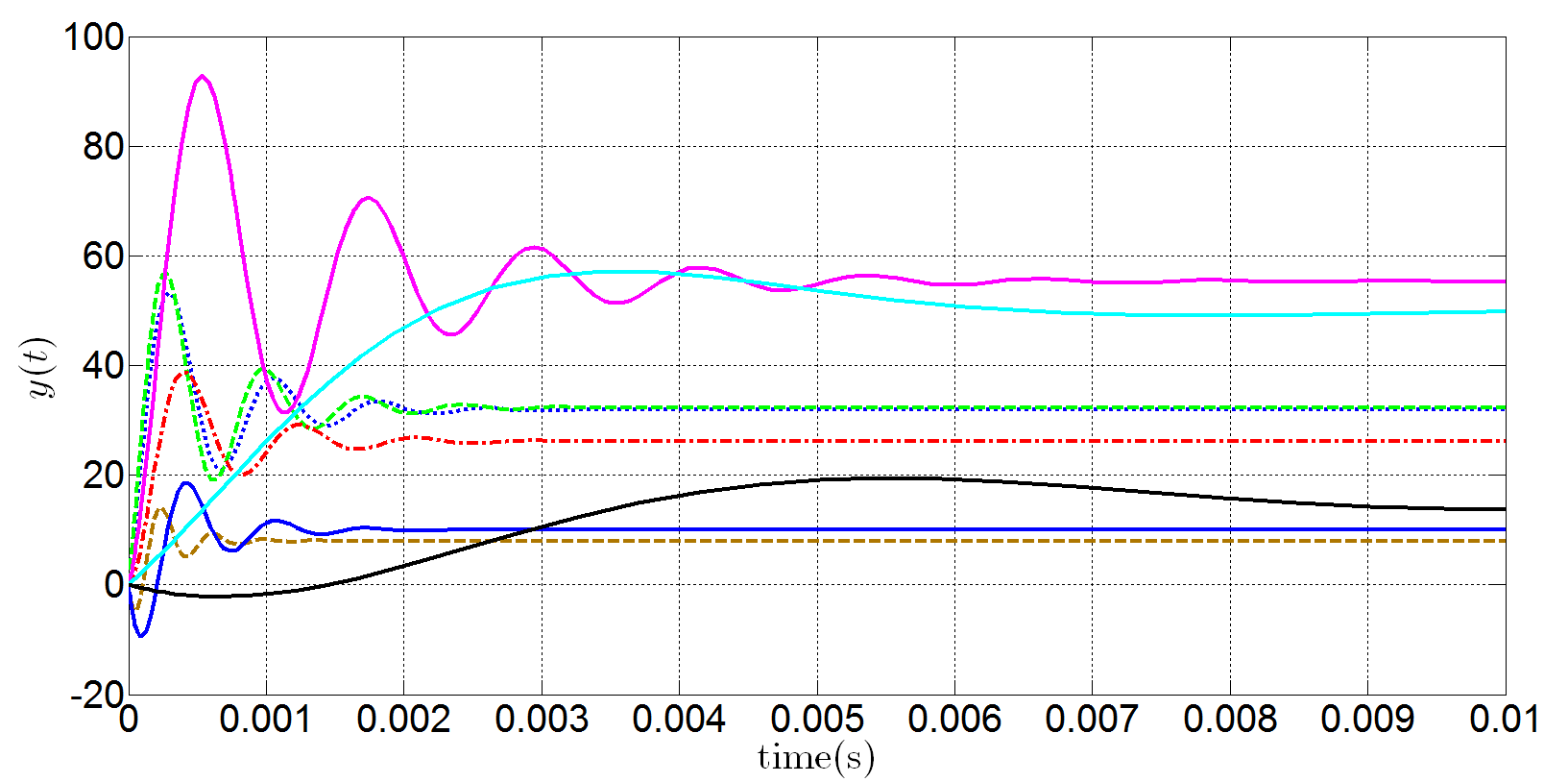}
\caption{Step responses of each system $\Sigma_i$.}
\label{fig:sys_nom}
\end{figure}

Consider the set $\Sigma$ of stable linear systems such that  $\Sigma = \{\Sigma_i \}, i = 1,2, \cdots n$, which are minimum and non-minimum phase systems, 
and which are considered as {\it unknown} in the sense that no explicit model has been identified for control purposes.
\noindent
Assume now that for all systems, there exists an integer $p = \{1, \cdots , 8\}$, called the switching index, such that
during a short time window, we have:

\begin{equation}\label{eq:sys_linear}
\Sigma_p( u \longmapsto y) := \left\{ \begin{array}{l}
\dot{x}(t) = A_p x(t)  + B_p u(t) \\
 y = C_p x(t)
      \end{array} \right.
\end{equation}

\noindent
where $u$ and $y$ are respectively the input and the output of the system $\Sigma_p$ ($p$ is the switching index). The step responses of these $p$ systems are presented Fig. \ref{fig:sys_nom}.

\begin{figure*}[ht] 
\centering
\subfigure[\footnotesize Step responses of each system $\Sigma_i$.]{\includegraphics[width=16cm]{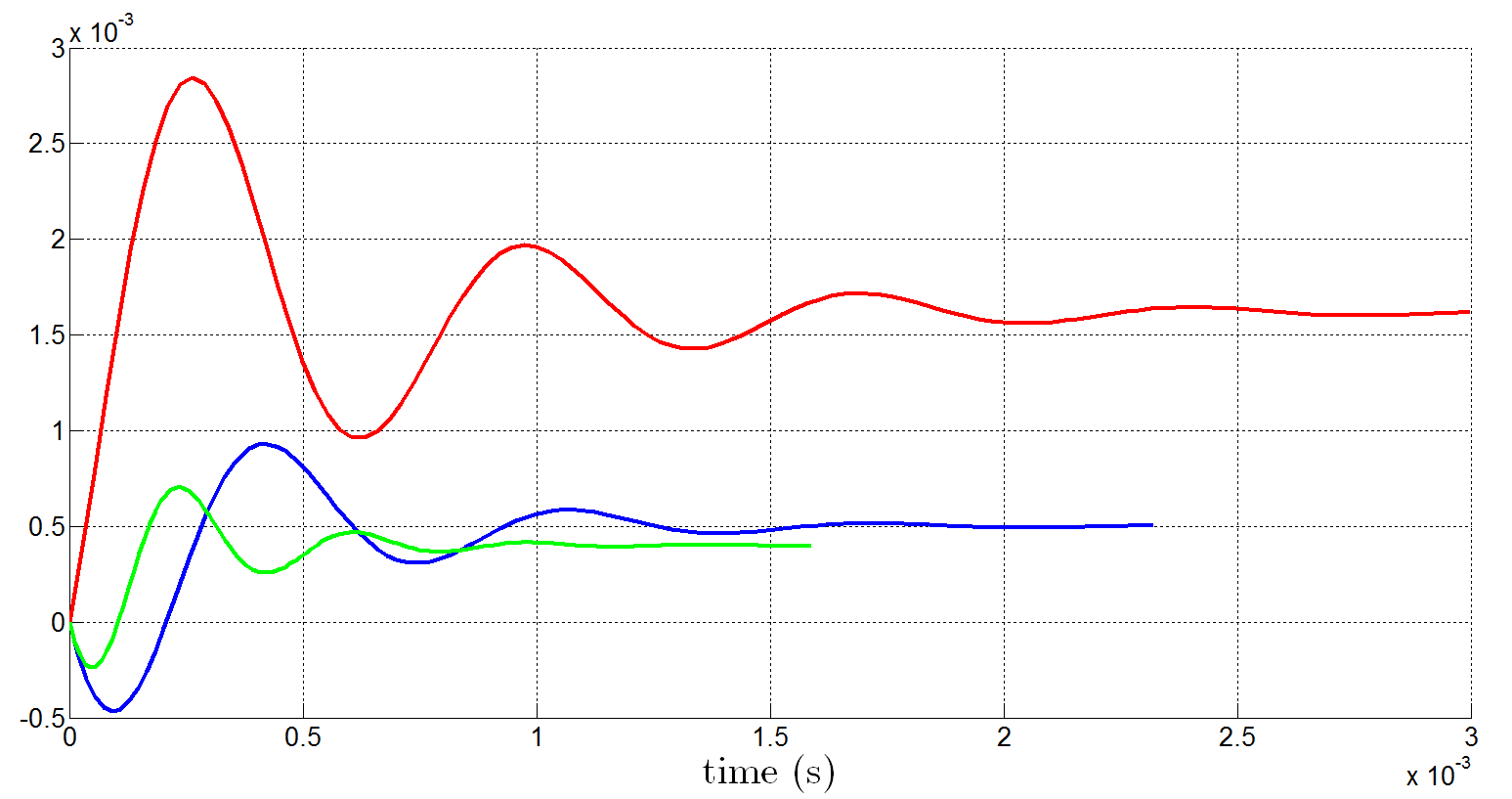} \label{fig:first}}
\subfigure[\footnotesize Controlled switched sequence.]{\includegraphics[width=16.2cm]{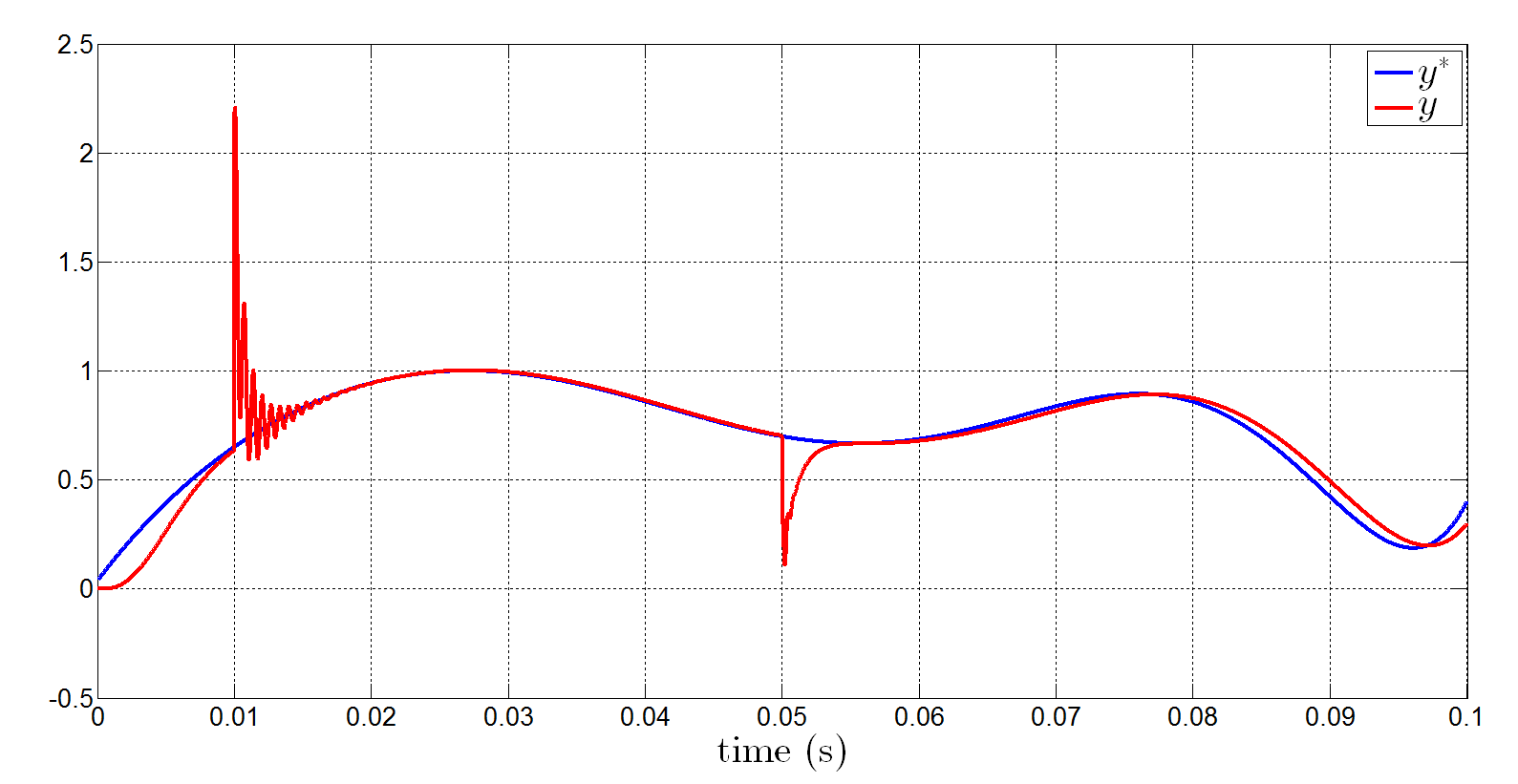} \label{fig:sec}}
\caption{Switched sequence $\# 1$ for $t_{1} = 0.01$ s and $t_2 = 0.05$ s.}
\label{fig:Fig_1_seq}
\end{figure*}

\begin{figure*}[ht] 
\centering
\subfigure[\footnotesize Step responses of each system $\Sigma_i$.]{\includegraphics[width=16cm]{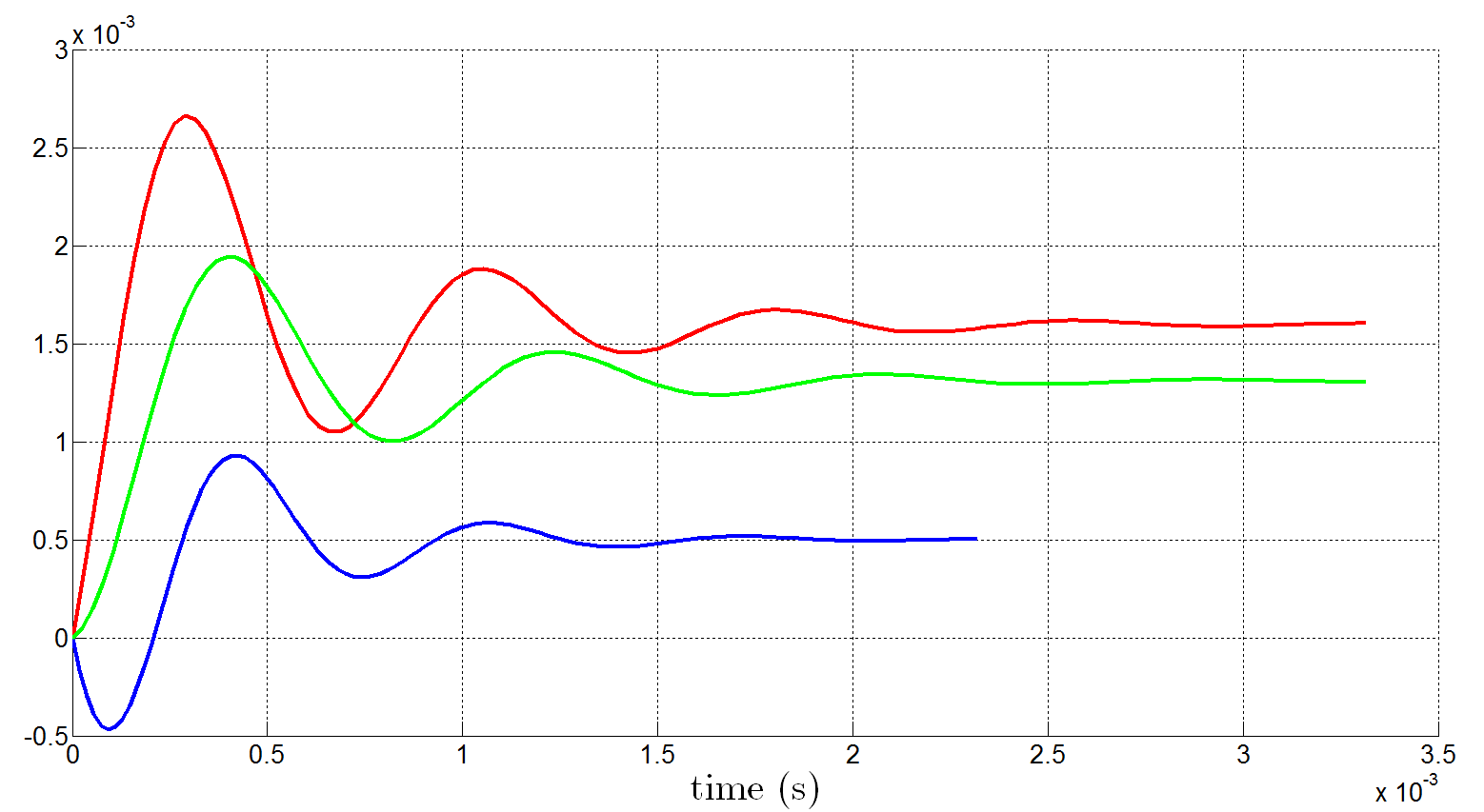} \label{fig:first}}
\subfigure[\footnotesize Controlled switched sequence.]{\includegraphics[width=16.2cm]{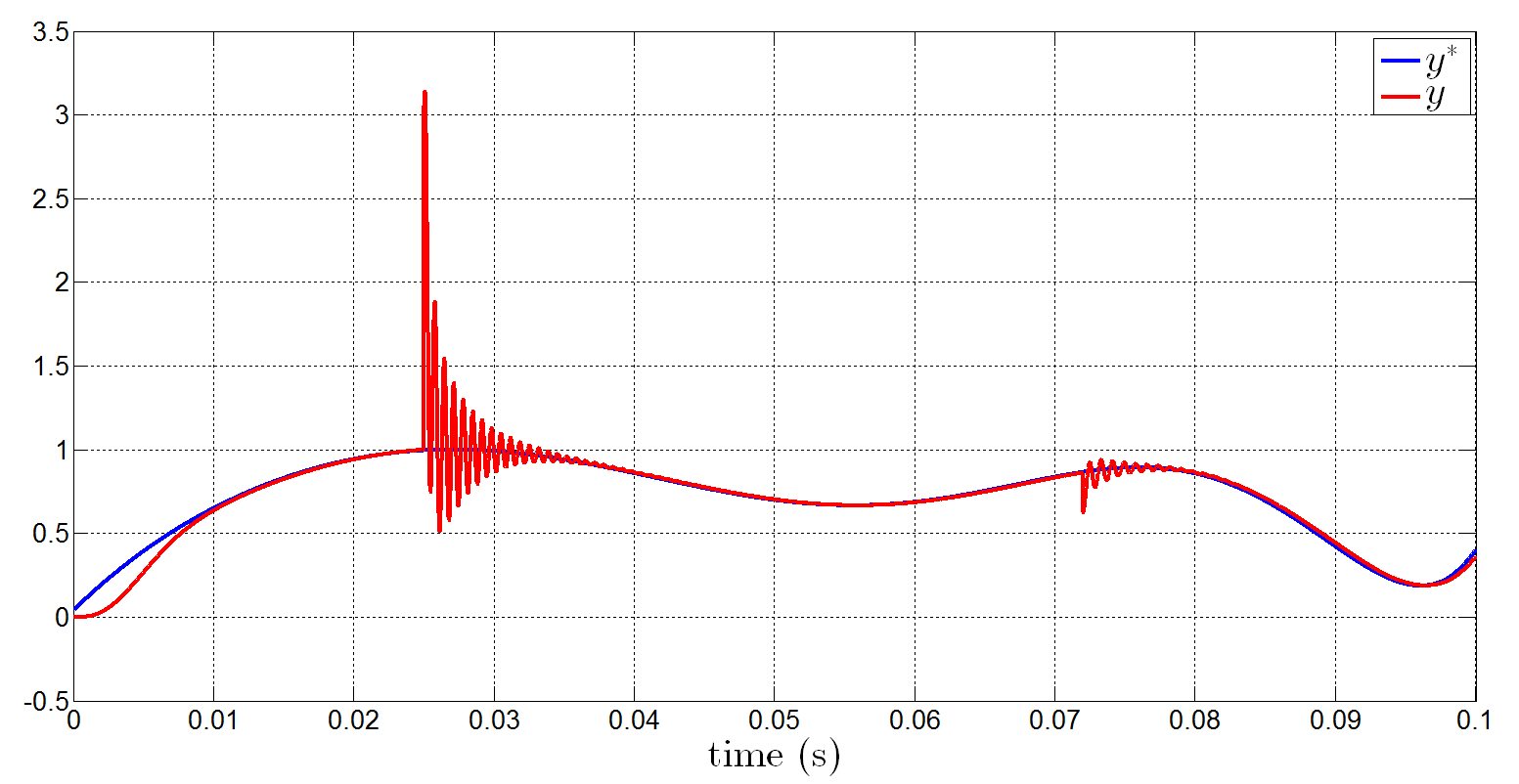} \label{fig:sec}}
\caption{Switched sequence $\# 2$ for $t_{1} = 0.025$ s and $t_2 = 0.072$ s.}
\label{fig:Fig_2_seq}
\end{figure*}

\begin{figure*}[!b] 
\centering
\subfigure[\footnotesize Step responses of each system $\Sigma_i$.]{\includegraphics[width=16cm]{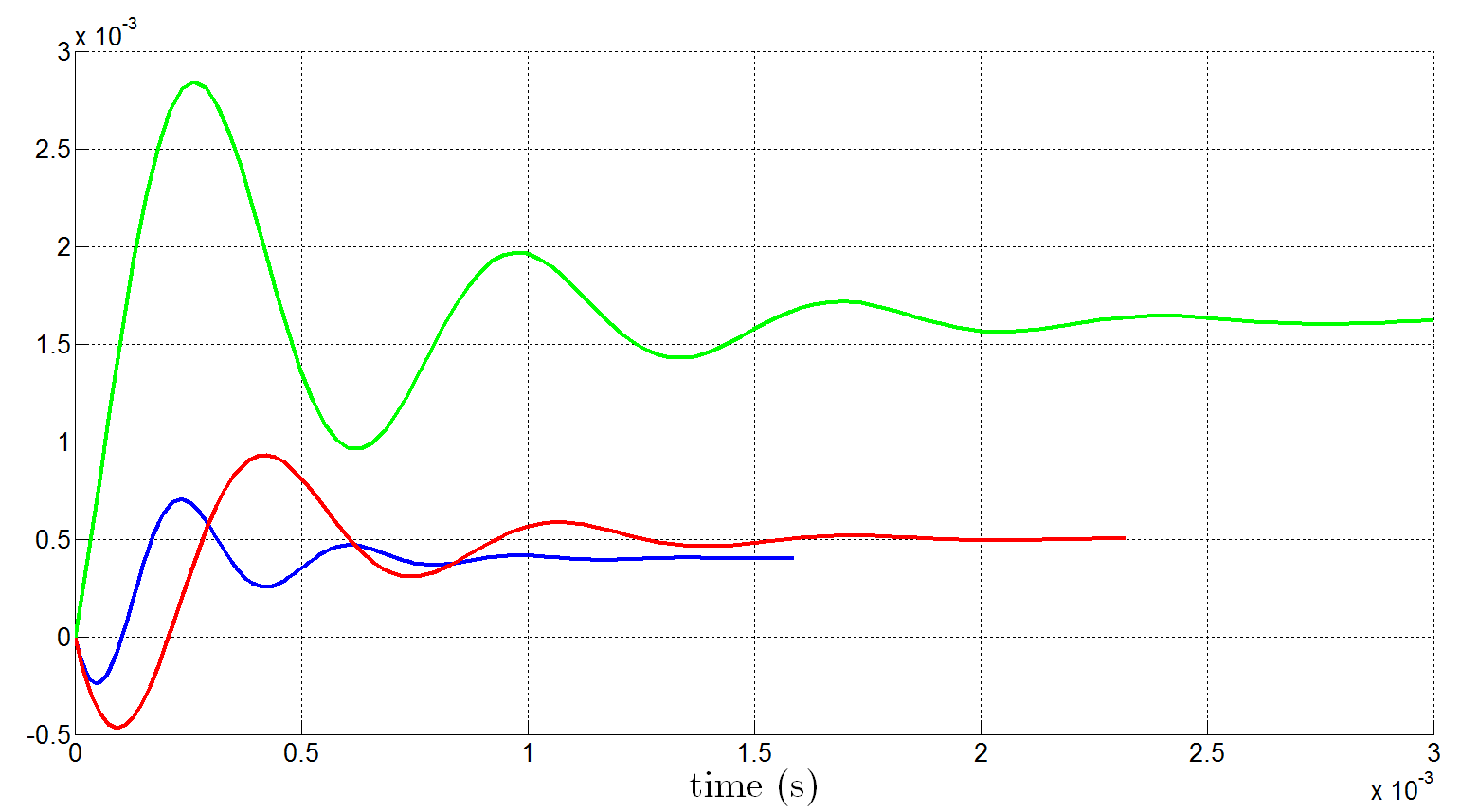} \label{fig:first}}
\subfigure[\footnotesize Controlled switched sequence.]{\includegraphics[width=16.2cm]{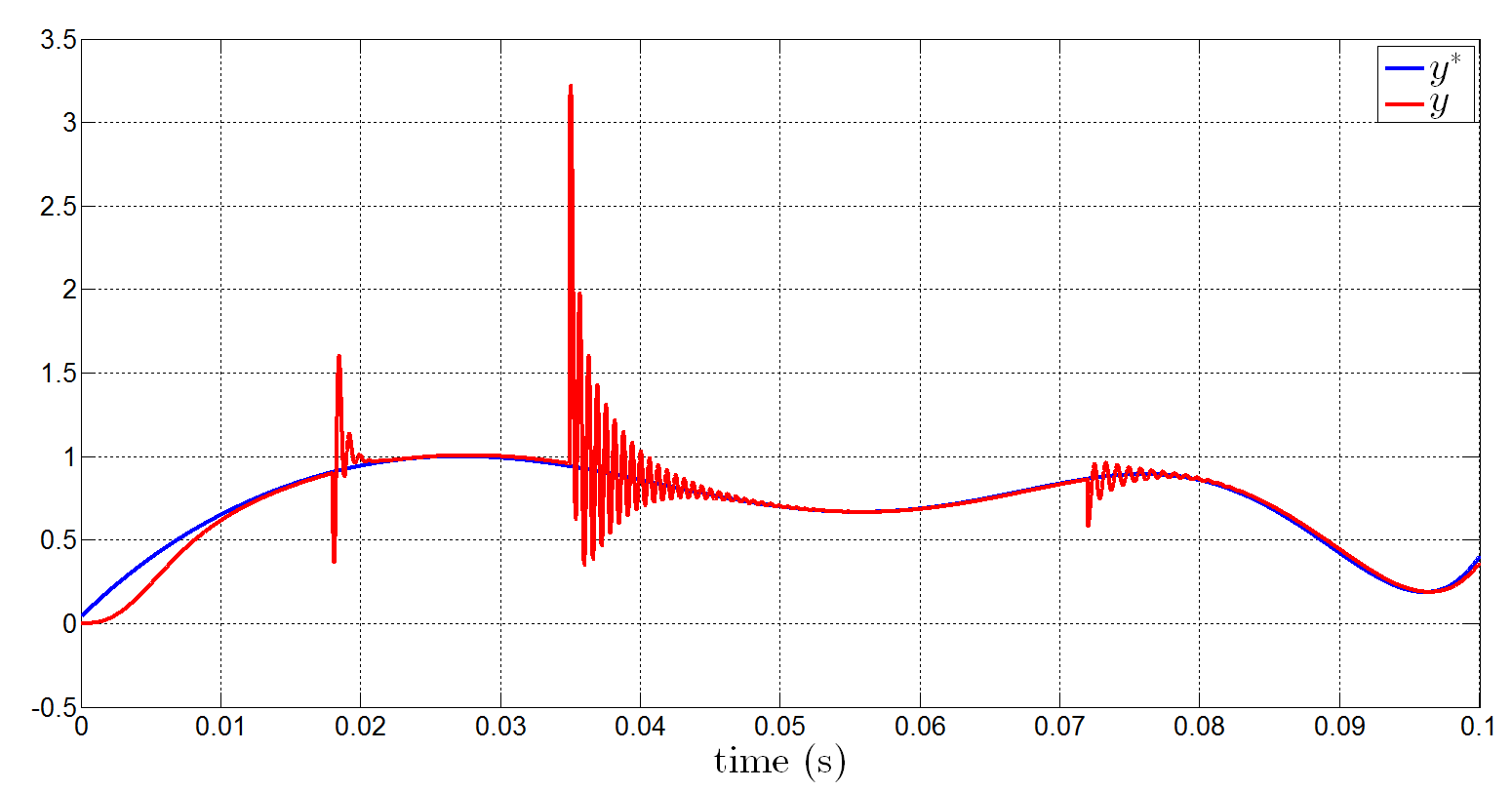} \label{fig:sec}}
\caption{Switched sequence $\# 3$ for $t_{1} = 0.018$ s, $t_2 = 0.035$ s and $t_3 = 0.072$ s.}
\label{fig:Fig_3_seq}
\end{figure*}

\begin{figure*}[!h] 
\centering
\subfigure[\footnotesize Step responses of each system $\Sigma_i$.]{\includegraphics[width=16cm]{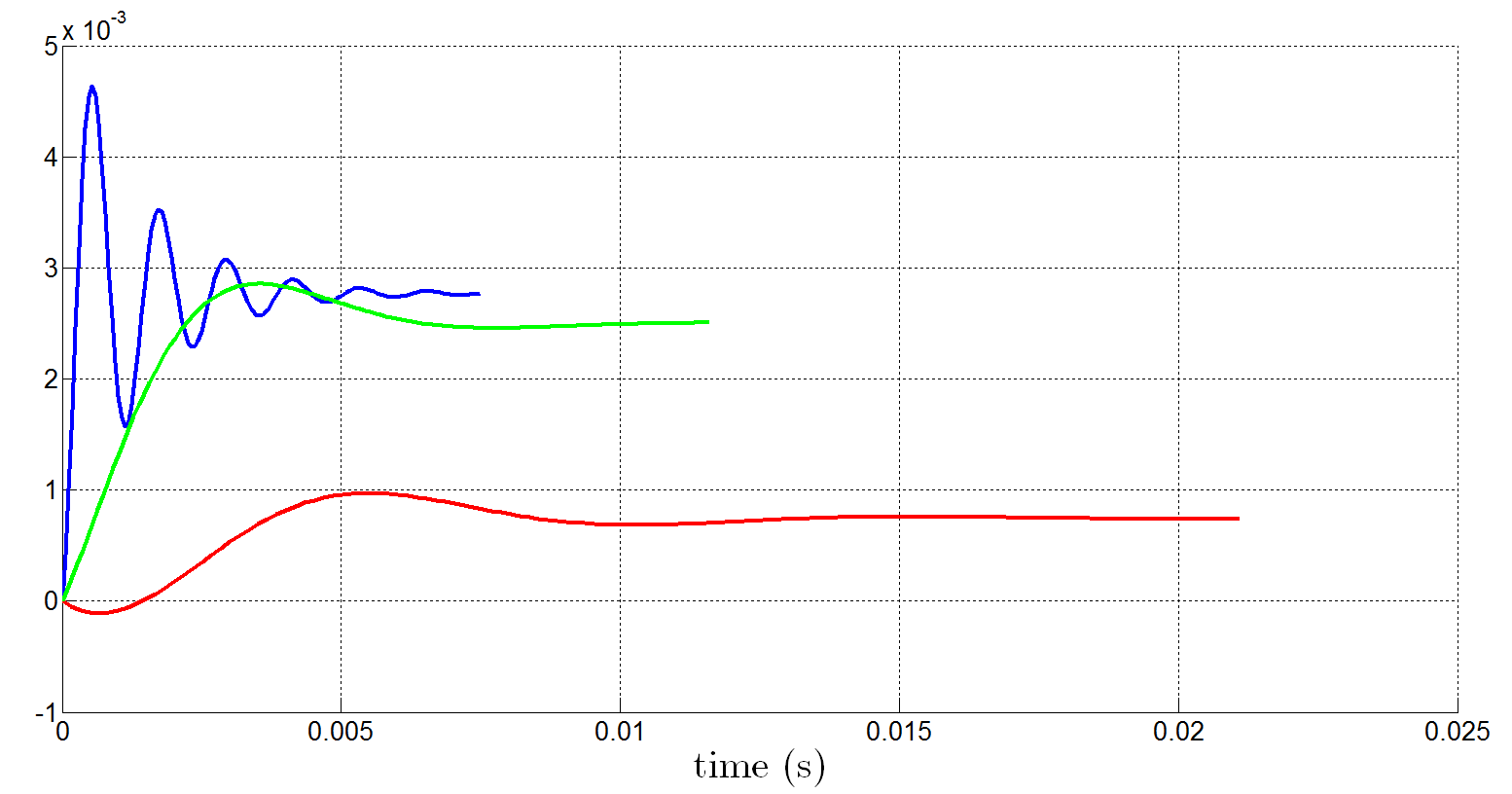} \label{fig:first}}
\subfigure[\footnotesize Controlled switched sequence.]{\includegraphics[width=16.2cm]{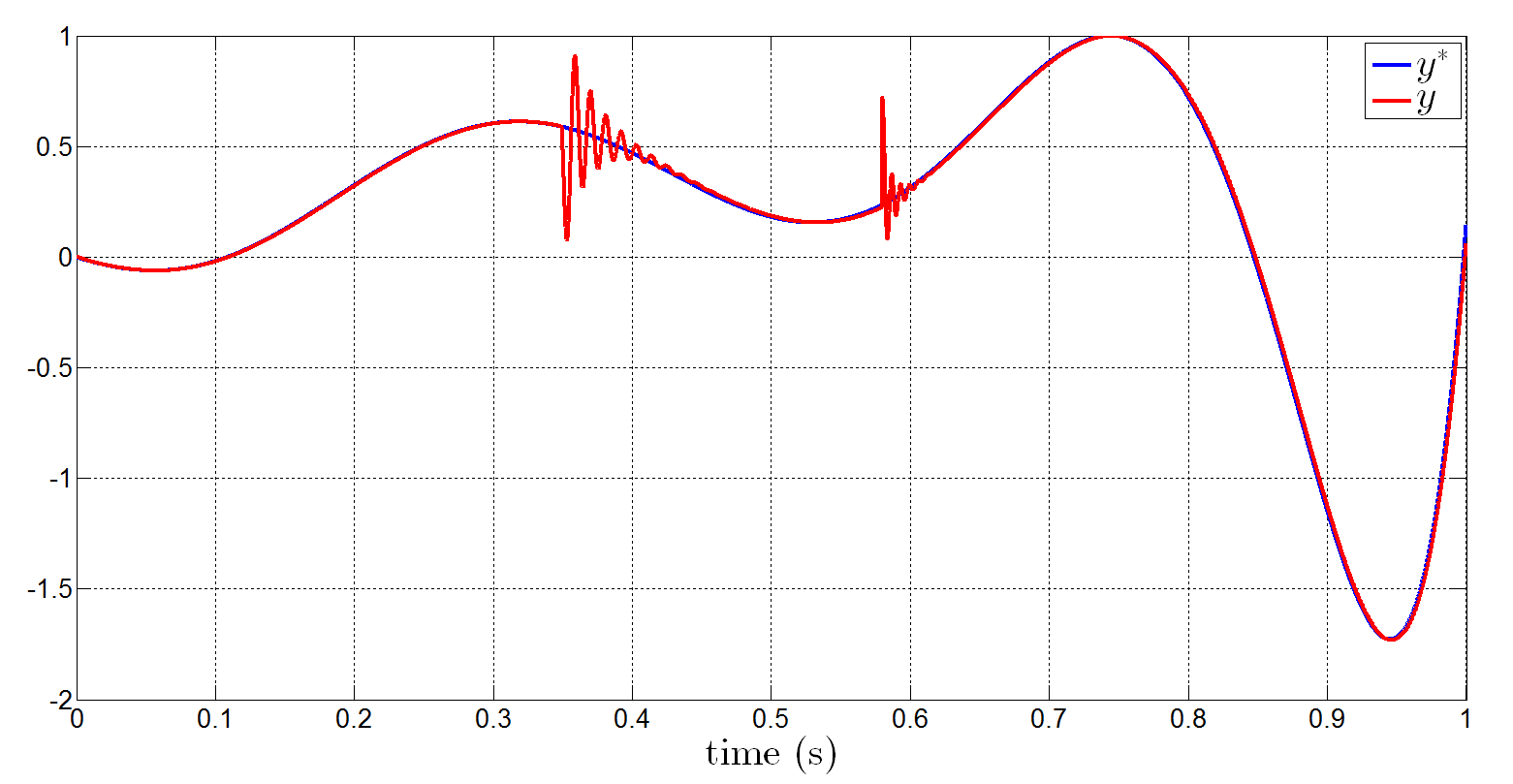} \label{fig:sec}}
\caption{Switched sequence $\# 4$ for $t_{1} = 0.35$ s, $t_2 = 0.58$ s.}
\label{fig:Fig_4_seq}
\vspace{0.5cm}
\end{figure*}

\noindent
Figures \ref{fig:Fig_1_seq} \ref{fig:Fig_2_seq} \ref{fig:Fig_3_seq} \ref{fig:Fig_4_seq} present some examples of the application of the $\mathcal{C}_{\pi}$-control under different 
arbitrary switching sequences that involve 
both minimum and non-minimum phase systems. The first switching time is $t_{1}$, the second is $t_{2}$ and the third is $t_{3}$.

Consider now the existence of a delay $\tau$ on $y$ that modify (\ref{eq:sys_linear}) such that:

\begin{equation}\label{eq:sys_linear_LD}
\Sigma_p( u \longmapsto y) := \left\{ \begin{array}{l}
\dot{x}(t) = A_p x(t)  + B_p u(t) \\
 y = C_p x(t- \tau)
      \end{array} \right.
\end{equation}

\begin{figure*}[!h] 
\centering
\subfigure[\footnotesize Step responses of each system $\Sigma_i$.]{\includegraphics[width=16cm]{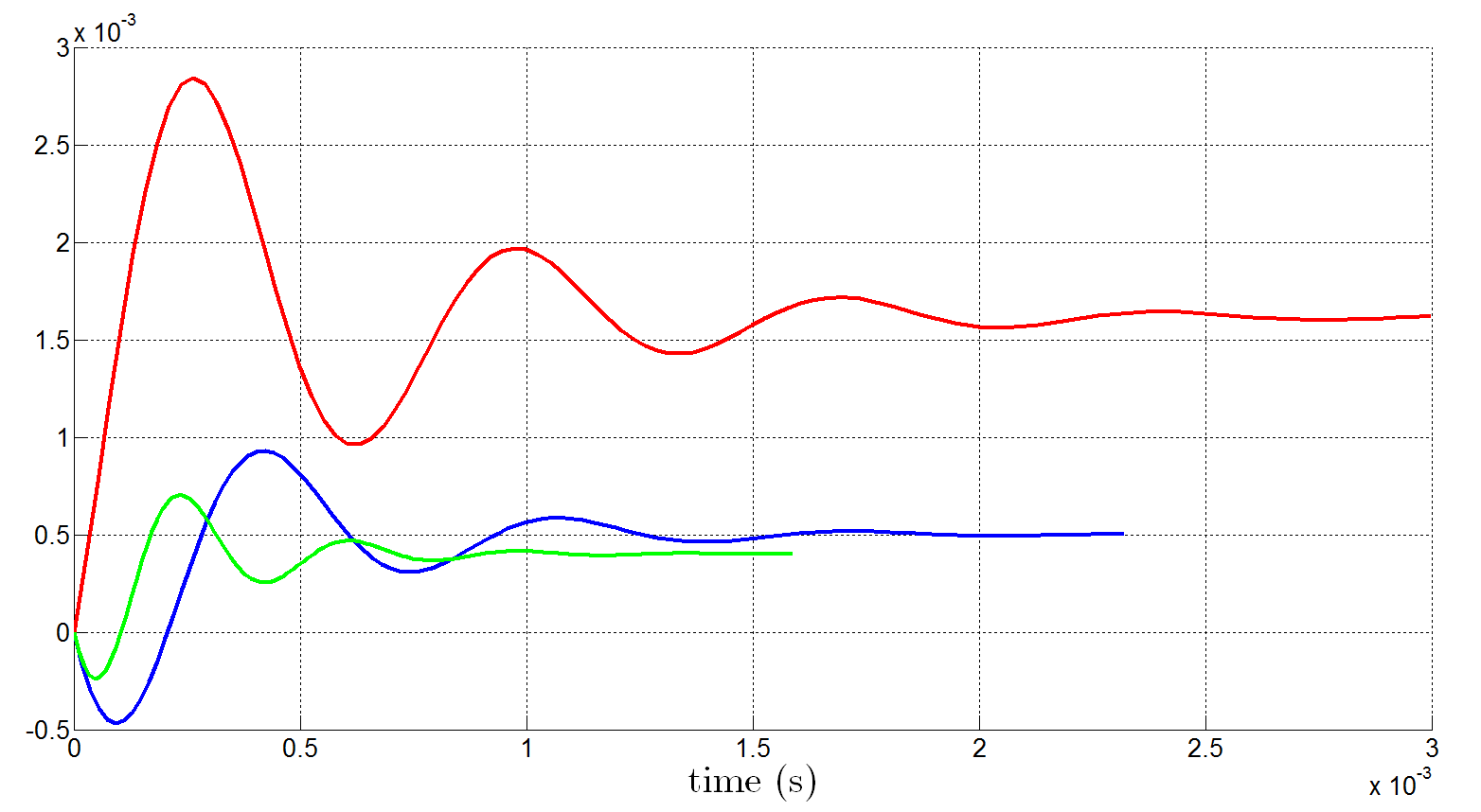} \label{fig:first}}
\subfigure[\footnotesize Controlled switched sequence.]{\includegraphics[width=16.2cm]{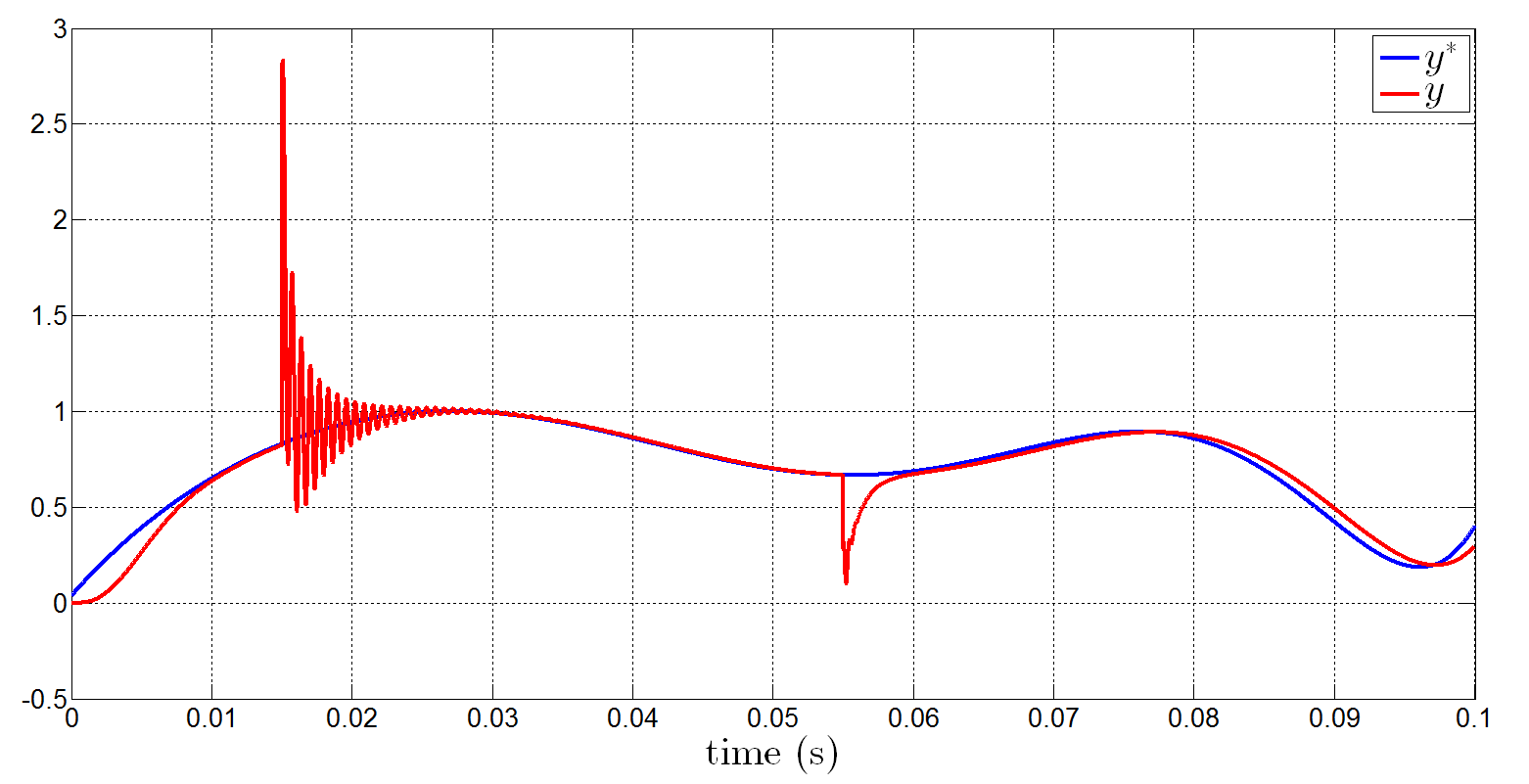} \label{fig:sec}}
\caption{Switched sequence $\# 5$ for $t_{1} = 0.015$ s, $t_2 = 0.055$ s. An addition of a time-delay occurs at $t = 0.06$ s.}
\label{fig:Fig_5_seq}
\end{figure*}

\begin{figure*}[!b] 
\centering
\subfigure[\footnotesize Step responses of each system $\Sigma_i$.]{\includegraphics[width=16cm]{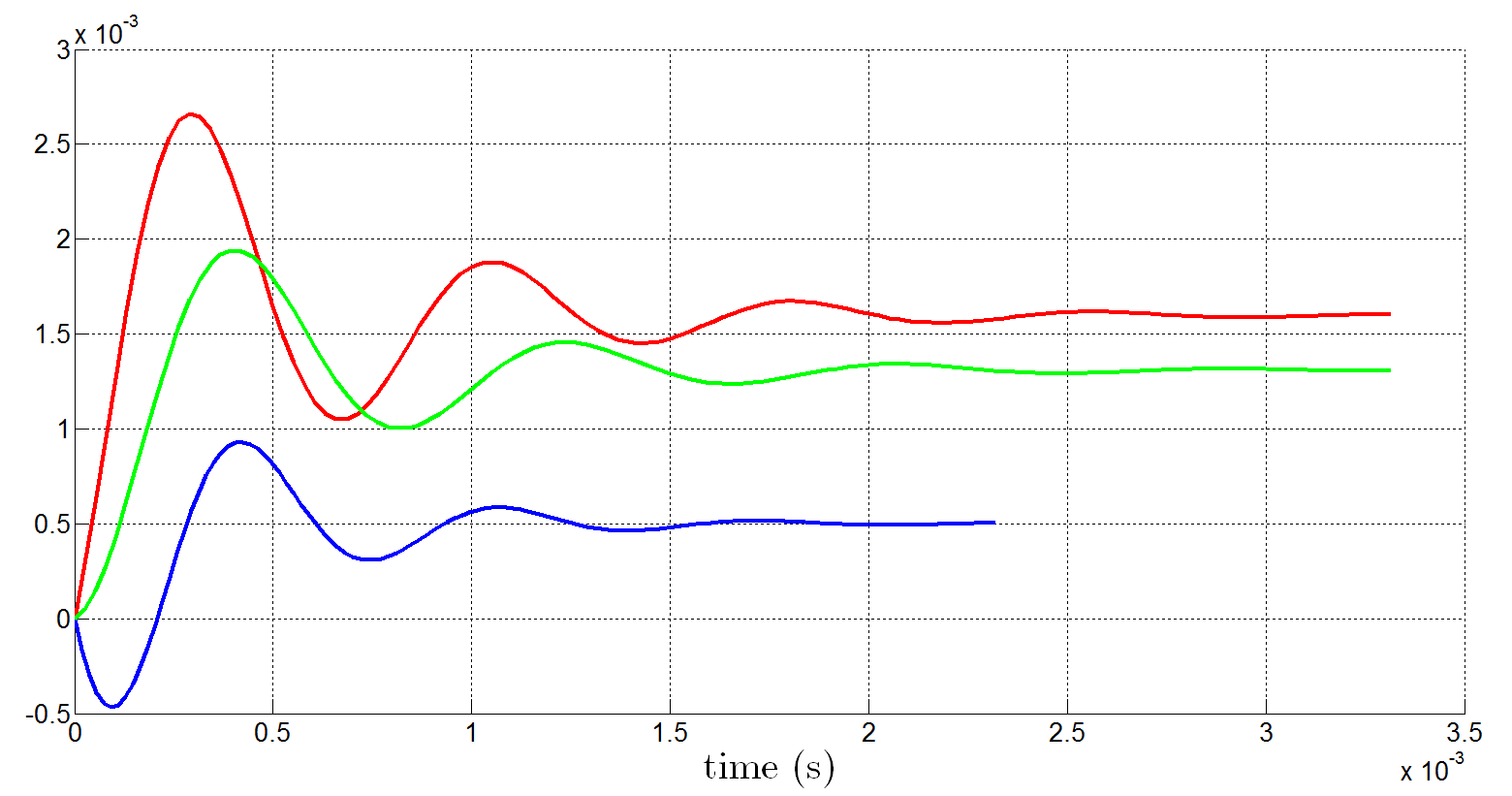} \label{fig:first}}
\subfigure[\footnotesize Controlled switched sequence.]{\includegraphics[width=16.2cm]{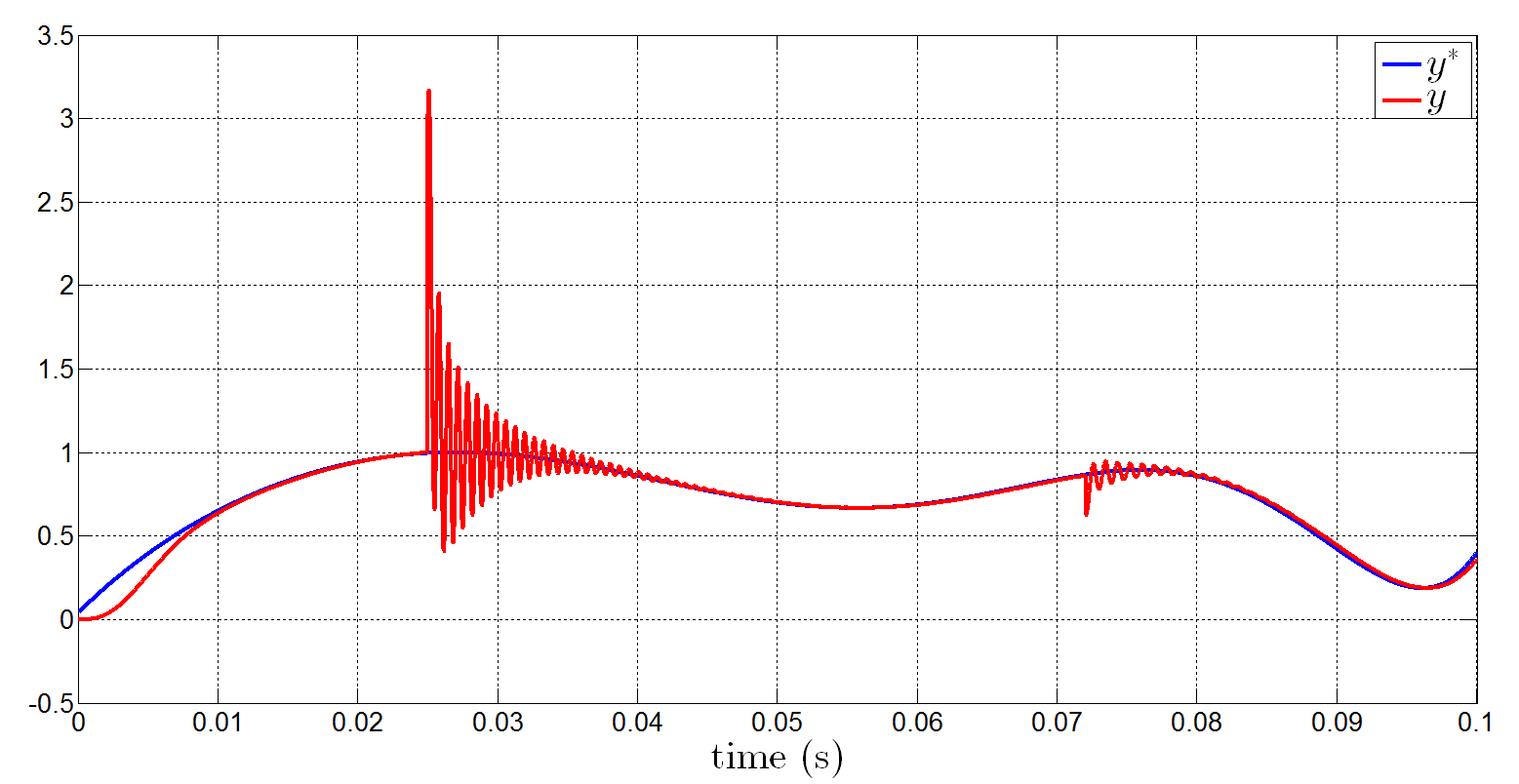} \label{fig:sec}}
\caption{Switched sequence $\# 6$ for $t_{1} = 0.025$ s, $t_2 = 0.072$ s. An addition of a time-delay occurs at $t = 0.06$ s.}
\label{fig:Fig_6_seq}
\end{figure*}

\noindent
This delay can e.g. simulate the propagation delay inside a sensor network. Figures \ref{fig:Fig_5_seq} and \ref{fig:Fig_6_seq} present two examples of the application of the $\mathcal{C}_{\pi}$-control 
under different switching sequences that involve both minimum and non-minimum phase systems.


\clearpage

\subsubsection*{Discussions}

These results have been obtained with a single specific set of the $\mathcal{C}_{\pi}$-parameters; to improve the tracking performances, one may consider an on-line adjustment 
of the $\mathcal{C}_{\pi}$-parameters. Although resonances 
occur at the instants of switches, the stability of the control is preserved (regarding the studied cases) when switching from the different types of systems. In the same manner, 
a direct tuning of the $\mathcal{C}_{\pi}$-parameters could damp (ideally, could cancel) the resonant effects.

\noindent
The presented simulation results show that the proposed control law is robust to "strong" model variations and in particular when the model is a switching non-minimum phase or
minimum phase system that include eventually time-delay. Moreover, the proposed control law seems to have the same properties than the original model-free control \cite{esta} \cite{fliess2} 
for which its performances have been successfully verified especially in simulation. 

\subsection{Ballistic and the fire control}

If one fire a projectile at an initial angle and an initial speed, then general physics allows calculating how far it will travel... but would it be possible to {\it control} the initial speed 
needed in such manner that the projectile reaches a precise distance? That's what we propose to do using the $\mathcal{C}_{\pi}$-control.

\subsubsection{Ballistic simplified model}

We define first a simple model of the trajectory $\mathbf{w}(t) = (w_x(t), \, w_z(t))$ of a projectile of mass {\it m} in the usual frame of reference $(Oxz)$ fired with an initial 
speed magnitude ${v_0}$ that makes a fire angle $\theta$ with the horizontal reference i.e. $\mathbf{v}(t) = v_0 \cos (\theta) \mathbf{e}_x + v_0 \sin (\theta) \mathbf{e}_z$. 
The origin $(0,0)$ of the frame reference is considered as the initial position of the projectile.

Denote $\mathbf{a}$ the acceleration vector and $\mathbf{v}$ the speed vector of the projectile. From Newton law, considering the action of the gravity $\mathbf{g}$ and 
the air resistance $c \mathbf{v}^2$, we have:

\begin{equation}\label{eq:dyn_gen}
  m \frac{d^2 \, \mathbf{w}(t) }{d \, t^2} = m \mathbf{a} = m \mathbf{g} - c \mathbf{v}^2
\end{equation}

\noindent
with the initial conditions:

\begin{equation}
\left. \frac{d \, w_x(t)}{d \, t} \right|_{(0,0)} = v_0 \cos \theta, \qquad \left. \frac{d \, w_z(t)}{d \, t} \right|_{(0,0)} = v_0 \sin \theta
\end{equation}

\noindent
Considering no air resistance i.e. $c = 0$, (\ref{eq:dyn_gen}) is simplified:

\begin{equation}\label{eq:sol_dyn}
\left\{\begin{array}{c}
\displaystyle{ m \frac{d \, w_x(t)}{d \, t}  = 0 } \\ 
\\[0.01cm]
\displaystyle{ m \frac{d \, w_z(t)}{d \, t}  = - m g }
\end{array} \right.
\end{equation}

\noindent
whose solution reads:

\begin{equation}\label{eq:sol_dyn_sol}
\left\{\begin{array}{c}
 \displaystyle{ w_x(t) = v_0 \cos \theta \, t + cte }\\ 
\\[0.01cm]
\displaystyle{ w_z(t) = v_0 \sin \theta \, t - \frac{1}{2} g t^2 + cte }
\end{array} \right.
\end{equation}

From (\ref{eq:sol_dyn_sol}), the range (or the target) of the projectile $x_d$ at $z = 0$ can be easily deduced. We have:

\begin{equation}\label{eq:target}
 d = \frac{2 v_0^2 \cos \theta \sin \theta}{g}.
\end{equation}

\subsubsection{Ballistic-fire control methodology}

\paragraph{Proposed strategy}

Consider by hypothesis that a "virtual" trajectory $\mathbf{w^*}$ of the projectile hits a target $x_d^*$ from an initial speed $v_0^*$ and an initial fire angle $\theta^*$.

\begin{figure}[!h]
  \centering
    \subfigure[\footnotesize Case $c = 0$] {\includegraphics[width=16cm]{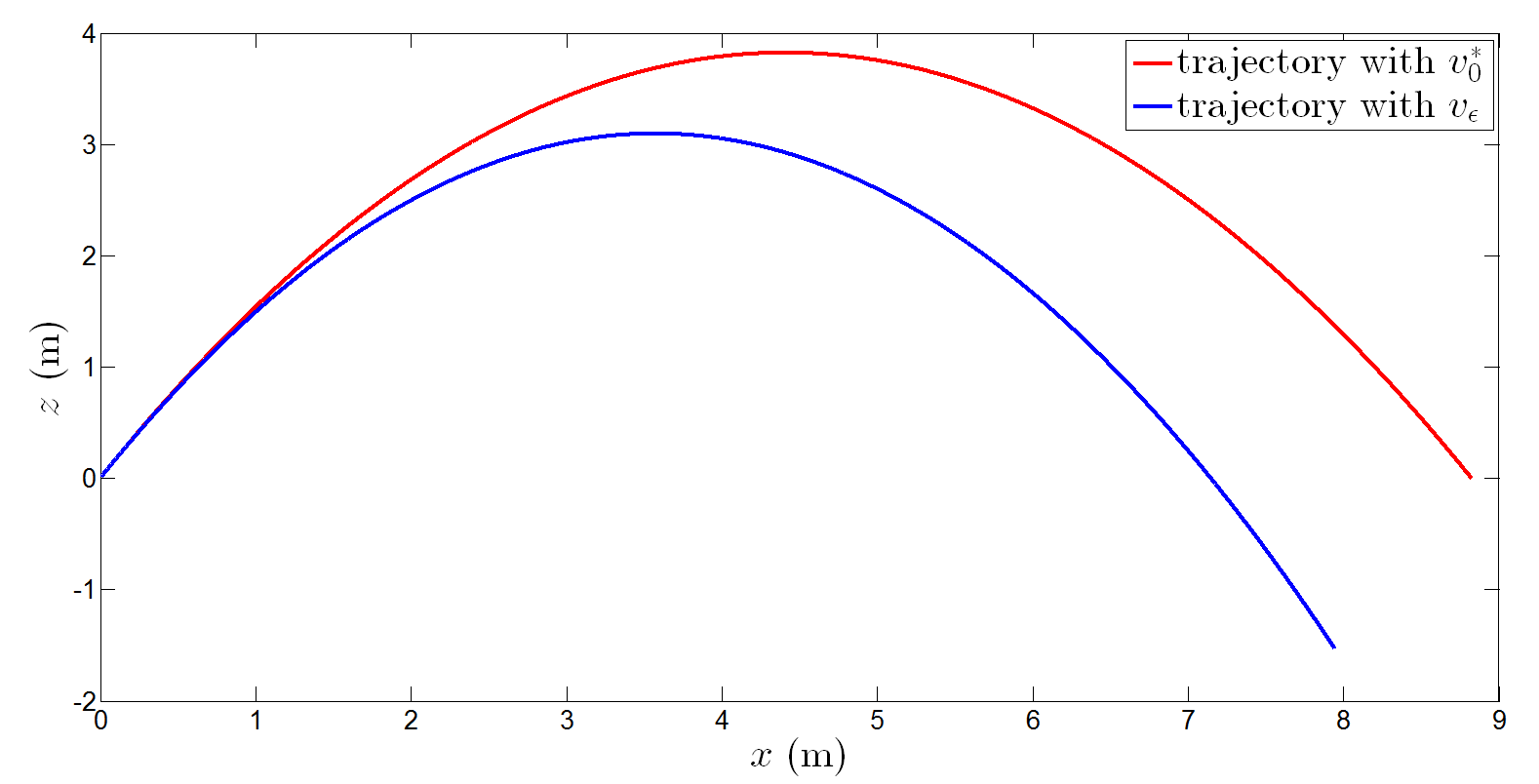} \label{fig:fig_1_b}}
    \subfigure[\footnotesize Case $c = 0$ and $c \neq 0$]{\includegraphics[width=16cm]{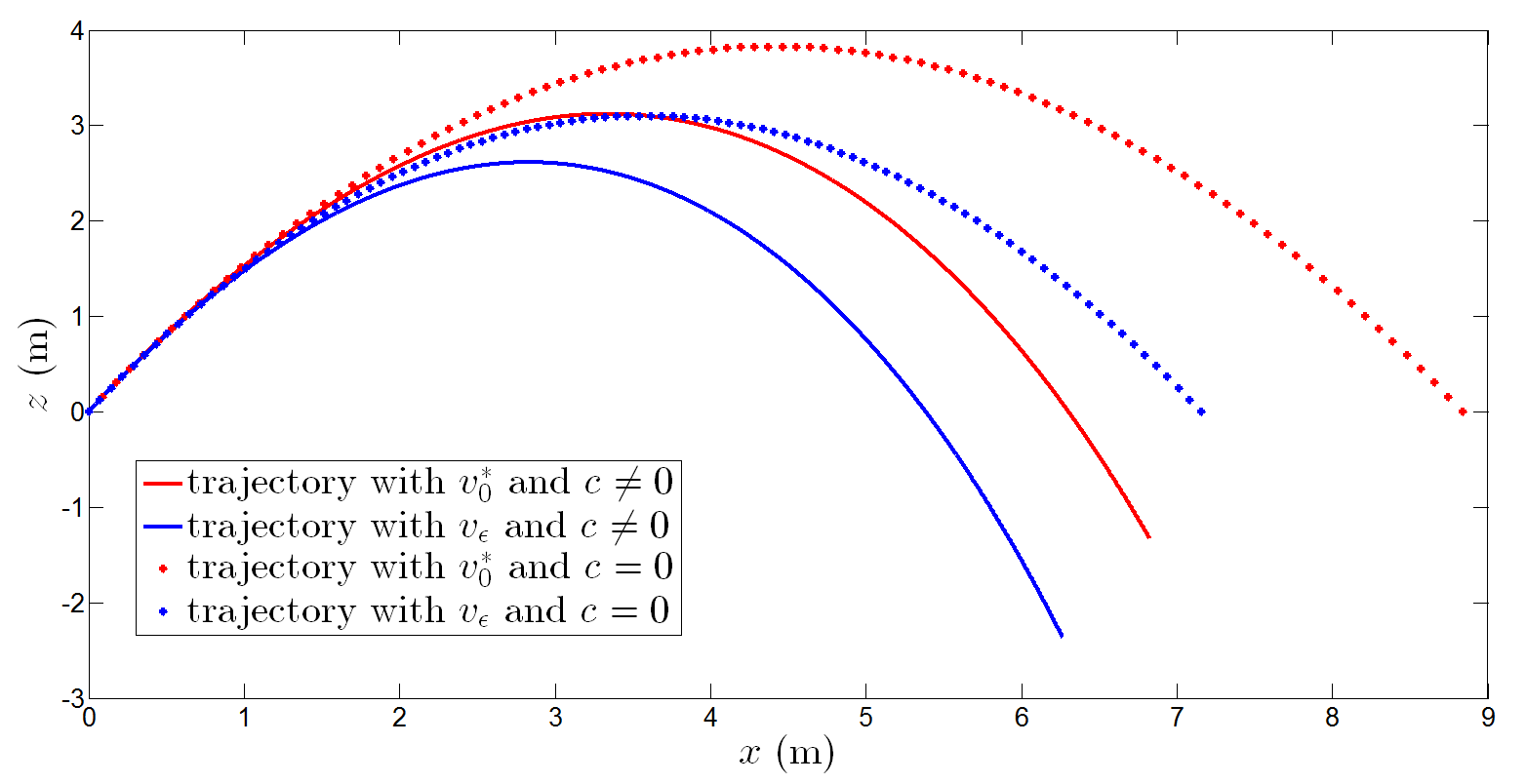} \label{fig:fig_2_b}}
    \caption{Example of comparison of the uncontrolled "true" projectile (subj. to $v_{\varepsilon}$ and $\theta^*$) and the "virtual" trajectory (subj. to $v_0^*$ and $\theta^*$),
    considering $c = 0$ and $c \neq 0$. We took $c = 0.05$, $\theta_{\varepsilon} = \theta^* = \pi/3$, $v_0 = 10$ m/s and $v_{\varepsilon} = 9$ m/s.}
    \label{fig:fig_12_b}
\end{figure}

\noindent
Consider now the "true" projectile to fire with a trajectory $\mathbf{w}$. To hit the target $x_d^*$ (at $z = 0$) from {\it a small initial speed} $v_{\varepsilon}$, one {\it controls } the projectile in 
such manner that the projectile reaches quickly
the initial speed $v_0^*$ and the initial fire angle $\theta^*$ required to hit the specified target $x_d^*$ according to (\ref{eq:target}). 
During such "launching" phase, that we define as the "launching" distance $\Delta x_0$ for which the trajectory of the projectile is fully controlled, we start from an initial condition that prevents 
the 
projectile to reach its target $x_d^*$ i.e. :

\begin{equation}
\left. \frac{d \, w_x(t)}{d \, t} \right|_{(0,0)} = v_{\varepsilon} \cos \theta_{\varepsilon}, \qquad \left. \frac{d \, w_z(t)}{d \, t} \right|_{(0,0)} = v_{\varepsilon} \sin \theta_{\varepsilon}
\end{equation}

\noindent
where $v_{\varepsilon} < v_0^*$ and $\theta_{\varepsilon} \leq \theta^*$ are positive resp. initial speed magnitude and fire angle. 
Figure \ref{fig:fig_12_b} illustrates simulation examples of a "virtual" trajectory (subj. to $v_0^*$ and $\theta^*$) and a "true" trajectory (subj. to $v_{\varepsilon}$ and $\theta^*$) that is not 
controlled; the simulations of the trajectories considering $c \neq 0$ are presented in Fig. \ref{fig:fig_2_b}.

The goal is to accelerate the projectile using a specific mechanical device in such manner that the references $v_0^*$ and $\theta^*$ 
are reached quickly. {\it We consider therefore controlling the trajectory $\mathbf{w}(t)$ of the projectile  over the launching distance $\Delta x_0$.}

\subsubsection*{Controllable ballistic model}

To simulate a controllable model of the trajectory $\mathbf{w}$ of the projectile, we consider adding an external acceleration force in (\ref{eq:dyn_gen}) that represents the mechanical action of the 
specific mechanical device over the distance $\Delta x_0$. We have:

\begin{equation}\label{eq:dyn_gen_cont_app}
  m \frac{d^2 \, \mathbf{w}(t) }{d \, t^2}  = m \mathbf{g} - c \mathbf{v}^2 + \mathbf{a}^{ext}
\end{equation}

\noindent
where $\mathbf{a}^{ext}(t) = (a^{ext}_x(t), \, a^{ext}_z(t))$ is equivalent to the external force provided by the mechanical device. Since such device acts only over $\Delta x_0$, then we assume that
$\mathbf{a}^{ext} = \mathbf{0}$ for all $x > \Delta x_0$.

\subsubsection{Implementation of the $\mathcal{C}_{\pi}$-controller}

A possible control scheme is to consider controlling the trajectory $\mathbf{w}$ that must be "as close as possible" to the reference $\mathbf{w^*}$ over the distance $\Delta x_0$. 
Therefore, $\mathbf{w}$ is physically measured and the external acceleration $\mathbf{a}^{ext}$ is driven by the $\mathcal{C}_{\pi}$-controller, through the specific mechanical device.

\noindent
We build a closed-loop that creates a feedback between (\ref{eq:iPI_discret_nm_eq}) and (\ref{eq:dyn_gen_cont_app}). We have "symbolically", for all $x \leq \Delta x_0$:

\begin{equation} \label{eq:iPI_discret_nm_eq_app}
    \left\{  \begin{array}{c}
    \displaystyle{\mathbf{u}_k = \mathbf{a}^{ext}_k = \left. \int_0^t K_i (\mathbf{w}^*_{k-1} - \mathbf{w}_{k-1} ) d \, \tau \right|_{k-1} \underbrace{\left\{ \mathbf{u}_{k-1}^i + {K_p} ( \alpha e^{-\beta k} - \mathbf{w}_{k-1}) \right\}}_{\mathbf{u}_{k}^i} } \\
    \\[0.1cm]
    \displaystyle{m \frac{d^2 \, \mathbf{w}(t) }{d \, t^2}  = m \mathbf{g} - c \mathbf{v}^2 + \mathbf{u}_k}
     \end{array} \right.
\end{equation}

\subsubsection*{Determination of the distance $\mathbf{\Delta x_0}^*$} Since we expect that the projectile is fired from $\Delta x_0$ with a speed that is very close to $v_0^*$ (and follows, 
via the $\mathcal{C}_{\pi}$-control, the same trajectory i.e. $w \approx w^*$ over $\Delta x_0$), then, we propose a possible definition of the theoretical launching distance $\Delta x_0^*$,
(considered only over the $x$ axis) that corresponds to the solution in $w_x$ of:

\begin{equation}
\frac{d \, w_x (t)}{d \, t} = v_{0 \, x}^*
\end{equation}

\noindent
Geometrically, the theoretical launching distance $\Delta x_0^*$ is associated to the speed $v_x$ that is reached by the projectile (launched with the initial speed $v_{\varepsilon} < v_0^*$) 
when $v_x$ is close to $v_{0 \, x}^*$.

\newpage

\subsubsection{Numerical simulations}

\subsubsection*{Case $\mathbf{\Delta x_0} > \mathbf{\Delta x_0}^*$}

\begin{figure}[!h]
  \begin{center}
    \subfigure[At the top, controlled trajectory ${w}$ relating to the reference $w^*$; at the bottom,
    calculated speed $d \, {w_x} / d \, t$ relating to the initial speed $v_0 \cos \theta$.]{ \includegraphics[width=16cm]{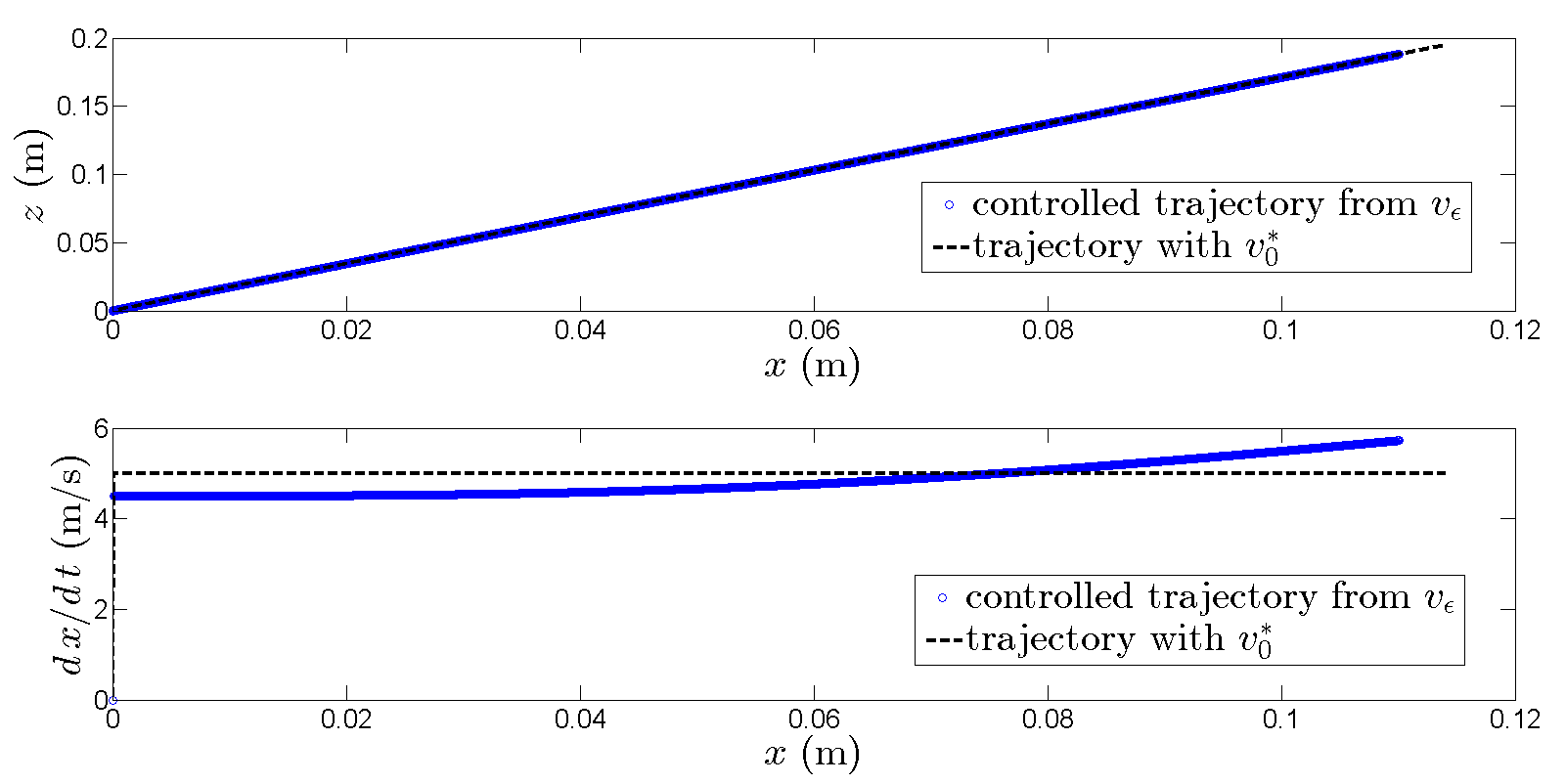}  \label{fig:fig_3_b}} 
    \subfigure[Complete controlled "true" trajectory in comparison with the virtual trajectory.]{\includegraphics[width=16cm]{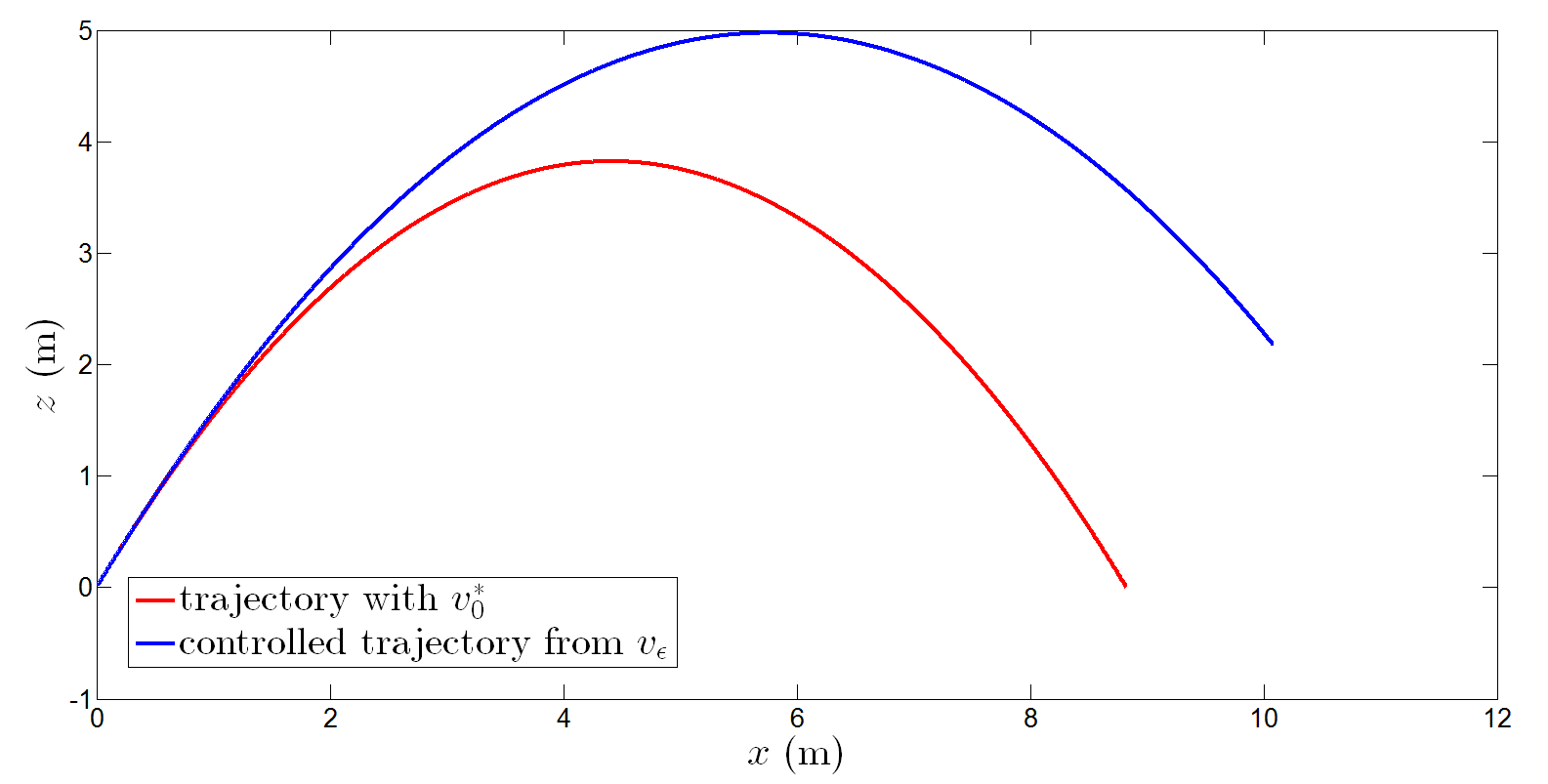} \label{fig:fig_4_b}}
    \caption{Example of controlled trajectory considering $c = 0$ over $\Delta x_0 = 0.11$ m. We took $\theta_{\varepsilon} = \theta^* = \pi/3$, $v_0 = 10$ m/s and $v_{\varepsilon} = 9$ m/s.}
    \label{fig:fig_34_b}
  \end{center}
\end{figure}

\begin{figure}[!b]
  \begin{center}
    \subfigure[At the top, controlled trajectory ${w}$ relating to the reference $w^*$; at the bottom,
    calculated speed $d \, {w_x} / d \, t$ relating to the initial speed $v_0 \cos \theta$.] {\includegraphics[width=16cm]{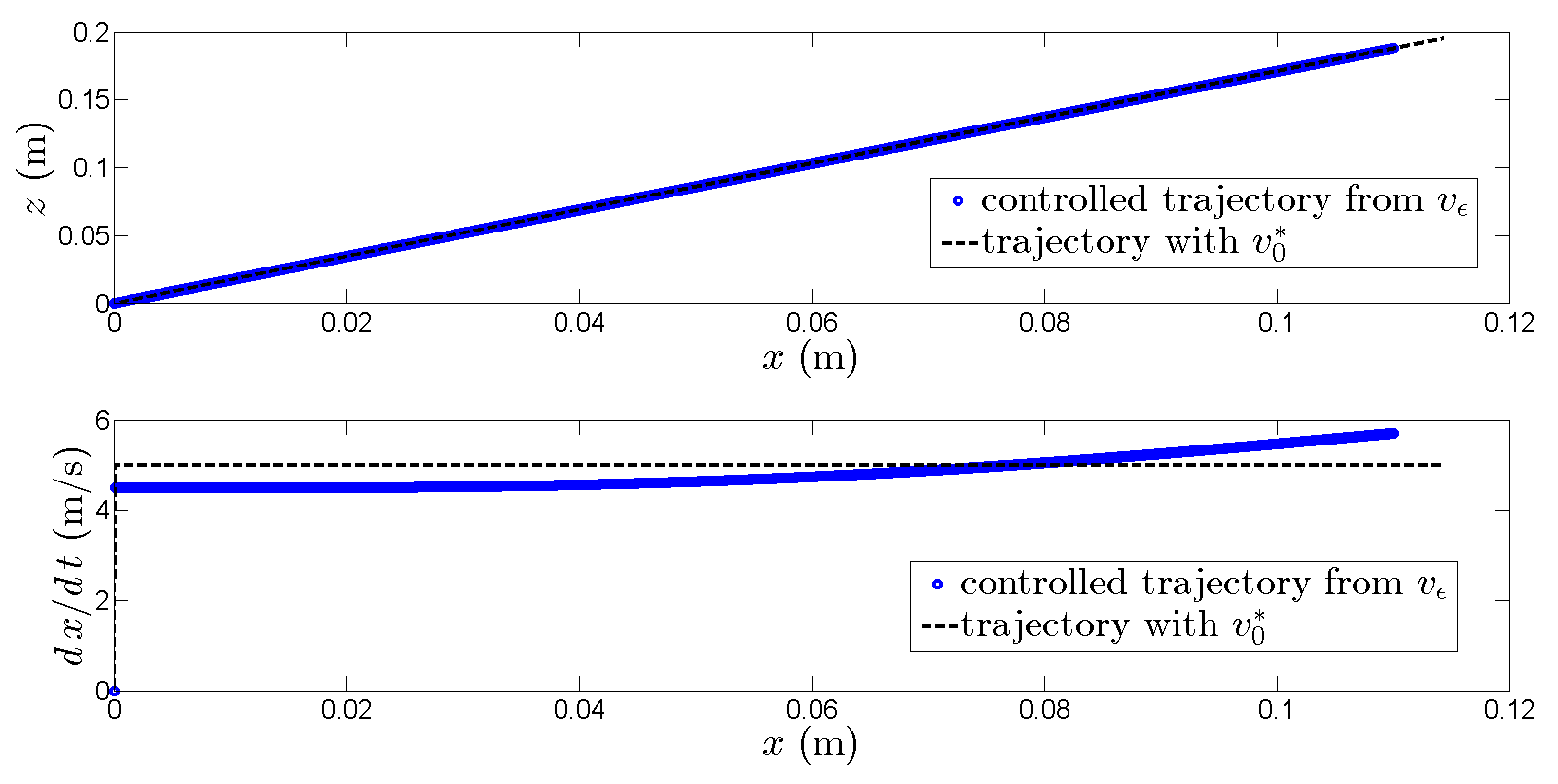} \label{fig:fig_5_b}}
    \subfigure[Complete controlled "true" trajectory in comparison with the virtual trajectory.]{\includegraphics[width=16cm]{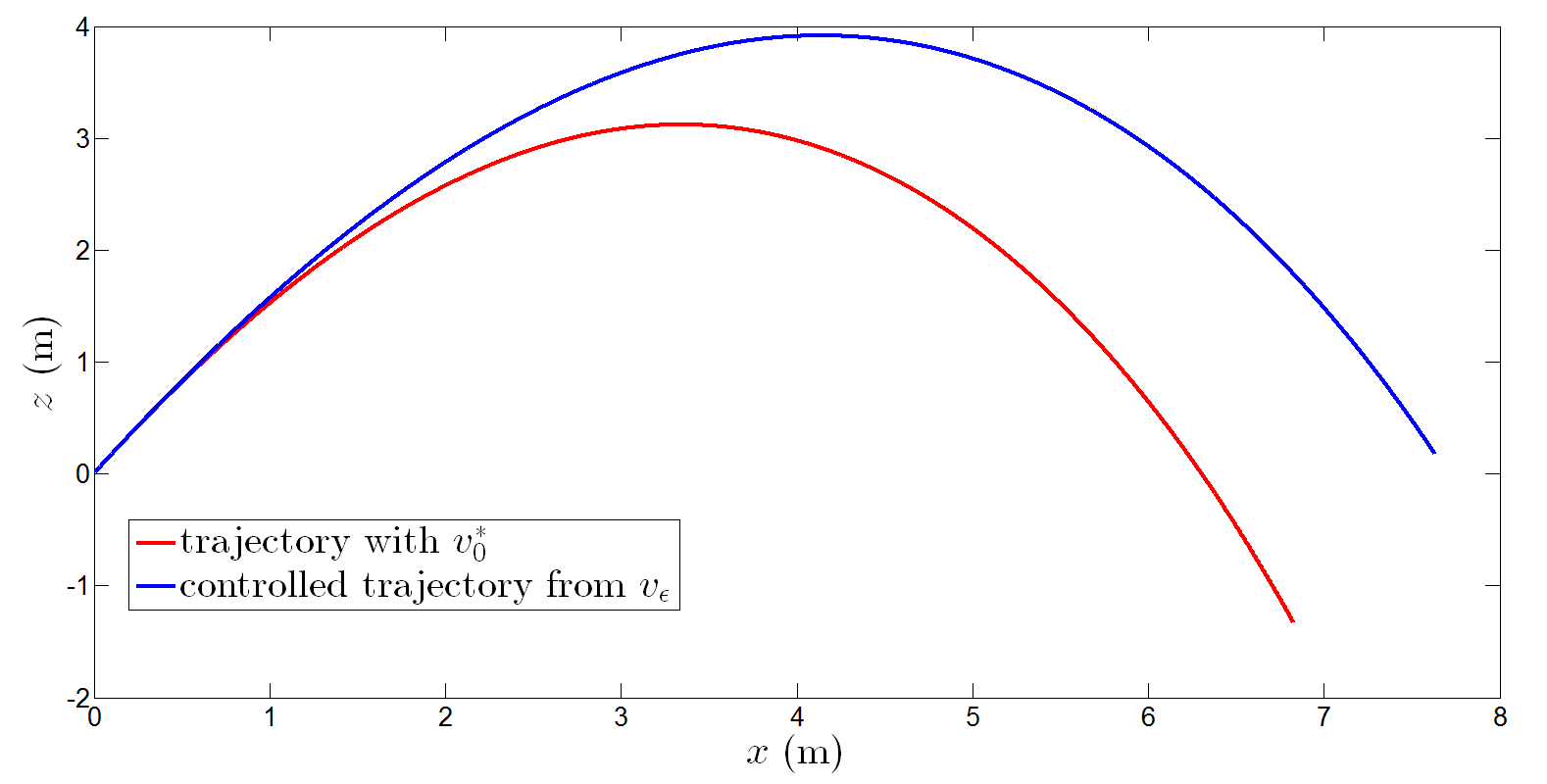}\label{fig:fig_6_b}}
    \caption{Example of controlled trajectory considering $c \neq 0$ over $\Delta x_0 = 0.11$ m. We took $c = 0.05$, $\theta_{\varepsilon} = \theta^* = \pi/3$, $v_0 = 10$ m/s and $v_{\varepsilon} = 9$ m/s.}
    \label{fig:fig_56_b}
  \end{center}
\end{figure}

Consider the simulated "virtual" trajectory, presented in Fig. \ref{fig:fig_12_b}, as the control reference ${w^*}$; to simplify, we consider $\theta$ as constant. 
Figure \ref{fig:fig_34_b} presents the case where the fire is controlled 
considering $c = 0$ 
over $\Delta x_0 = 0.11$ m. In particular, Fig. \ref{fig:fig_3_b} presents, at the top, the evolution of the controlled trajectory ${w}$ in comparison with the reference $w^*$, 
and, at the bottom, the calculated speed $d \, {w_x} / d \, t$ in comparison with the initial speed $v_0 \cos \theta$. Figure \ref{fig:fig_4_b} presents the complete "true" 
controlled trajectory in comparison with the virtual trajectory. Figure \ref{fig:fig_56_b} presents the same simulations in the case $c \neq 0$. 

\subsubsection*{Case $\mathbf{\Delta x_0} \sim \mathbf{\Delta x_0}^*$} 

\begin{figure}[!h]
  \begin{center}
    \subfigure[At the top, controlled trajectory ${w}$ relating to the reference $w^*$; at the bottom,
    calculated speed $d \, {w_x} / d \, t$ relating to the initial speed $v_0 \cos \theta$.]{\includegraphics[width=16cm]{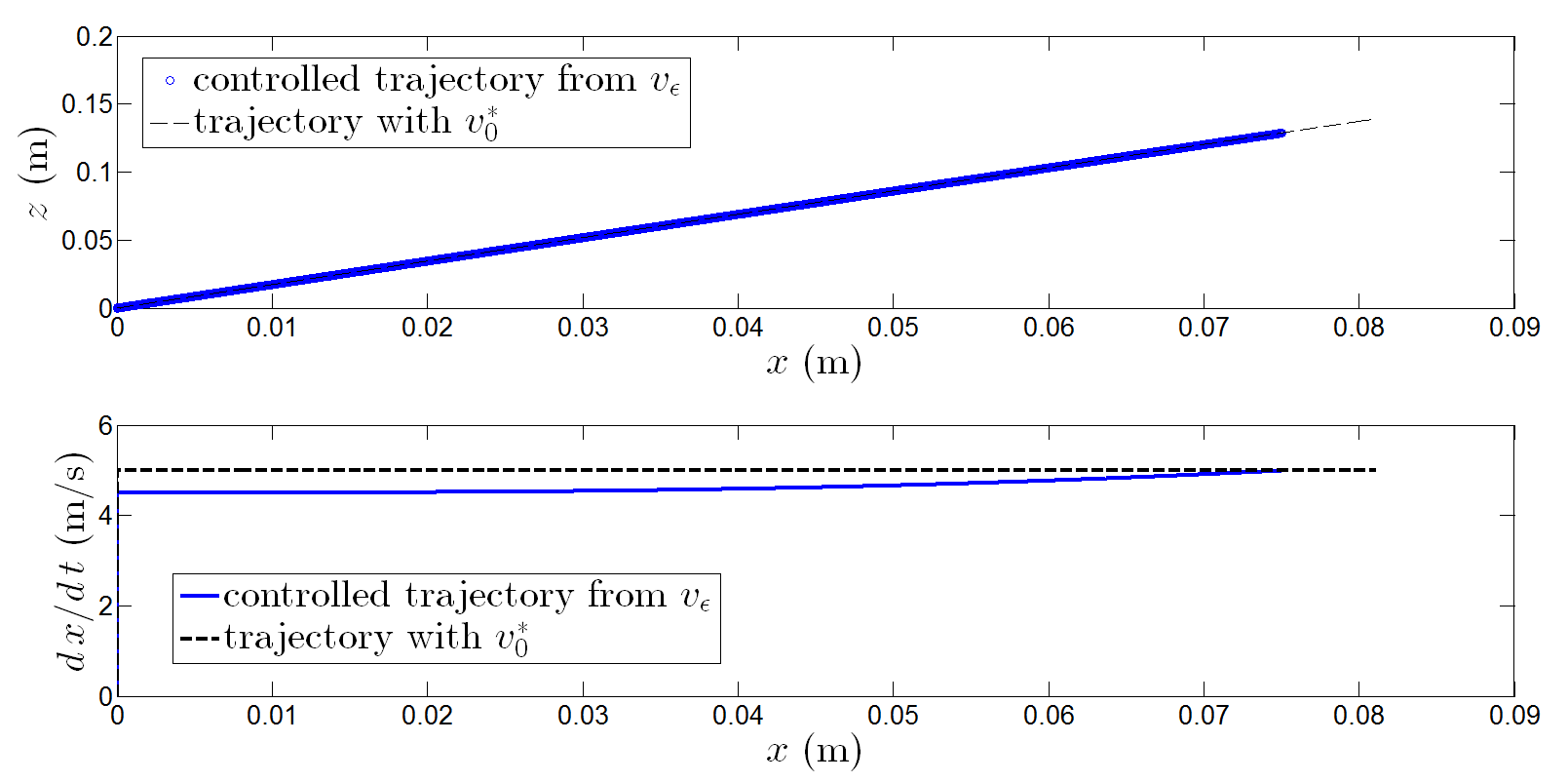} \label{fig:fig_7_b}}
    \subfigure[Complete controlled "true" trajectory in comparison with the virtual trajectory.]{\includegraphics[width=16cm]{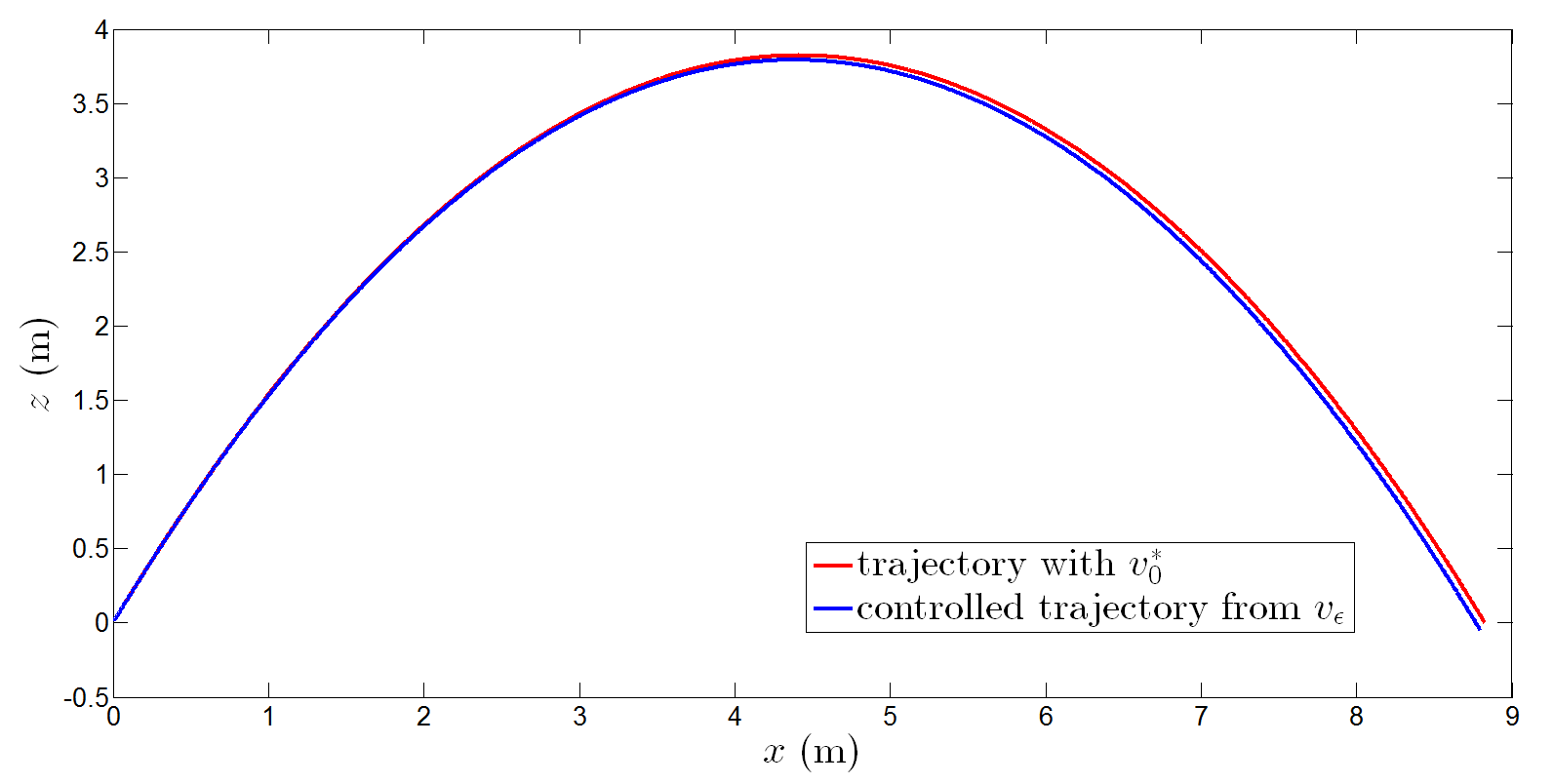}\label{fig:fig_8_b}}
    \caption{Example of controlled trajectory considering $c = 0$ over $\Delta x_0^*$. We took $\theta_{\varepsilon} = \theta^* = \pi/3$, $v_0 = 10$ m/s and $v_{\varepsilon} = 9$ m/s.}
    \label{fig:fig_78_b}
  \end{center}
\end{figure}

\begin{figure}[!b]
  \begin{center}
    \subfigure[At the top, controlled trajectory ${w}$ relating to the reference $w^*$; at the bottom,
    calculated speed $d \, {w_x} / d \, t$ relating to the initial speed $v_0 \cos \theta$.]{\includegraphics[width=16cm]{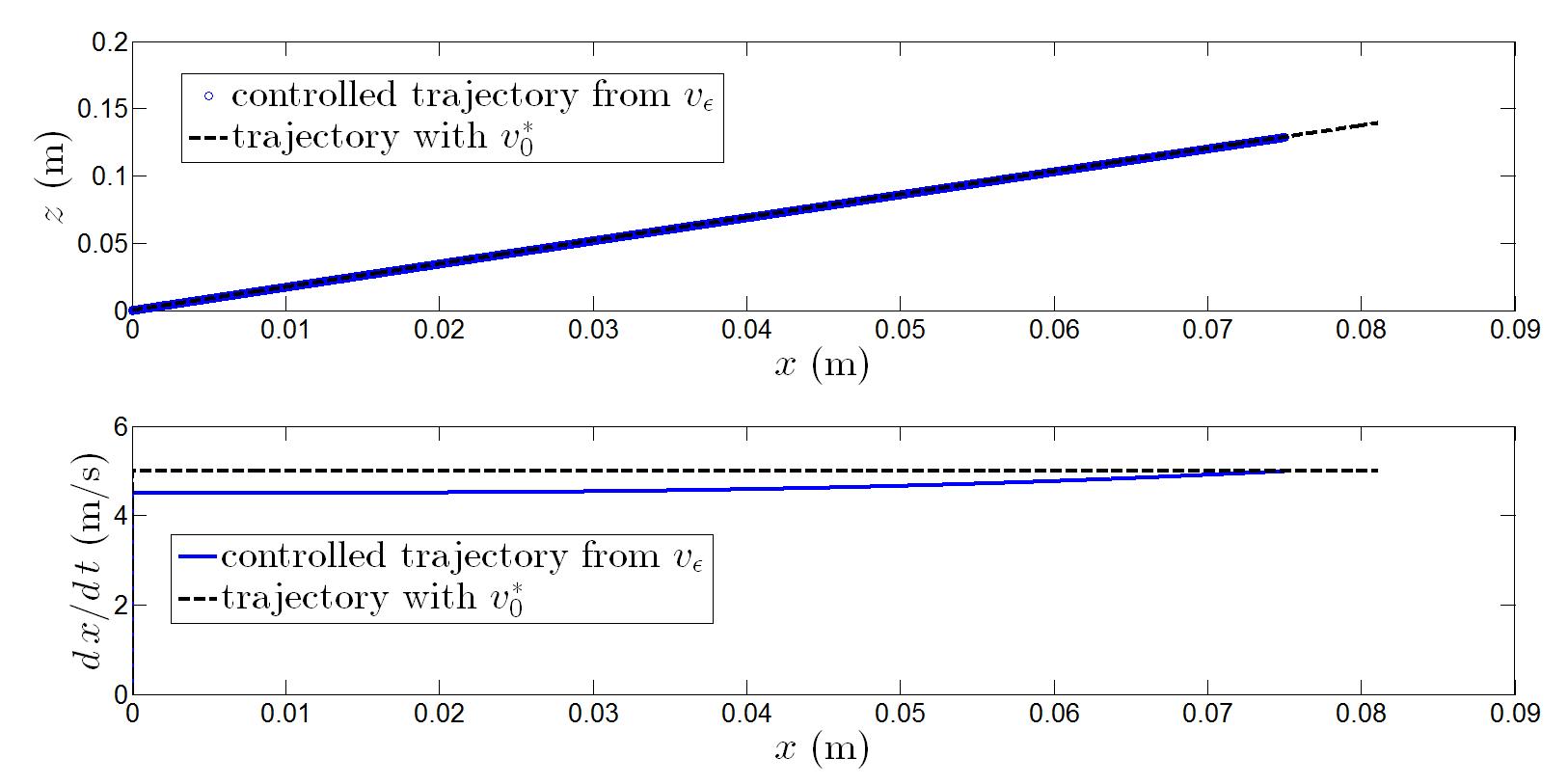} \label{fig:fig_9_b}}
    \subfigure[Complete controlled "true" trajectory in comparison with the virtual trajectory.]{\includegraphics[width=16cm]{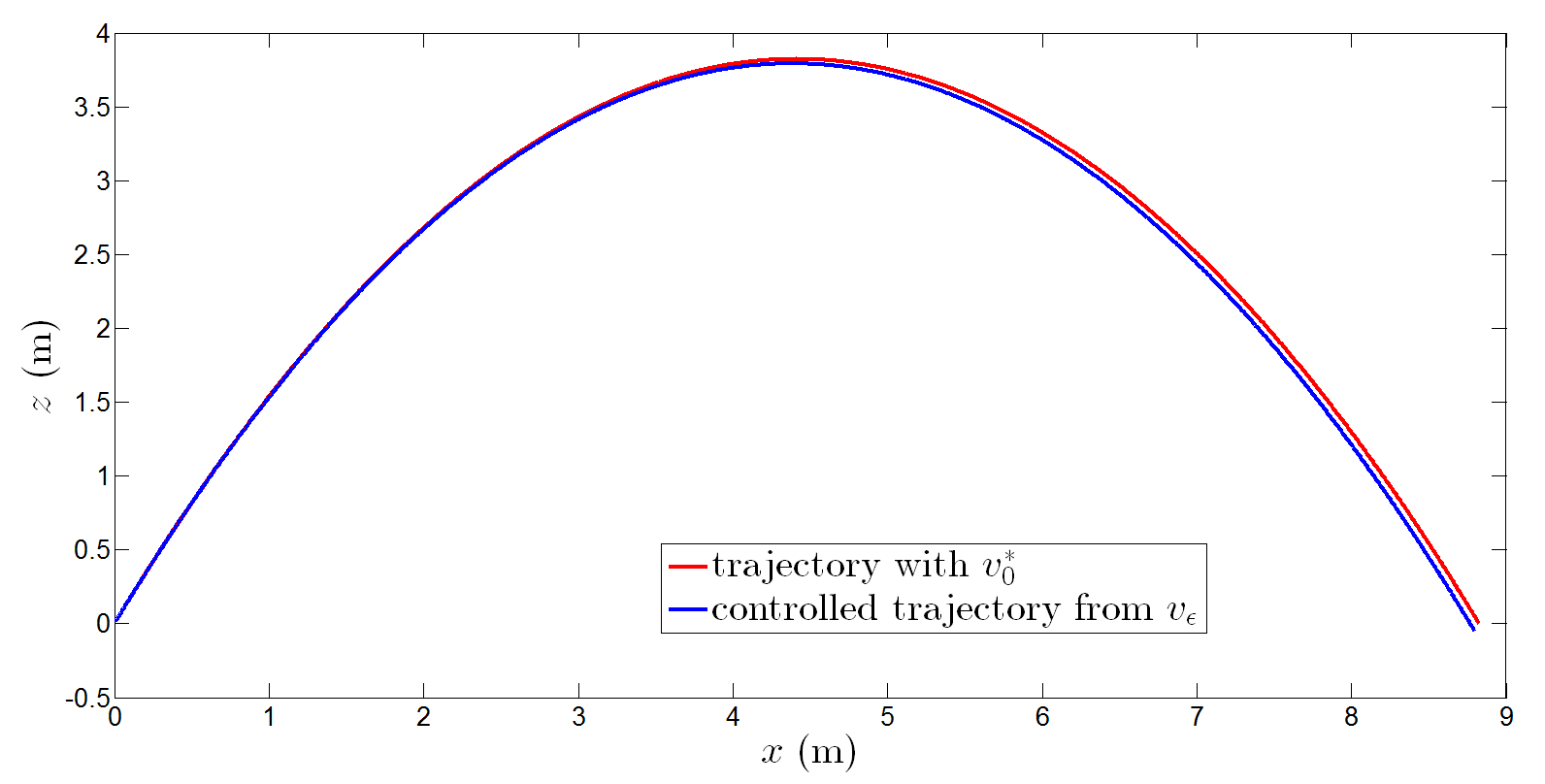} \label{fig:fig_10_b}}
    \caption{Example of controlled trajectory considering $c \neq 0$ over $\Delta x_0^*$. We took $c = 0.05$, $\theta_{\varepsilon} = \theta^* = \pi/3$, $v_0 = 10$ m/s and $v_{\varepsilon} = 9$ m/s.}
    \label{fig:fig_910_b}
  \end{center}
\end{figure}

Consider the simulated "virtual" trajectory, presented in Fig. \ref{fig:fig_12_b}, as the control reference ${w^*}$; to simplify, 
we consider $\theta$ as constant. 
Figure \ref{fig:fig_78_b} presents the case where the fire is controlled 
considering $c = 0$ 
over $\Delta x_0^*$. In particular, Fig. \ref{fig:fig_7_b} presents, at the top, the evolution of the controlled trajectory ${w}$ in comparison with the reference $w^*$, 
and, at the bottom, the calculated speed $d \, {w_x} / d \, t$ in comparison with the initial speed $v_0 \cos \theta$. Figure \ref{fig:fig_8_b} presents the complete "true" 
controlled trajectory in comparison with the virtual trajectory. Figure \ref{fig:fig_910_b} presents the same simulations in the case $c \neq 0$.

\subsubsection*{Discussions} These results have been obtained using the same set of the $\mathcal{C}_{\pi}$-parameters. The properties of stabilization of the control law, like in the previous case when 
dealing with switching systems (\S \ref{para_sw_sys}) , seem to be preserved and ensure good tracking performances in particular when considering $c = 0$ and $c \neq 0$.

Further generalizations would allow using multiple and parallel $\mathcal{C}_{\pi}$-controllers in order to control simultaneous physical quantities. In particular, the 
speed profile $\mathbf{v}$ could be controlled simultaneously with the trajectory $\mathbf{w}$...


\clearpage

\clearpage


\subsection{Control of the HIV-1 model}

\begin{figure}[!h]
\centering
\includegraphics[width=13cm]{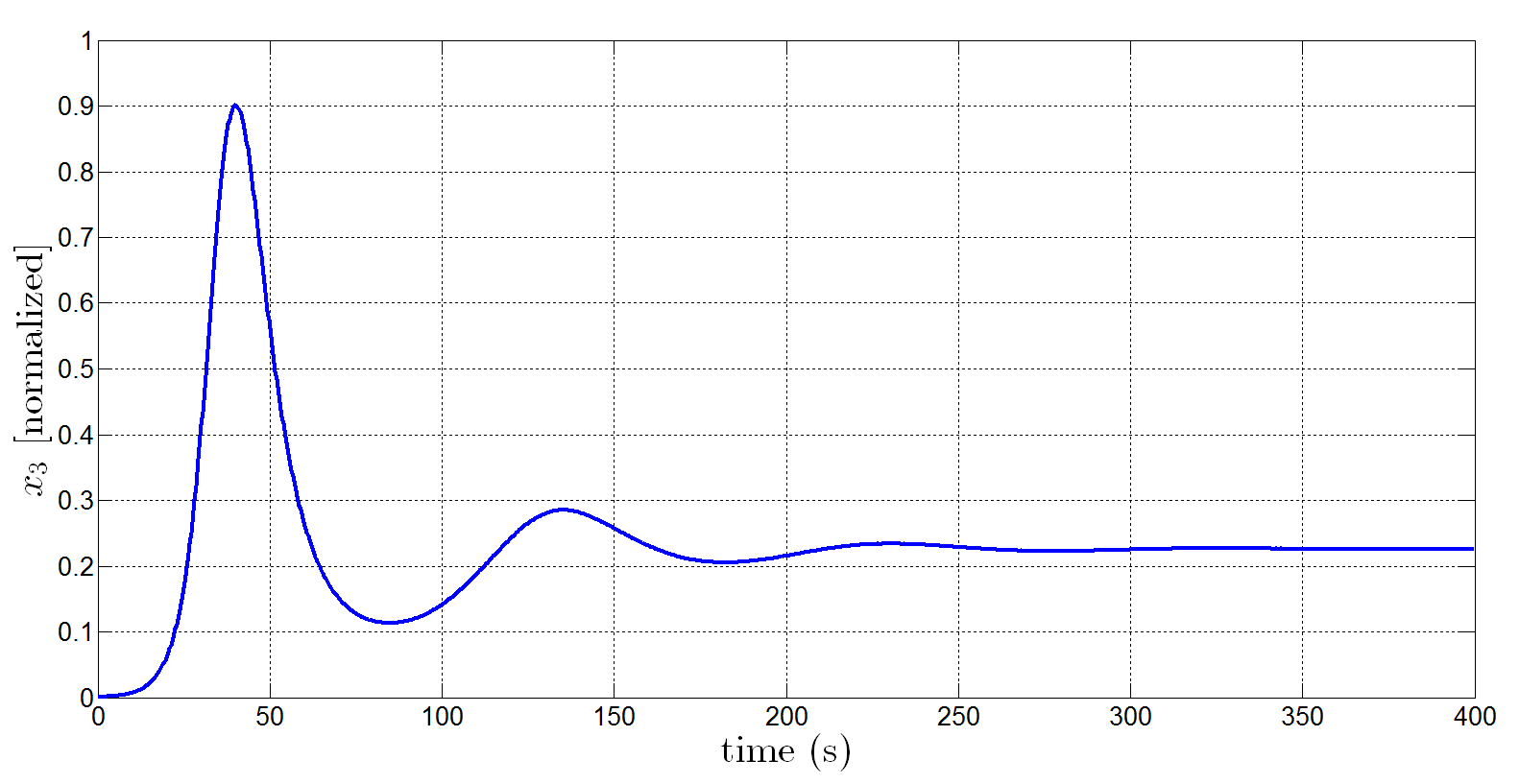}
\caption{Transient response of $y$ considering $u_1 = u_2 = 0$.}
\label{fig_HIV_BO}
\end{figure}

\begin{figure}[!h]
  \begin{center}
    \subfigure[Output $y$]{ \includegraphics[width=16cm]{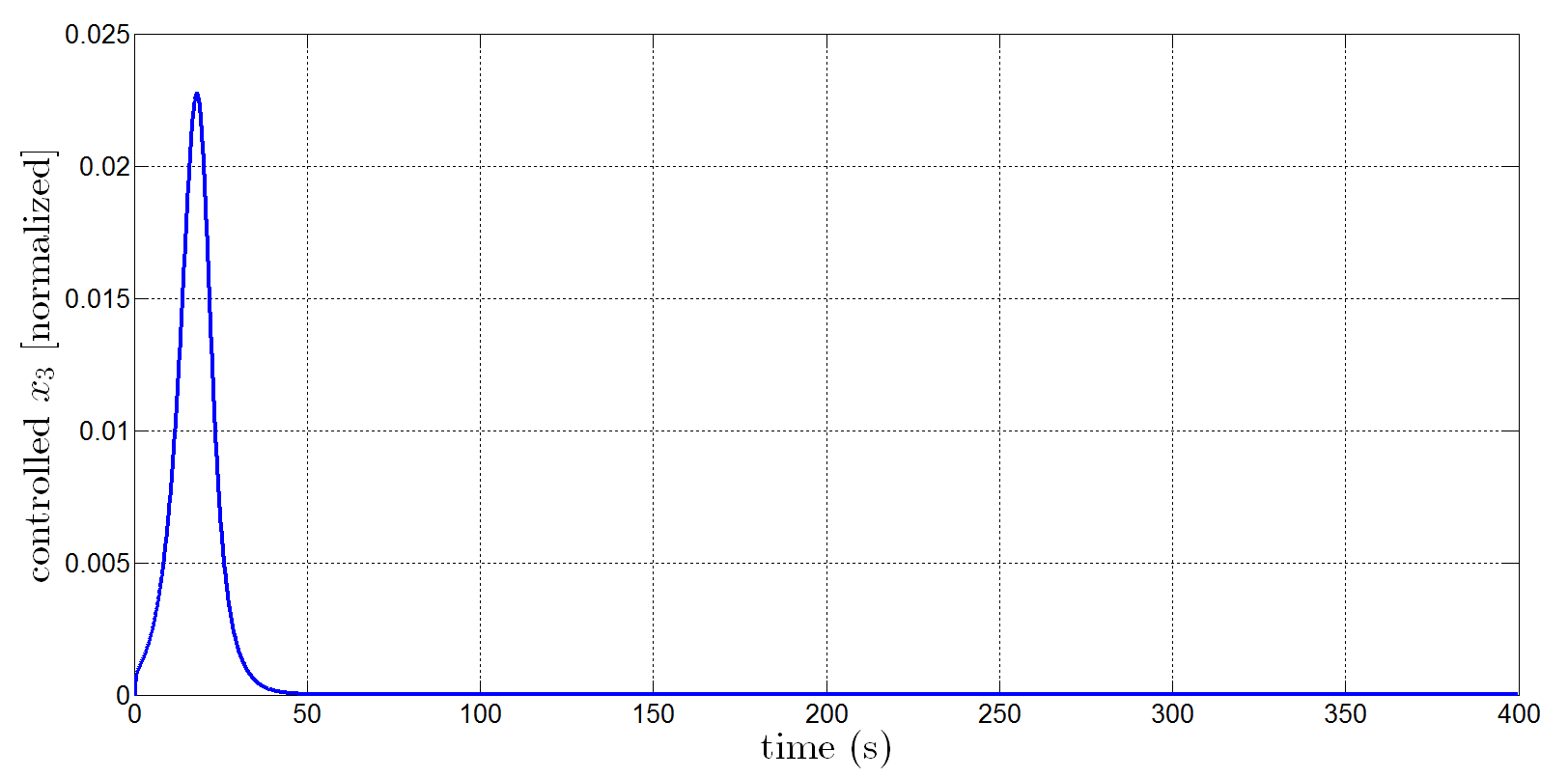}  \label{fig_HIV_BF_11}} 
    \subfigure[Input $u$]{\includegraphics[width=16cm]{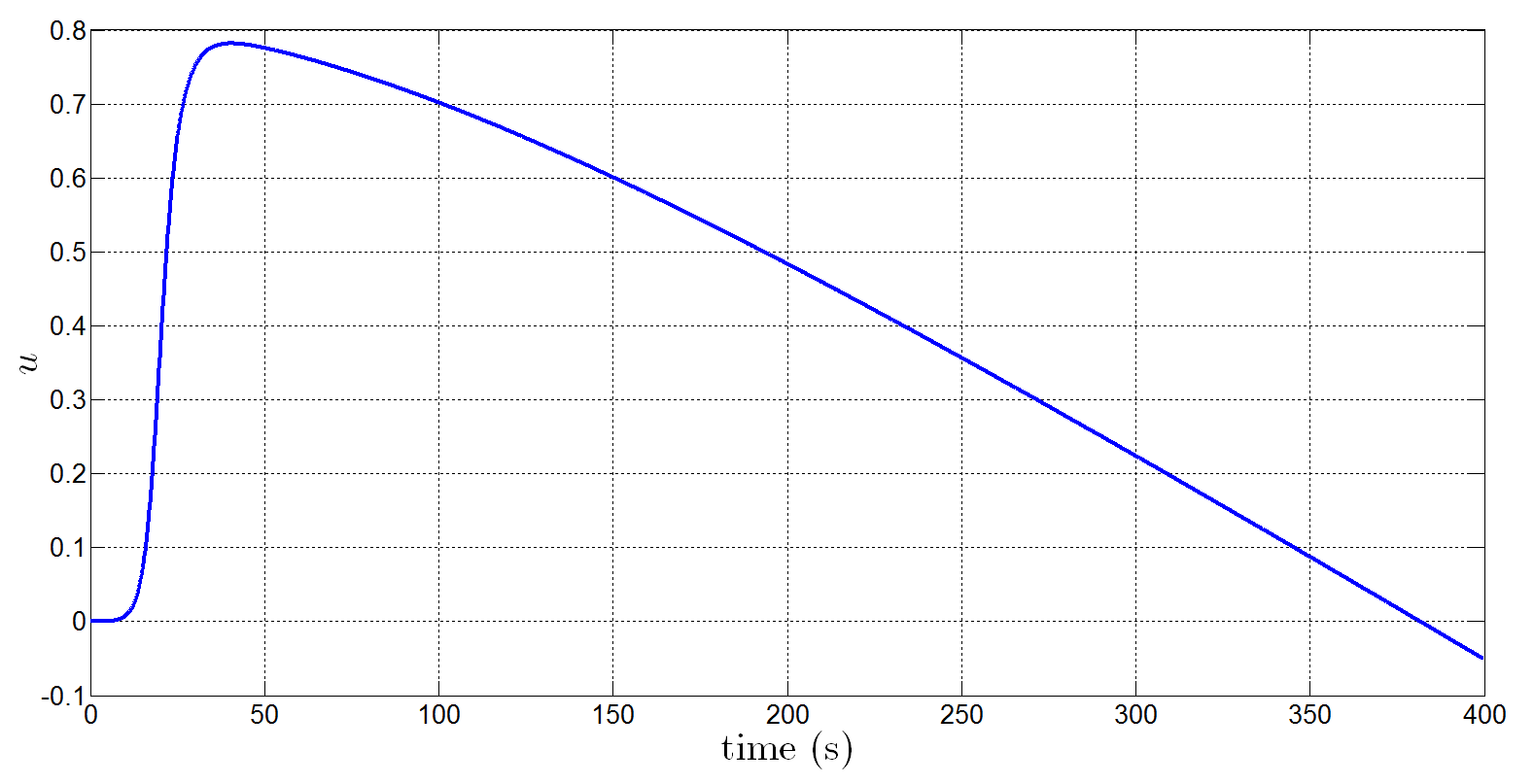} \label{fig_HIV_BF_12}}
    \caption{Controlled output $y$ (viral load) in correspondence with $u_1 = u_2 = u$.}
    \label{fig_HIV_BF_1}
  \end{center}
\end{figure}

\begin{figure}[!h]
  \begin{center}
    \subfigure[Output $y$]{ \includegraphics[width=16cm]{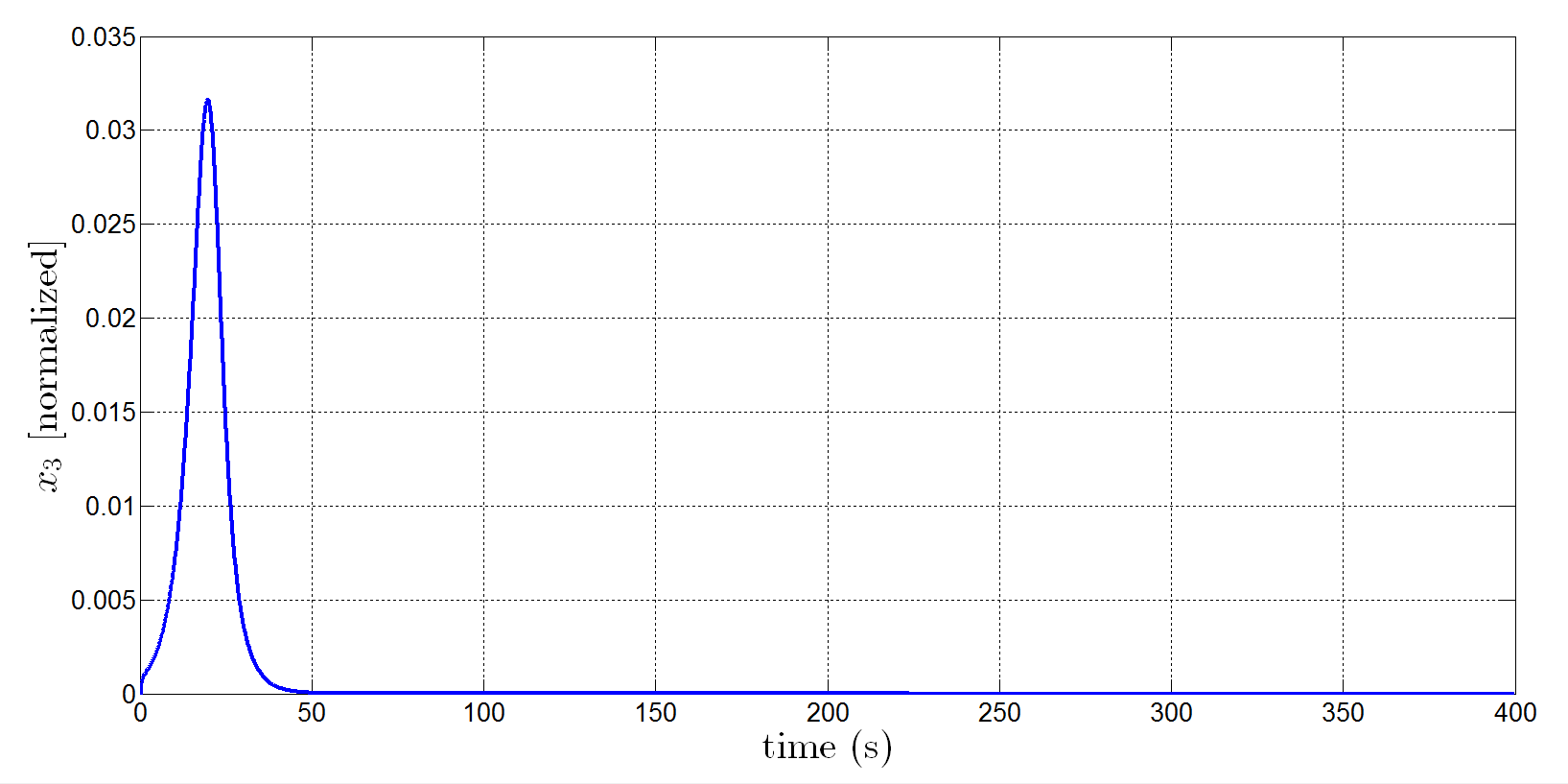}  \label{fig_HIV_BF_21}} 
    \subfigure[Input $u$]{\includegraphics[width=16cm]{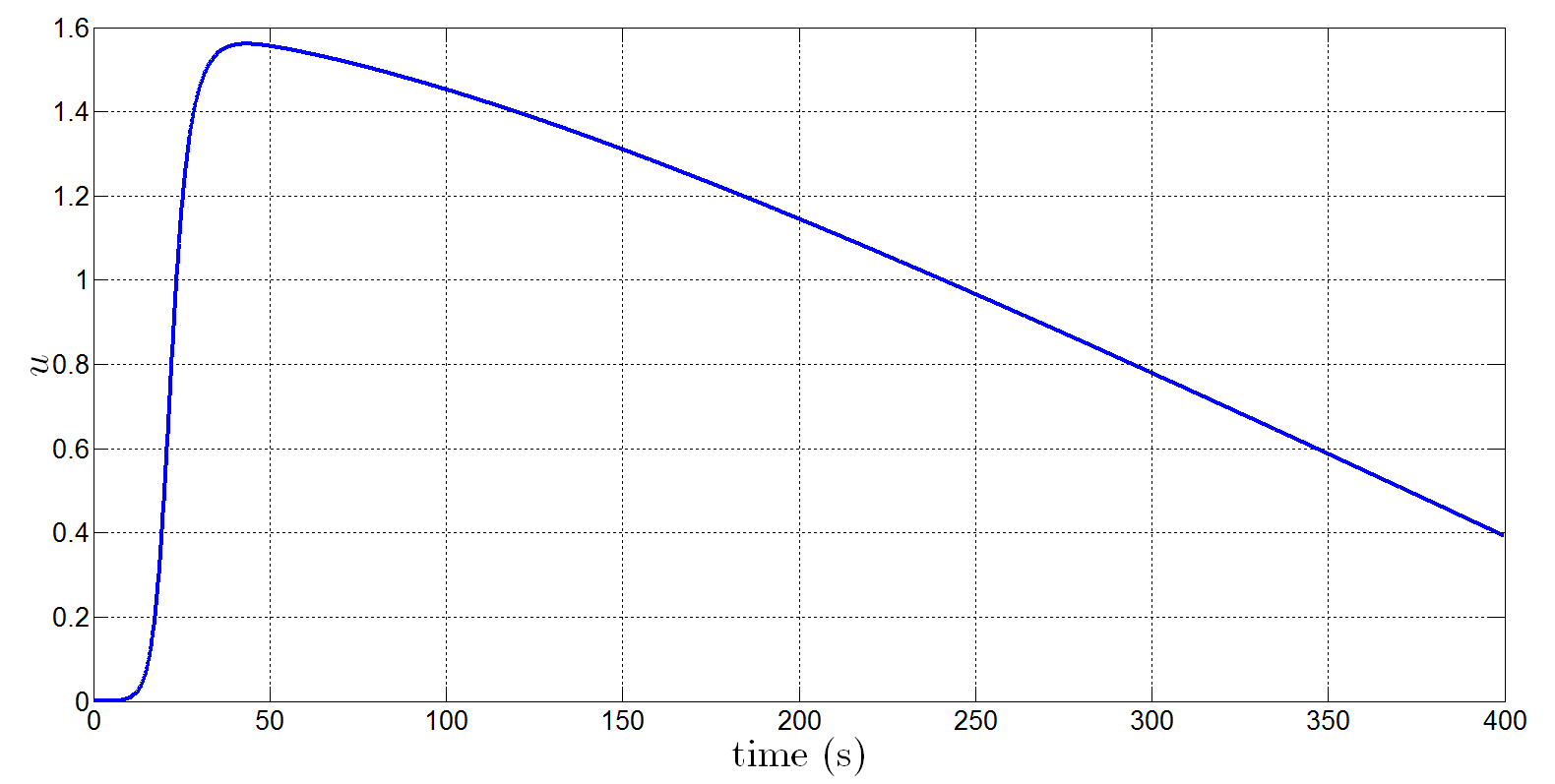} \label{fig_HIV_BF_22}}
    \caption{Controlled output $y$ (viral load) in correspondence with $u_1 = u_2 = \frac{1}{2}u$.}
    \label{fig_HIV_BF_2}
  \end{center}
\end{figure}

The problem is to control the predator-prey like model that describes the evolution of the HIV-1 when subjected to an external "medical agent". From a mathematical point of view, 
we study the possibility of controlling the model (\ref{eq:HIV_form}) for which the purpose is to control the output $y$ (corresponding to the viral load) using the double 
inputs $u_1$ and $u_2$ in such manner that $y$ converges rapidly to zero \footnote{Since we are trying to control this model only from the mathematical point of view, we do not take into account 
the constraints that are medically imposed.} \cite{Barao} \cite{Craig}.

\begin{equation}
\left\{ \begin{array}{l}
\dot{x}_1 =  s - d x_1 - (1 - u_1)  \beta x_1  x_3 \\
\dot{x}_2 = (1 -u_1)  \beta  x_1  x_3 - \mu  x_2 \\
\dot{x}_3 = (1 - u_2)  k  x_2 - c  x_3 \\
y = \begin{pmatrix}
    0 & 0 & \gamma
    \end{pmatrix} x
\end{array}
\right.
\label{eq:HIV_form}
\end{equation}

\noindent
where (mathematically) :  $d = 0.02, k = 100, s = 10, \beta = 2.4.10^{-5}, \mu = 0.24, c = 2.4$. $\gamma$ is a scaling factor that allows normalizing the output.
Figure \ref{fig_HIV_BO} presents the evolution of the output $y$ in open-loop when $u_1 = u_2 = 0$ i.e. when no medical drug is considered.
Figures \ref{fig_HIV_BF_1} and \ref{fig_HIV_BF_2} illustrate the control of $y$ considering two different ratios between $u_1$ and $u_2$.

\clearpage

\clearpage
\subsection{Control of the Epstein frame}

\begin{figure}[!h]
\centering
\includegraphics[width=10cm]{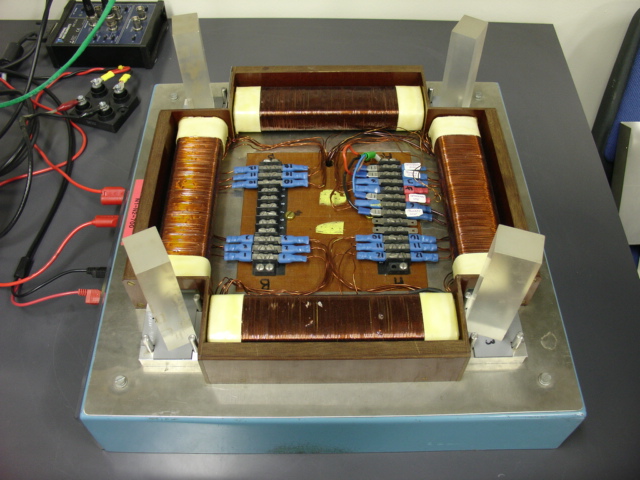}
\caption{An experimental Epstein frame to characterize magnetic materials.}
\label{fig:fig_0__}
\end{figure}

The Epstein frame (see Fig. \ref{fig:fig_0__} \footnote{Picture taken from Wikipédia \url{http://upload.wikimedia.org/wikipedia/commons/5/51/Epstein_frame.jpg}.}) aims to characterize a magnetic material 
by determining its $B-H$ hysteresis curve. The principle is to create a {\it magnetic field} $H$ inside 
the material using a "magnetizing" current $i_H$. The material gives a response to the field $H$ that physically corresponds to the 
measurable {\it magnetic induction field} $B$. This $B$ field creates a voltage $v_B$ through magnetic induction and the quantities $v_B$ and $i_H$ is a representation of the magnetic hysteresis 
curve $B-H$. To describe experimentally the major $B-H$ hysteresis loop, the material under study has to be magnetized using a current $i_H$ that is alternative and of enough magnitude 
in order to describe the magnetic behavior at saturation. {\it The purpose of the control law implementation is to control $i_H$ such that the output voltage of the Epstein frame $v_B$ remains "as close 
as possible" to a desired reference waveform.}

\subsubsection{Epstein frame control}

\paragraph{Proposed $\mathcal{C}_{\pi}$-control scheme} 

\begin{figure}[!h]
\centering
\includegraphics[width=15cm]{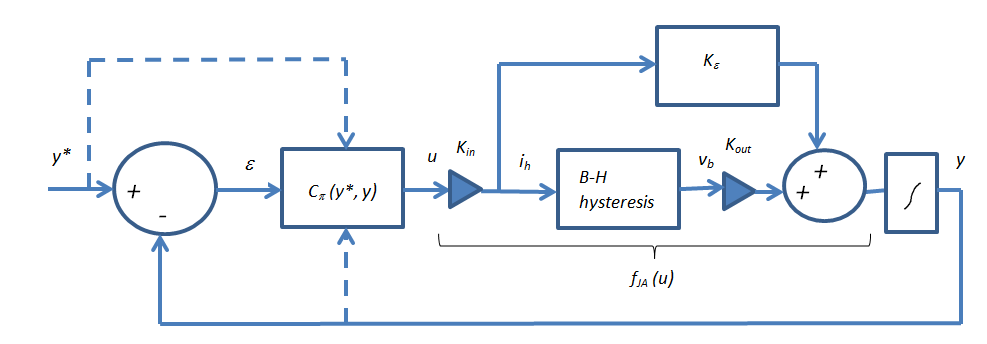}
\caption{Proposed PMA scheme to control the electrical waveforms measured from a magnetic hysteresis.}
\label{fig:fig_1__}
\end{figure}

Consider the control scheme depicted in Fig. \ref{fig:fig_1__} where $\mathcal{C}_{\pi}$ is the proposed PMA controller. $K_{in}$ and $K_{out}$ are positive real gains. Denote $f_{BH}$
the numerical Jiles-Atherton model that is associated to the magnetic hysteresis $B-H$ and $f_{JA}$ is the complete hysteresis to control. 

\paragraph{Jiles-Atherton based hysteresis model}

The Jiles-Atherton model \cite{Jiles} describes a magnetic hysteresis cycle $B - H$. It reads:

\begin{equation}\label{eq:JA_model}
 \frac{d \, M}{d \, H} = \frac{1}{1+c} \frac{M_{an} - M}{\delta k - \alpha (M_{an} - M)} + \frac{c}{1+c} \frac{d \, M_{an}}{d \, H}
\end{equation}

\noindent
where $c$, $\delta k$, $M_{an}$, $\alpha$ are physical coefficients well identified from magnetic hysteresis measurement and we assume that the current $i_H$ corresponds to the magnetization $H$ i.e.
$i_H \propto H$ and the voltage $v_B$ corresponds to the derivative of the magnetic induction field response $B = \mu_0 H + J_{BH}(H)$ (where $J_{BH}$ describes the $B-H$ hysteresis via (\ref{eq:JA_model})). 

\paragraph{Simplified model of the Epstein frame} The Epstein frame admits a complex model based on the Jiles-Atherton model that represents all electric phenomena that occur inside the Epstein frame\footnote{The 
Epstein frame is equivalent to a transformer and the "mutual interactions" between primary and secondary coils must be taken into account in addition to the hysteresis behavior.}. 
To simplify the model of the Epstein frame to control, we consider controlling a nonlinear function $f_{JA}$, which represents a modified Jiles-Atherton model. Denote $v_H = f_{JA}(i_H)$ the nonlinear dynamical 
system that describes the $B-H$ hysteresis as a function of $i_H$,
and consider controlling directly the hysteresis cycle in such manner that $\mathcal{C}_{\pi}$ controls $i_H$ in order to get $v_H = f_{JA}(i_H)$ as close as possible to a reference waveform.

To define the global $f_{JA}$ hysteresis function that is controlled by $\mathcal{C}_{\pi}$, which includes the scaling coefficients needed by the $\mathcal{C}_{\pi}$ corrector, 
consider $i_H = K_{in} \, u $, $y = K_{out} v_H$ and a coefficient $K_e$ such that:

\begin{equation}\label{eq:fJA}
y = f_{JA} (u) = K_{in} K_e u + K_{out} J_{BH} (K_{in} u).
\end{equation}

\noindent
The small coefficient $K_e \, \# \, 10^{-4}$ compensates the very small variations of $f_{BH}$ when $B$ is close to $B_{sat}$. Such variations may induce time-delays in the response of the 
$\mathcal{C}_{\pi}$-controller that induce some distortions of the output signal $y$. The model (\ref{eq:fJA}) could be seen as an "affine" derivation of the original Jiles-Atherton model.

\begin{figure}[!h].
\centering
\includegraphics[width=16cm]{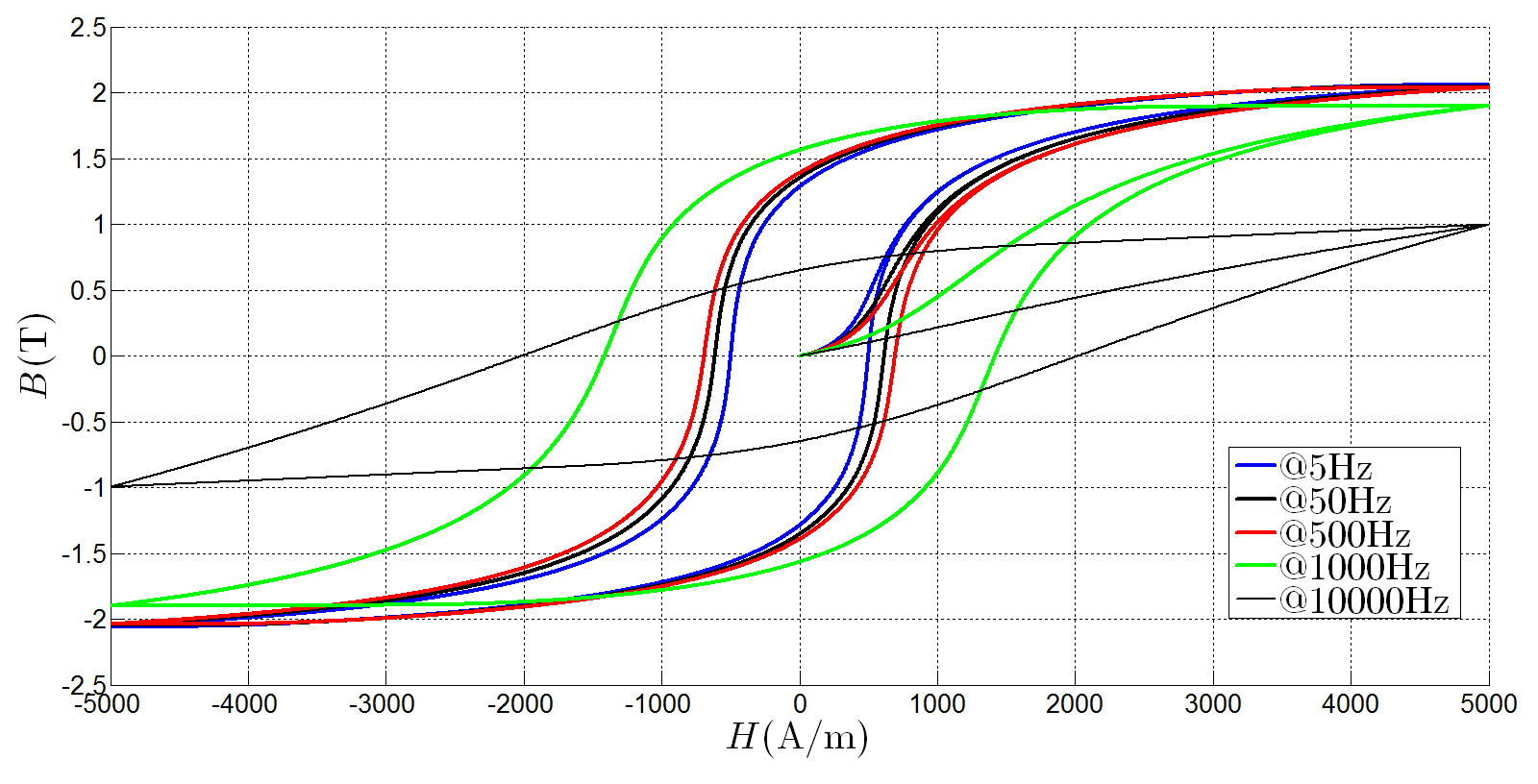}
\caption{Simulation of the Jiles-Atherton model for different operating frequencies.}
\label{fig:fig_2__}
\end{figure}

The  $B-H$ loops, obtained from the Jiles-Atherton model, are depicted in Fig. \ref{fig:fig_2__} considering the frequencies 5 Hz, 50 Hz, 500 Hz, 1 kHz and 10 kHz. Since this hysteresis model does not allow 
$H > H_{max}$, a limitation is necessary to bound the evolution of $H$, that may occur eventually during the transient of the dynamic stabilization of the control loop (ex. in Fig. \ref{fig:fig_6__}).

\subsubsection{Simulation results}

\begin{figure}[!h]
\centering
\includegraphics[width=15cm]{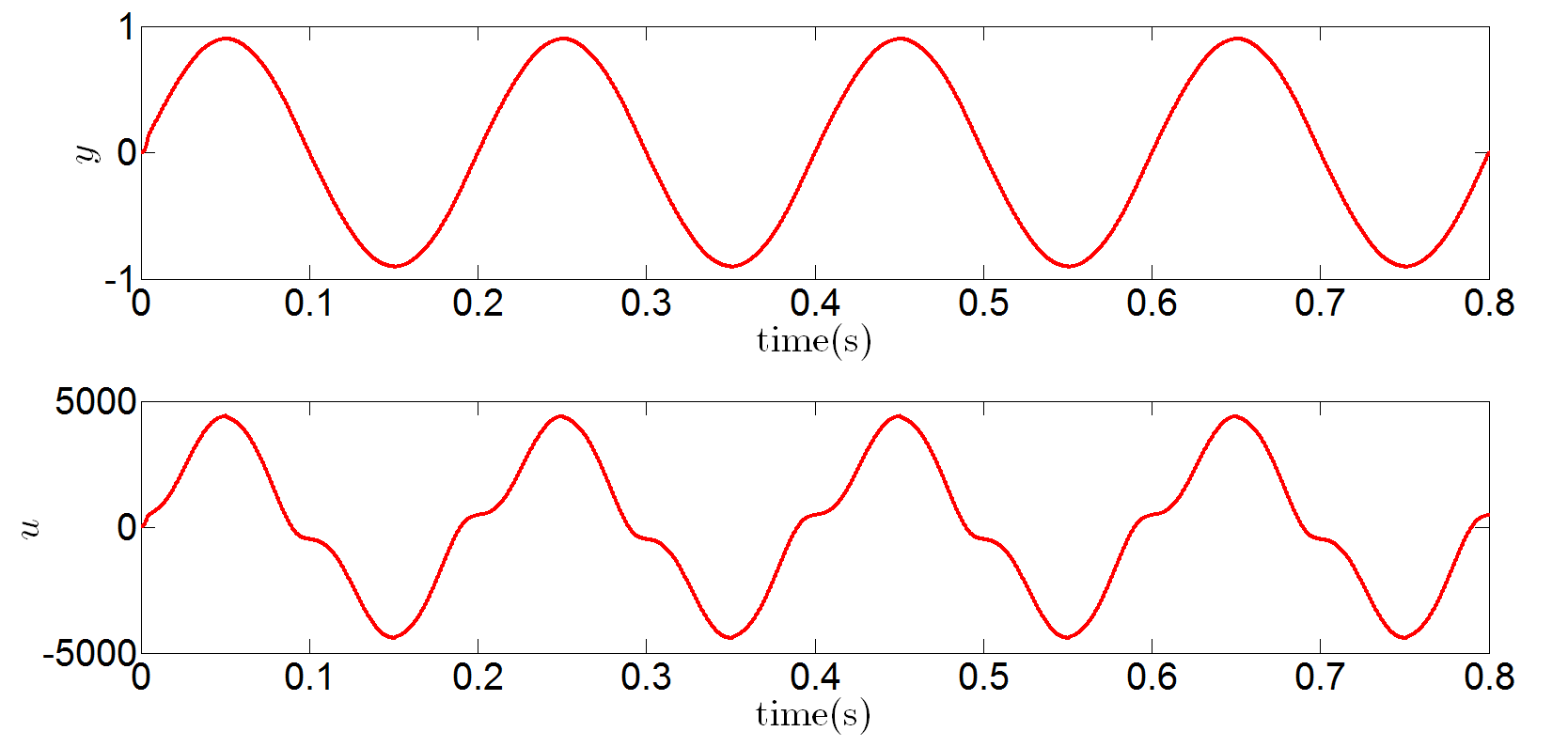}
\caption{Simulated $u$ and $y$ signals according to the time at 5 Hz.}
\label{fig:fig_3__}
\end{figure}

\begin{figure}[!b]
\centering
\includegraphics[width=15cm]{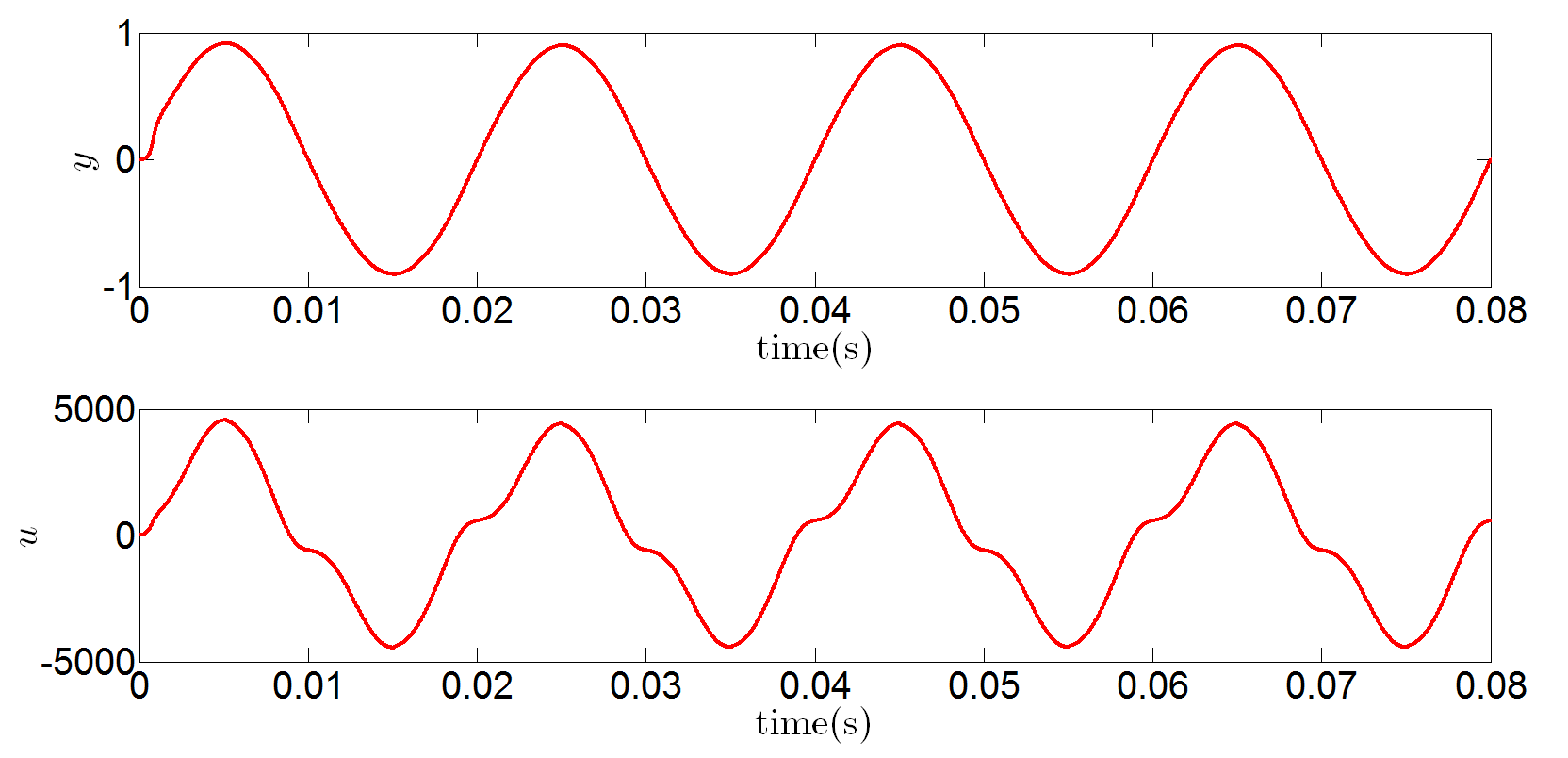}
\caption{Simulated $u$ and $y$ signals according to the time at 50 Hz.}
\label{fig:fig_4__}
\vspace{2cm}
\end{figure}

\begin{figure}[!h]
\centering
\includegraphics[width=15cm]{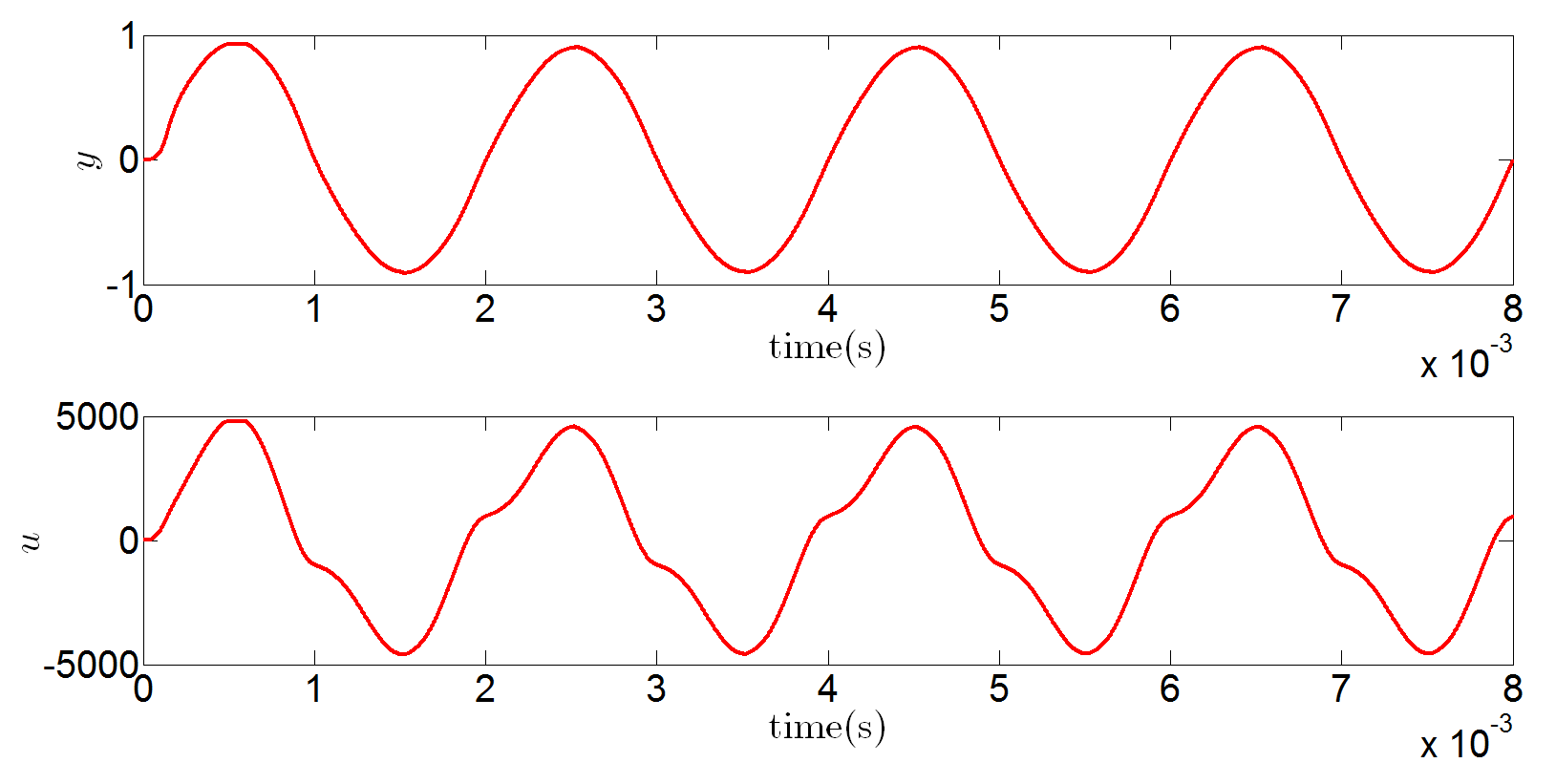}
\caption{Simulated $u$ and $y$ signals according to the time at 500 Hz.}
\label{fig:fig_5__}
\end{figure}

\begin{figure}[!b]
\centering
\includegraphics[width=15cm]{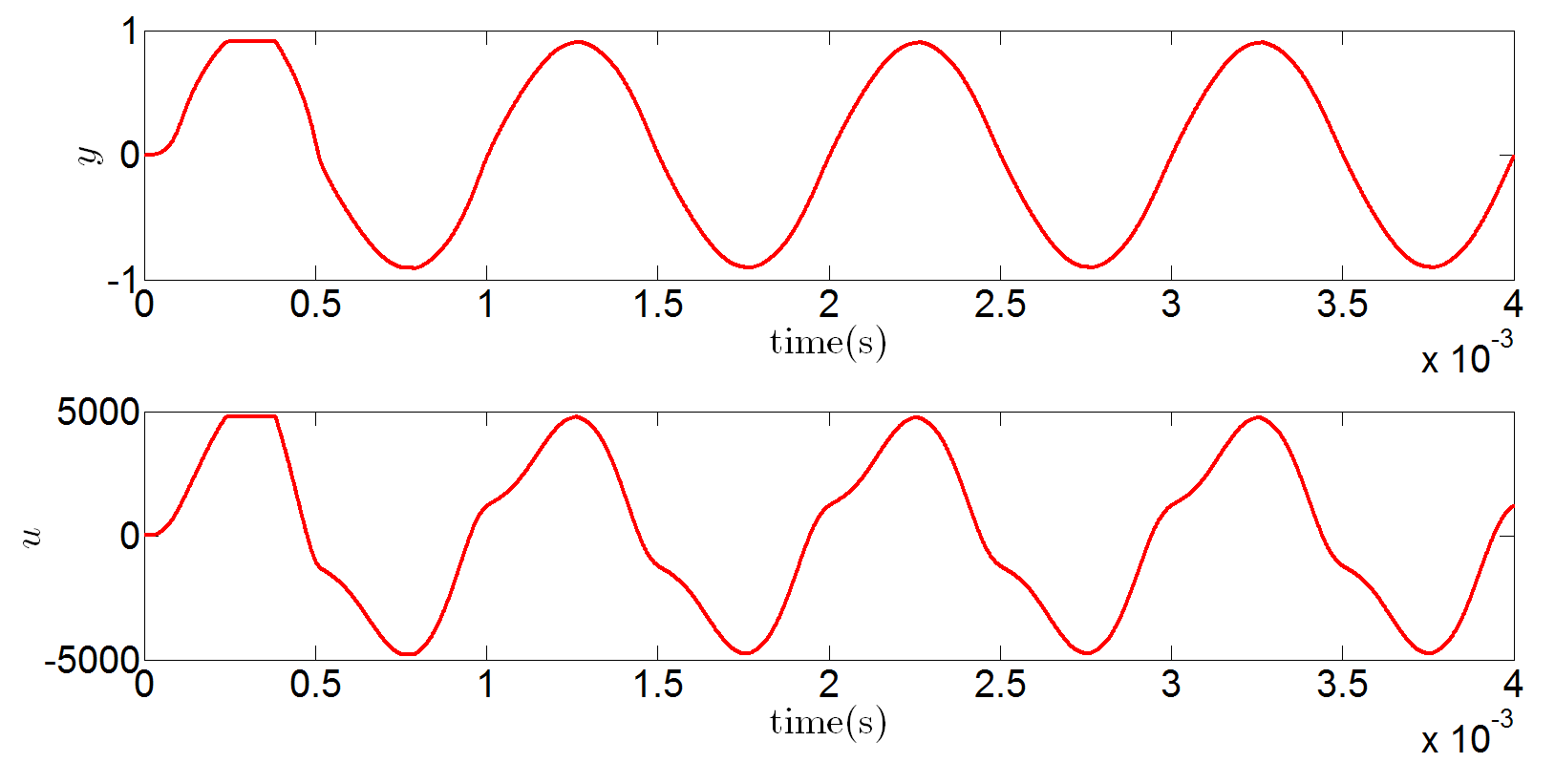}
\caption{Simulated $u$ and $y$ signals according to the time at 1 kHz.}
\label{fig:fig_6__}
\vspace{2cm}
\end{figure}

\begin{figure}[!h]
\centering
\includegraphics[width=15cm]{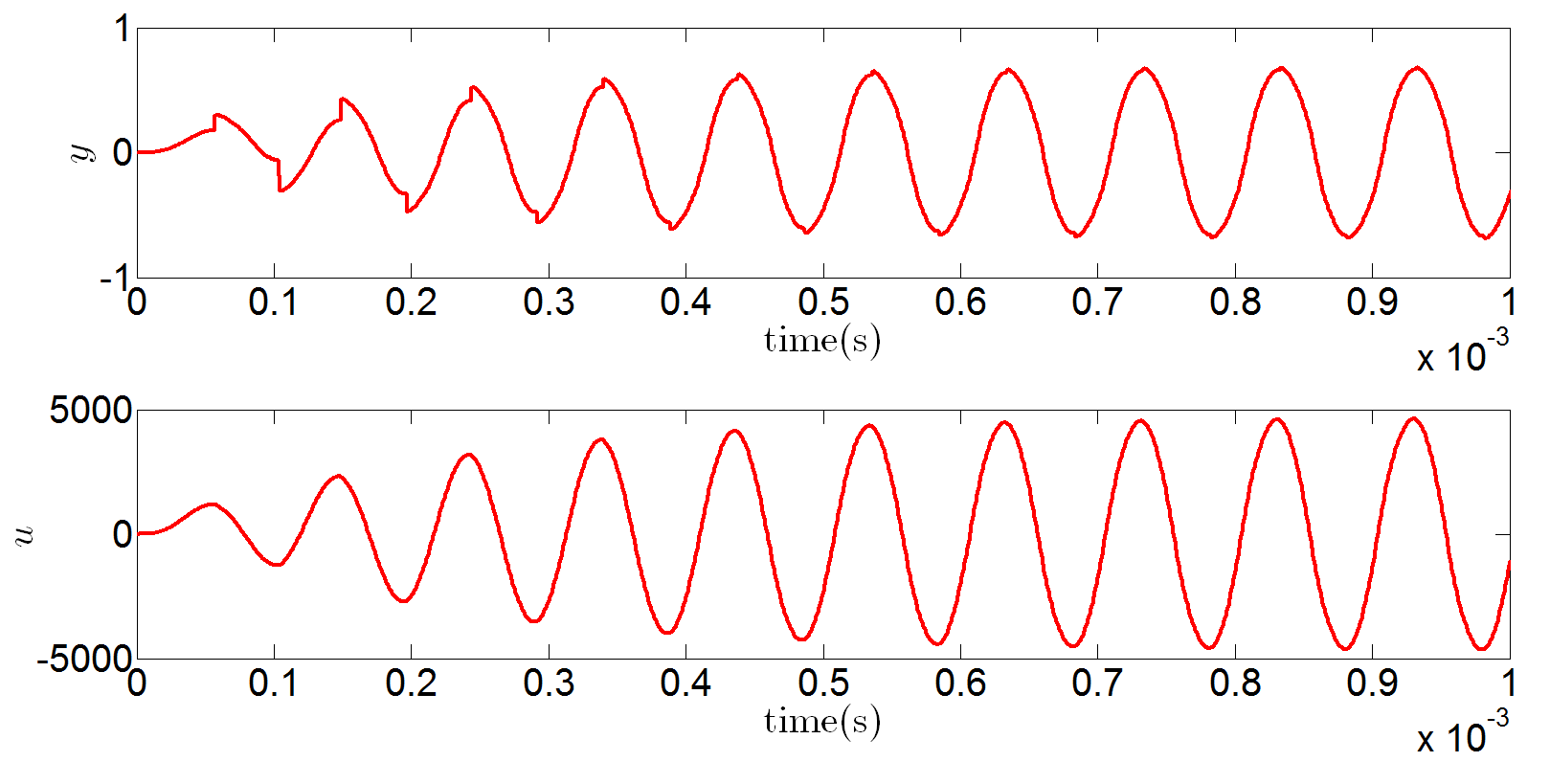}
\caption{Simulated $u$ and $y$ signals according to the time at 10 kHz.}
\label{fig:fig_7__}
\end{figure}

Figures \ref{fig:fig_3__}, \ref{fig:fig_4__}, \ref{fig:fig_5__}, \ref{fig:fig_6__} and \ref{fig:fig_7__} depict the input $u$ and the (rescaled) output $y$ of the control loop according to the time. 
Given a particular operating frequency, for which a particular $H-B$ hysteresis is studied (see Fig. \ref{fig:fig_2__}), and assuming that the output reference $y^*$ is a 
sinusoid, whose magnitude corresponds to the theoretical $H_{max}$ of the $H-B$ hysteresis, different frequencies are considered (5 Hz, 50 Hz, 500 Hz, $\,$ 1 kHz and 10 kHz) in order to highlight the 
behavior of the controlled voltage $v_b$ when the frequency changes. In particular, high frequencies (e.g. Figures \ref{fig:fig_6__} and \ref{fig:fig_7__}) introduce an important transient response on $y$ due to the fact that the variations of $y^*$
are too fast to get an immediate stabilization to the dynamic working point of the hysteresis. The simulation show that the $\mathcal{C}_{\pi}$-controller gives very interesting dynamic performances over a wide range of dynamic
working points relating to the operating frequency. An optimization
algorithm (see \S \ref{BFO_paragraph}) has been used to adjust the parameters of the $\mathcal{C}_{\pi}$-controller in such manner that the shape of the output response $y$ is "as close to" a 
sine shape\footnote{Remember that the purpose of the optimization procedure is to minimize the tracking error $y^* - y$ in such manner that ideally $y \equiv y^*$ for the particular sine output reference $y^*$.}.

\begin{figure}[!b]
\centering
\includegraphics[width=15cm]{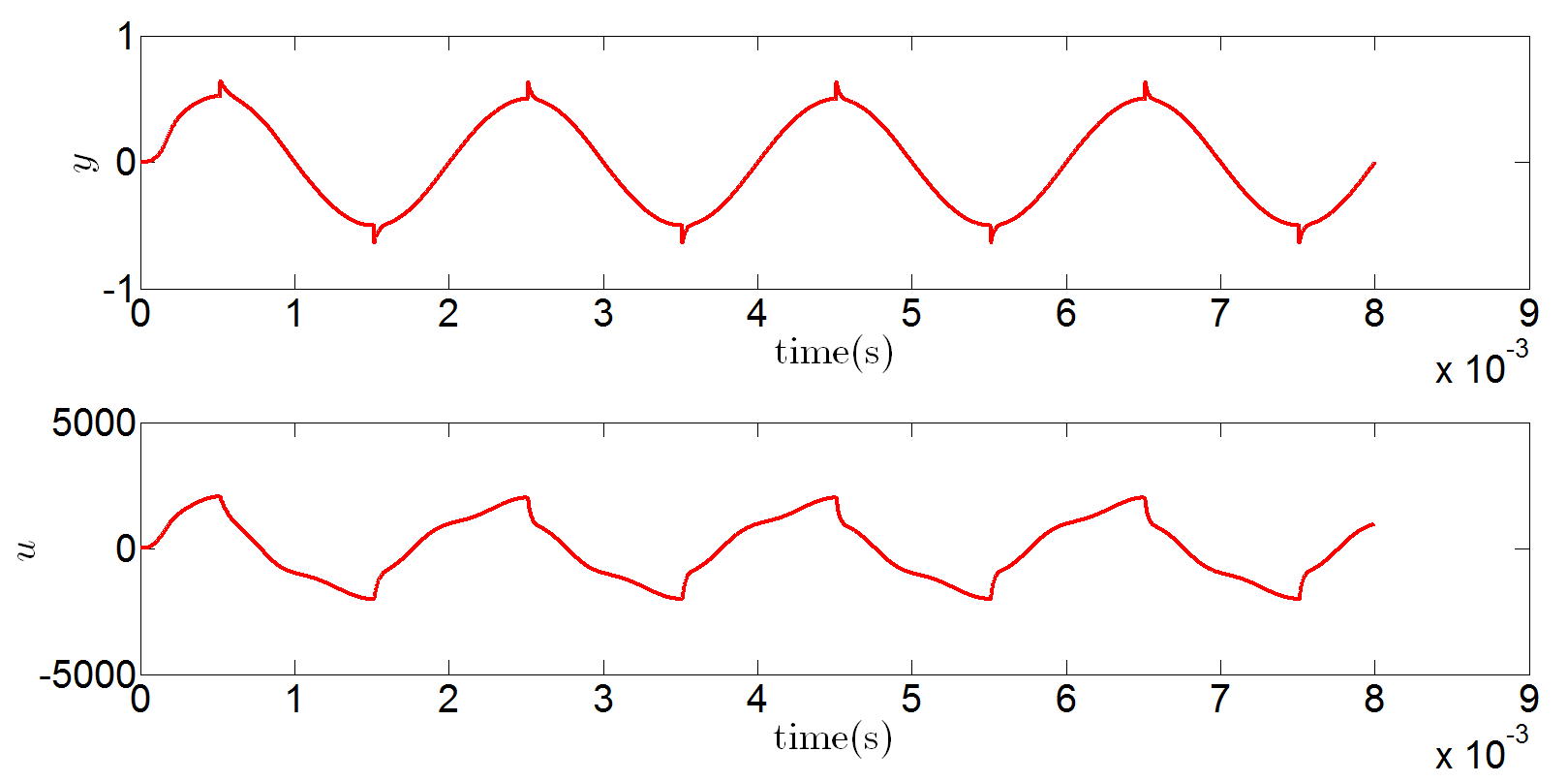}
\caption{Simulated $u$ and $y$ signals according to the time at 500 Hz $(H << H_{max})$.}
\label{fig:fig_8__}
\vspace{2cm}
\end{figure}

\paragraph{Remarks} This hysteresis model is composed of three subsystems (the first magnetization branch, the increasing and decreasing branches) that switch depending on the value of $d H / d t$. When $H << H_{max}$, the switch between
the branches may not be smooth and such "connection" may induce a small transient on $y$. An illustration is presented in Fig. \ref{fig:fig_8__} at a low frequency in 
comparison with Fig. \ref{fig:fig_5__}.

\begin{figure}[!h]
\centering
\includegraphics[width=15cm]{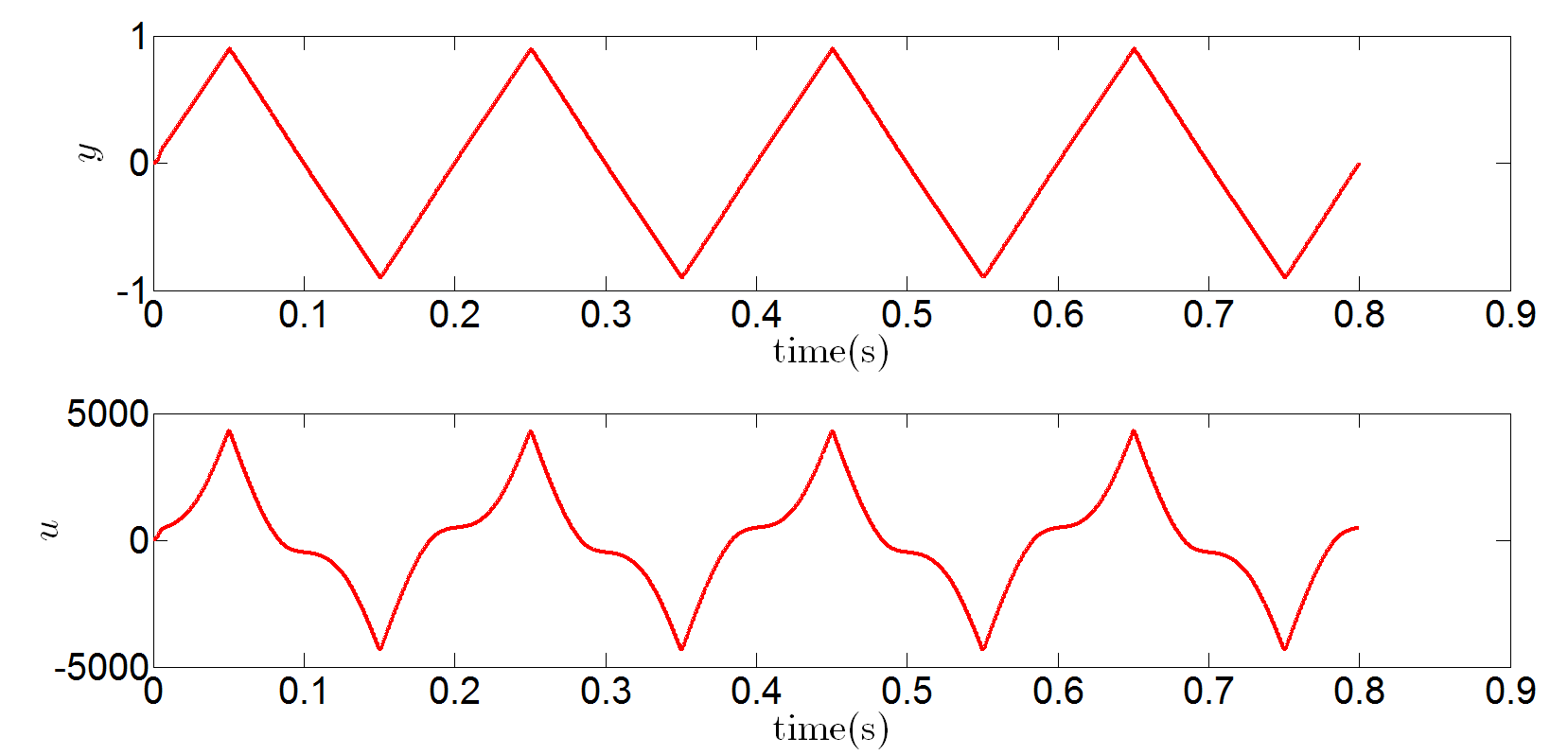}
\caption{Simulated $H$ and $y$ signals according to the time at 5 Hz.}
\label{fig:fig_9__}
\end{figure}

The closed-loop has been also tested using a triangular shaped output reference $y^*$. Figure \ref{fig:fig_9__} depicts the magnitude of the magnetic field $H$ and the corresponding 
(rescaled) magnetic induction $B$ during the control loop process according to the time. The frequency of 5 Hz has been considered as an example holding the parameters 
of the simulation with the sine reference at 5 Hz.

\clearpage

\section{Derivative-free $\&$ "extremum-seeking" control}

To describe how the PMA could be used as a derivative-free optimization (DFO) algorithm (e.g. \cite{Vicente} \cite{global}) or as an "extremum-seeking" (ES) control scheme (e.g. \cite{Ari} \cite{seek} \cite{one}), 
we first define each element of the associated control 
scheme and then, we derive the operating conditions that would allow to minimize nonlinear functions. {\it We assume that it is possible to derive a control scheme such that the PMA 
can be used to minimize nonlinear functions.}

\subsection{Proposed $\mathcal{C}_{\pi}$-control scheme} 

\paragraph{Definition of the closed loop}

\begin{figure}[!h]
\centering
\includegraphics[width=15cm]{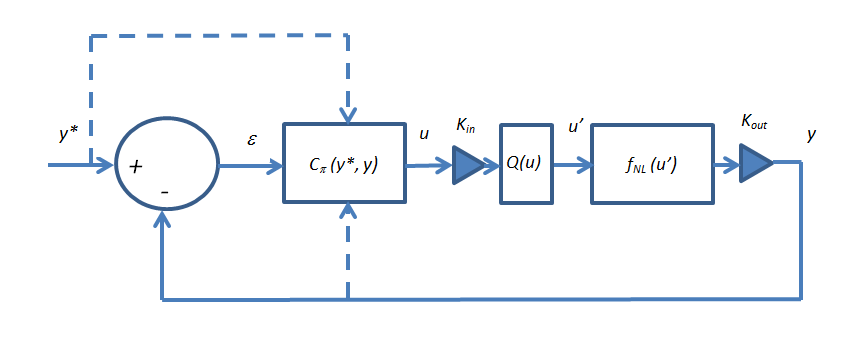}
\caption{Proposed PMA scheme to minimize a nonlinear function $f_{nl}$.}
\label{fig:CSM_DFO}
\end{figure}

Consider the control scheme depicted in Fig. \ref{fig:CSM_DFO} where $\mathcal{C}_{\pi}$ is the proposed PMA "extremum-seeking" controller. $K_{in}$ and $K_{out}$ are positive real gains. 
We consider either a static nonlinear function $f_{nl}^s$ (regarding DFO), which does not have any internal dynamical properties, or a nonlinear SISO dynamical system (\ref{eq:gen_sys}) 
(regarding ES), to minimize. 
The function $Q$ is e.g. a basic first order transfer function.

\paragraph{Function to control} 

\begin{itemize}

\item Define the $f_{nl}^s$ static function to optimize such that: 

\begin{equation}
f_{nl}^s : \begin{array}{c}
          \mathbb{R}^n \rightarrow \mathbb{R} \\
          u' \mapsto y
         \end{array} 
\end{equation}

or consider a nonlinear SISO dynamical system (\ref{eq:gen_sys}).

This function represents the "nonlinear optimization problem". Currently, the assumption $n \leq 2$ is considered and we denote $u'_x$ the input variable 
for $n = 1$.

\item $Q$ is a standard linear transfer function (typically of first order), such that:

\begin{equation}
Q : \begin{array}{c}
     \mathbb{R} \rightarrow \mathbb{R} \\
   \displaystyle{   \left. \frac{d u'}{d t} \right|_k + \gamma  \left. u'\right|_k = u_k }
    \end{array} 
\end{equation}

\noindent
where $\gamma$ is a time-constant, chosen in such manner that the step response of $Q$ is very fast.
As presented in the Fig. \ref{fig:CSM_DFO}, $Q$ is associated with $f_{nl}^s$ in order to provide some minimal dynamical properties regarding the system $f_{nl}^s$ to control.
Obviously, the $Q$ function is not necessary when $f_{nl}$ is already a dynamical function like (\ref{eq:gen_sys})\footnote{Last investigations suggest that the linear transfer function $Q$
may be not necessary even for a $f_{nl}^s$ function to control. The properties of the para-model algorithm are currently under study considering  nonlinear systems that are "non-dynamical".}.

\end{itemize}

A single $\mathcal{C}_{\pi}$-controller\footnote{To extend this scheme to 
multi-input variables $(n > 2)$ of $f_{nl}$, one may consider the use of a $\mathcal{C}_{\pi}$-controller per input variable.} drives a single input of $f_{nl}$ 
(eventually through $Q$), as presented in Fig. \ref{fig:CSM_DFO}.

\subsection{Numerical applications}
 
Since the PMA is designed for nonlinear systems and does not contain any derivatives, it is assumed that the "extremum-seeking" control is possible considering
a specific definition of $y^*$ in order to reach and stabilize $f_{nl}$ to its minimum.
 
Let us assume that the following (eventually constrained) minimization problem (described for a single variable):
 
 \begin{equation}\label{eq:min_goal}
  \min_{x \in \mathbb{R}} \, f_{nl}(x), \quad ( x = u' \, \hbox{identically inside the control scheme})
 \end{equation}
 
 \noindent
 is equivalent to the control scheme described in Fig. \ref{fig:CSM_DFO}, for which the output reference $y^*$ "follows" the minimum value of $f_{nl}$. 
 We denote $x = x_{opt}$ the value that gives the minimum of $f_{nl}$.

\begin{figure}[!h]
\centering
\includegraphics[width=15cm]{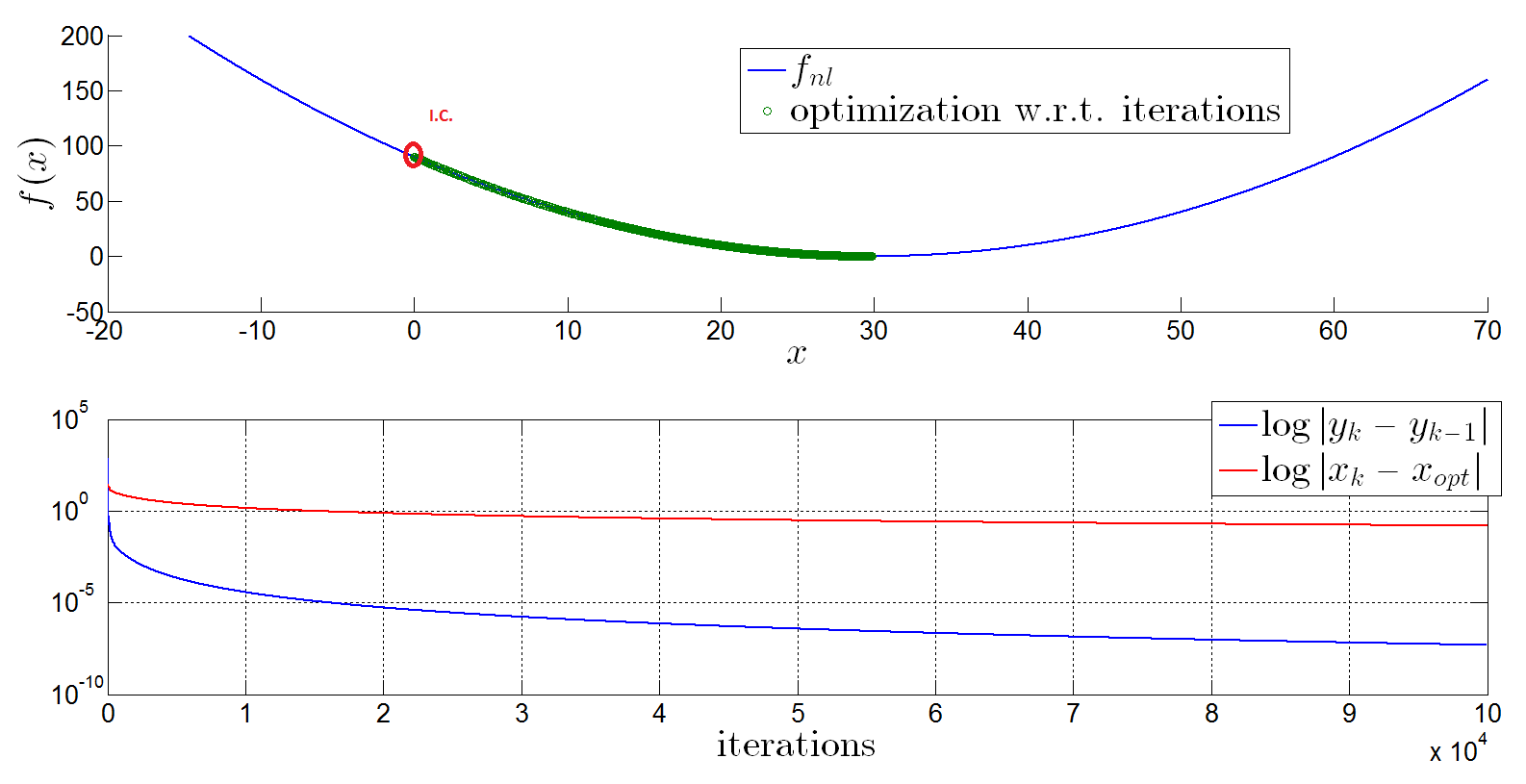}
\caption{$\displaystyle{\min_{x} \, (x - 30)^2}$ (initial condition in red spot).}
\label{fig:fig_1_DFO}
\end{figure}

\begin{figure}[!b]
\centering
\includegraphics[width=15cm]{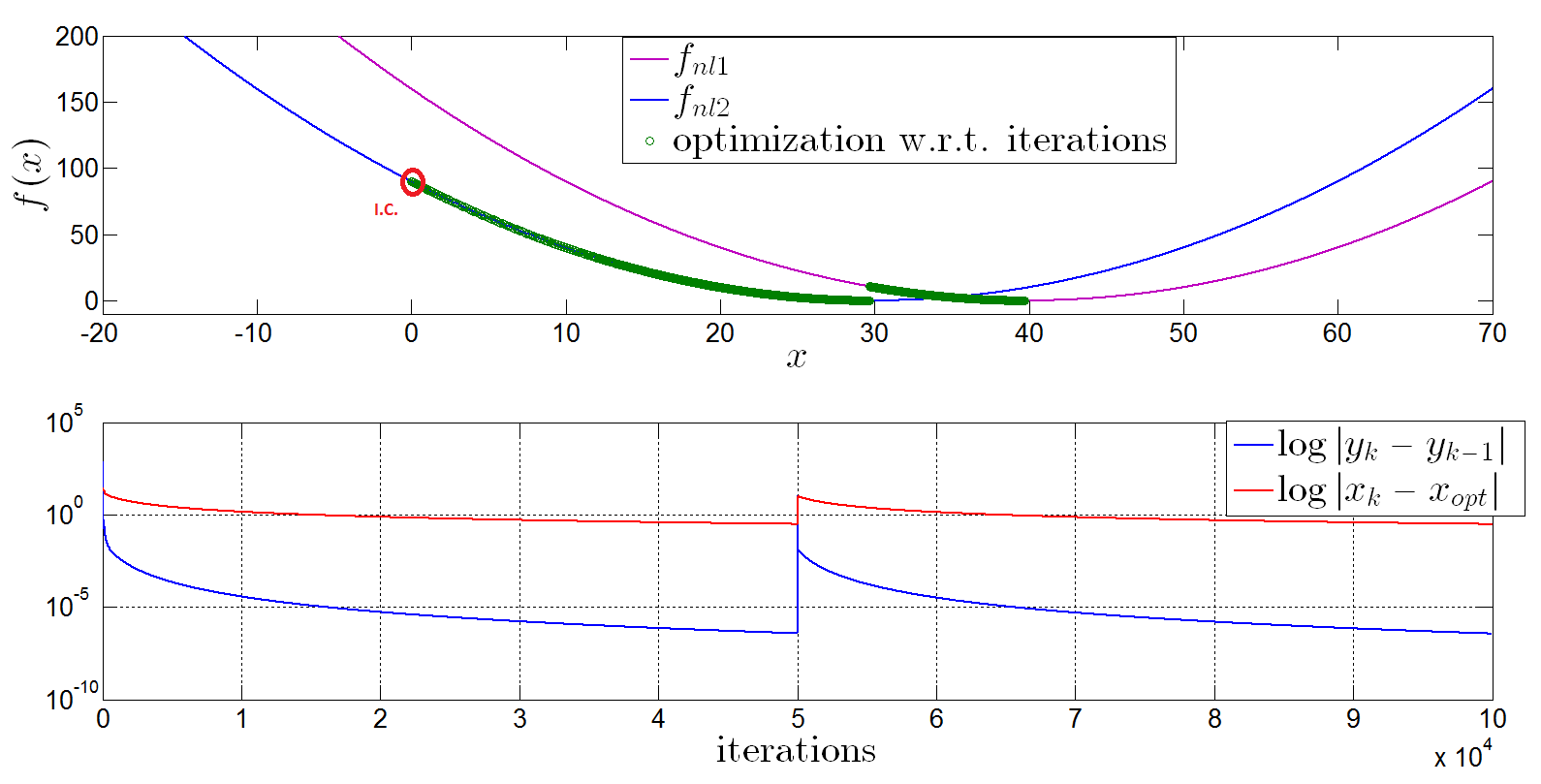}
\caption{$\displaystyle{\min_{x} \, f_{nl2} = (x - 30)^2 \stackrel{t_1 \, ?}{\rightarrow} \min_{x} f_{nl1} =  (x - 40)^2}$ (initial condition in red spot).}
\label{fig:fig_2_DFO}
\end{figure}

\begin{figure}[!h]
\centering
\includegraphics[width=15cm]{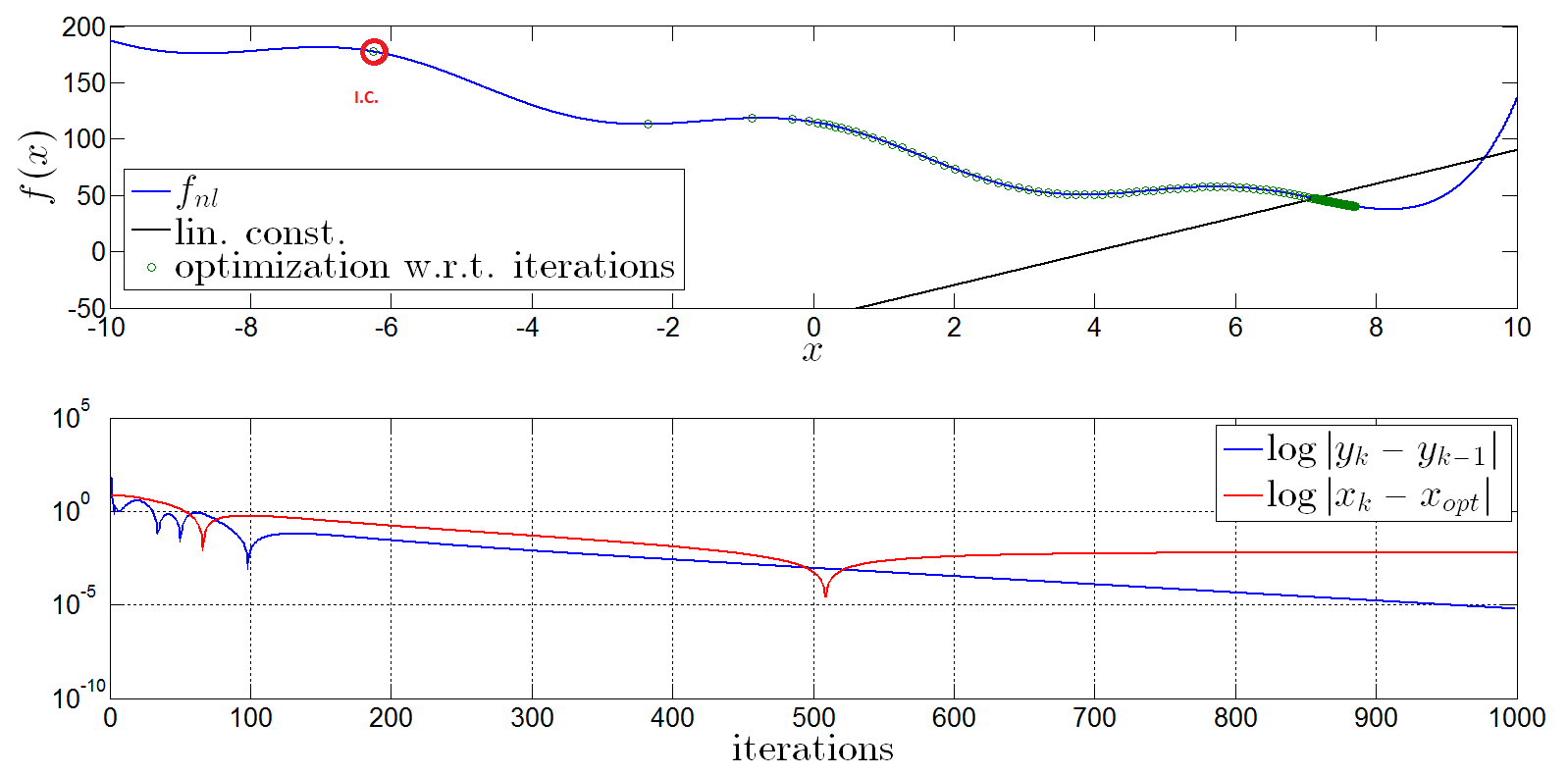}
\caption{$\displaystyle{\min_{x, y} \, 10(1.5 \cos(x) - x) + \exp(x-5)+100 \quad \hbox{subj. to :} \quad y \geq 15x - 60}$ (initial condition in red spot).}
\label{fig:fig_3_DFO}
\end{figure}

\paragraph{Results} 
For each case in 1D, are plotted: the difference between two iterations $y_{k}, \, y_{k-1}$ and
the error between $x_{opt}$ and the evolution of $x$ through the closed-loop.
In these cases, all the parameters of $\mathcal{C}_{\pi}$ have been set experimentally to give interesting performances but are not 
optimal (the choice of the $y^*$ function may influence the speed of the convergence). The following numerical cases are studied:

\begin{itemize}
\item See Fig. \ref{fig:fig_1_DFO} regarding the minimization of a 1D convex function such that: 
\begin{equation}
\min_{x} \, (x - 30)^2
\end{equation}

\item See Fig. \ref{fig:fig_2_DFO} regarding the minimization of a 1D convex function with a minimum that changes according to the time at an unknown instant $t_1$ such that:  
\begin{equation}
\min_{x} \, (x - 30)^2 \stackrel{t_1 \, ?}{\rightarrow} \min_{x} (x - 40)^2
\end{equation}

\item See Fig. \ref{fig:fig_3_DFO} regarding the minimization of a 1D non convex function such that: 

\begin{equation}
\min_{x, y} \, 10(1.5 \cos(x) - x) + \exp(x-5)+100 \quad \hbox{subj. to :} \quad y \geq 15x - 60
\end{equation}

 
 
 
  
\end{itemize}

\section{Concluding remarks}

We presented how the proposed para-model agent\footnote{Why "para-model"? Based on the "model-free" methodology, for which the model is not correlated to the controller in the sense that the controller
does not need an explicit definition of the model to be parametrized and can thus rebuild an ultra-local model from measurements, we propose the prefix "para" to highlight the fact that the 
agent is "in close proximity" 
to the model of the process that he controls. Although no specific information is needed for the PMA control part, some basic information may be needed 
to configure properly the PMA and the output reference $y^*$ in the case of extremum-seeking control approach. The stability study is an essential part that should formalize the proposed PMA approach in order to justify theoretically the 
operating conditions of the $\mathcal{C}_{\pi}$-control.}, as a model-free and derivative-free based controller, can be used to control nonlinear systems or perform optimization / "extremum-seeking" 
control.
Further investigations include extensive tests and applications to complex systems as well as a complete study of the stability.

\section*{Acknowledgement}

The author is sincerely grateful to Dr. Edouard Thomas for his strong guidance and his valuable comments that improved this paper. 
The author is also sincerely grateful to Olivier Ghibaudo, Ph.D. student at G2Elab-CNRS (Grenoble) France, for having worked on the numerical version of the hysteresis model.

\clearpage


\end{document}